\documentclass[preprint,review,11pt,authoryear]{elsarticle}

\usepackage{amsmath,amssymb}
 \usepackage{amsthm}
 \usepackage{bbm}
 \usepackage{todonotes}
 \usepackage{graphicx}
\usepackage{color}
\definecolor{cornell-red}{RGB}{179,27,27}
\usepackage{subcaption}
\usepackage{mathtools}
\usepackage{dsfont}
\usepackage{color}
\usepackage{natbib}
    \usepackage{multirow}
 
\usepackage[left=1in,right=1in,top=1in,bottom=1in]{geometry}

\usepackage{xcolor}
\usepackage[colorlinks = true,linkcolor = blue,urlcolor  = blue,citecolor = blue,anchorcolor = blue]{hyperref}
\usepackage[flushleft]{threeparttable}
\newtheorem{thm}{Theorem}

\newtheorem{prop}[thm]{Proposition}

\theoremstyle{definition}

\theoremstyle{remark}

\newcommand{\E}{\mathbb{E}}
\newcommand{\Prob}{\mathbb{P}}

\newcommand{\calS}{\mathcal{S}}

\newcommand{\calY}{\mathcal{Y}}

\newcommand{\ov}{\mbox{\tiny o}}
\newcommand{\calP}{\mathcal{P}}

\newcommand{\g}{\mbox{\tiny g}}

\newcommand{\f}{\mbox{\tiny f}}

\newcommand{\calF}{\mathcal{F}}

\newcommand{\norms}[1]{\lVert#1\rVert}

\newcommand{\du}{\pmb{d}}
\newcommand{\e}{\pmb{e}}
\newcommand{\dL}{\underline{d}}
\newcommand{\dU}{\overline{d}}
\newcommand{\eL}{\underline{e}}
\newcommand{\eU}{\overline{e}}

\newcommand{\dLb}{\pmb{\dL}}
\newcommand{\dUb}{\pmb{\dU}}
\newcommand{\eLb}{\pmb{\eL}}
\newcommand{\eUb}{\pmb{\eU}}

\newcommand{\db}{\pmb{d}}
\newcommand{\eb}{\pmb{e}}
\newcommand{\xib}{\pmb{\xi}}
\newcommand{\dhat}{\hat{d}}
\newcommand{\ehat}{\hat{e}}

\newcommand{\yhat}{\hat{y}}
\newcommand{\ob}{\pmb{o}}
\newcommand{\gb}{\pmb{g}}

\newcommand{\yb}{\pmb{y}}
\newcommand{\xb}{\pmb{x}}

\usepackage{natbib}
\usepackage[ruled, commentsnumbered]{algorithm2e}

  \journal{ArXiv (under journal review)}

\begin{document}

 \begin{frontmatter}
\title{ \textcolor{black}{Data-driven} Distributionally Robust Surgery Planning in Flexible Operating Rooms Over a Wasserstein Ambiguity}




\author[mymainaddress1]{Karmel S. Shehadeh\corref{cor1}}
\cortext[cor1]{Corresponding author. }
 \ead{kas720[at]lehigh.edu, karmelshehadeh[at]gmail.com}

 \address[mymainaddress1]{Department of Industrial and Systems Engineering, Lehigh University, Bethlehem, PA,  USA}

\begin{abstract}

\noindent We study elective surgery planning in flexible operating rooms (ORs) where emergency patients are accommodated in the existing elective surgery schedule.  \textcolor{black}{Specifically, elective surgeries can be scheduled weeks or months in advance. In contrast, an emergency surgery arrives randomly and must be performed on the day of arrival. Probability distributions of the actual durations of elective and emergency surgeries are unknown, and only a possibly small set of historical realizations may be available.} To address distributional uncertainty, we first construct an ambiguity set that encompasses all possible distributions of surgery durations within a 1-Wasserstein distance from the empirical distribution. We then define  a \textit{\underline{d}istributionally \textit{\textit{r}}obust \underline{s}urgery \underline{a}ssignment} (DSA) problem  \textcolor{black}{to determine optimal elective surgery assignment decisions to available surgical blocks in multiple ORs, considering the capacity needed for emergency cases. The objective is to minimize the total cost consisting of the fixed cost related to scheduling or rejecting elective surgery plus the maximum expected cost associated with OR overtime and idle time over all distributions defined in the ambiguity set.} Using the DSA  model's structural properties, we derive an equivalent mixed-integer linear programming (MILP) reformulation that can be implemented and solved efficiently using off-the-shelf optimization software. \textcolor{black}{In addition, we extend the proposed model to determine the number of ORs needed to serve the two competing surgery classes and derive a MILP reformulation of this extension.  We conduct extensive numerical experiments based on real-world surgery data, demonstrating our proposed model's computational efficiency and superior out-of-sample operational performance over two state-of-the-art approaches. In addition, we derive insights into surgery scheduling in flexible ORs.}


\end{abstract}
\end{frontmatter}

\section{Introduction}


\noindent Operating Rooms (OR) generates about 40--70\% and 20--40\% of hospitals' revenues and operating costs, respectively \citep{bovim2020stochastic, jackson2002business,  li2016improving, viapiano2000operating}.  As an essential area for cost management, planning and scheduling OR activities have received intense research attention \citep{cardoen2010operating, hof2017case, may2011surgical, samudra2016scheduling, Survery2020, zhu2019operating}. In this paper, we focus on elective surgery planning in flexible operating rooms where emergency patients are accommodated in the existing elective surgery schedule. Elective cases can be scheduled weeks or months in advance. In contrast, the arrival of emergency surgeries is random and they must be performed on the day of arrival The probability distributions of the durations of elective and emergency cases are unknown, and only small data on their durations is available.  \textcolor{black}{The goal is to construct a plan that specifies the assignments of a subset of elective cases from a waiting list to available surgical blocks in multiple ORs. The plan's quality is a function of costs related to performing (scheduling) or delaying elective surgeries and costs related to OR overtime and idle time.}

 \color{black}

The concept of flexible ORs differs from dedicated ORs. In the latter type, one or more ORs are solely dedicated to accommodating emergency surgeries, and elective surgery can only be scheduled and performed in any ORs dedicated to its type. Therefore, we do not consider emergency surgeries when planning for elective surgeries in dedicated ORs. In contrast, in flexible ORs, the capacity of the ORs is shared among the two competing surgery classes. Thus, when planning for elective surgery in flexible ORs, we need to account for emergency surgery. Note that hospitals have either flexible or dedicated ORs.  However, most of the existing elective surgery scheduling approaches focus on elective surgery scheduling in dedicated ORs (see Section~\ref{sec2:Literature}).

 \color{black}

\begin{figure}[t!]
 \centering
\begin{subfigure}[b]{0.5\textwidth}
          \centering
        \includegraphics[width=\textwidth]{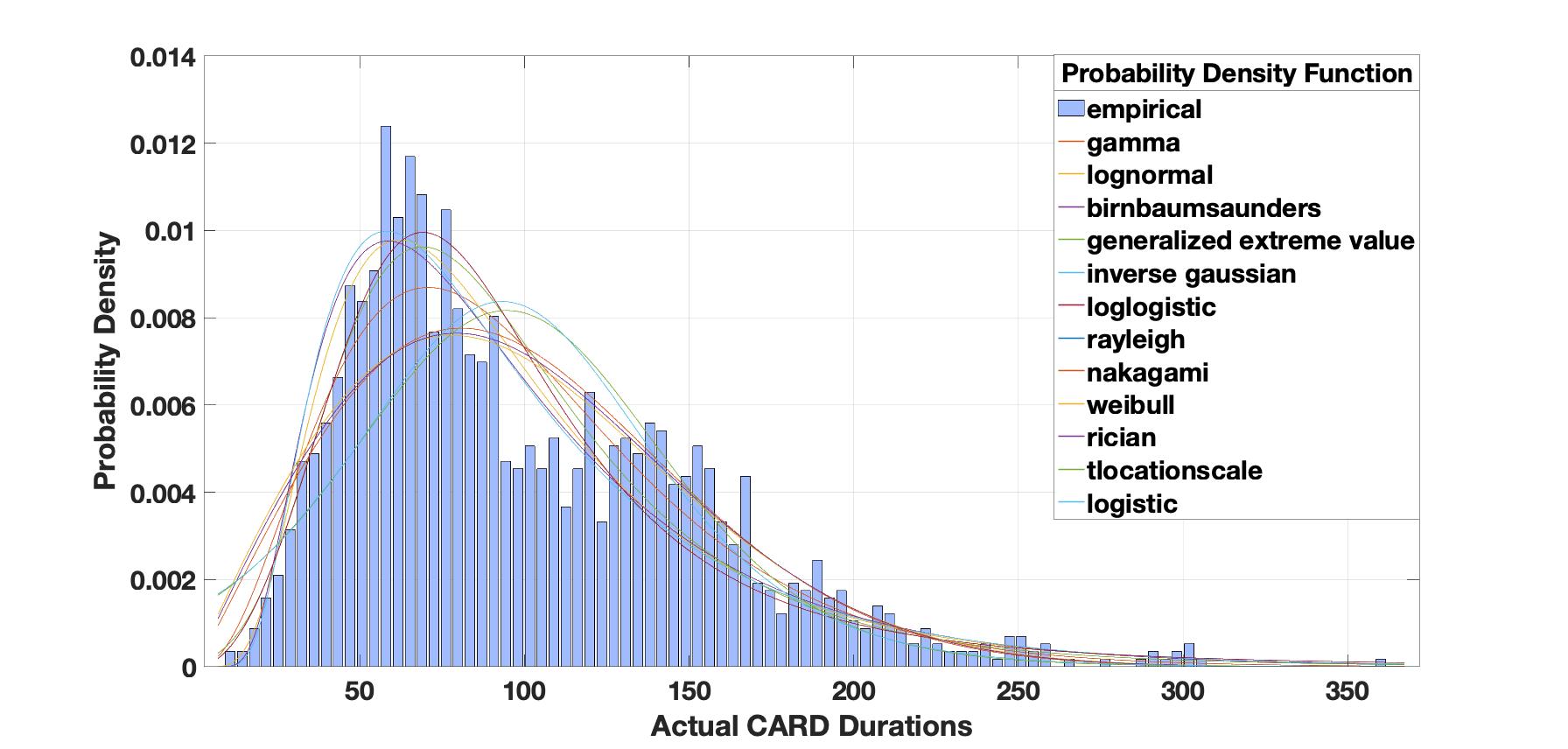}
        \caption{CARD}
    \end{subfigure}%
  \begin{subfigure}[b]{0.5\textwidth}
          \centering
        \includegraphics[width=\textwidth]{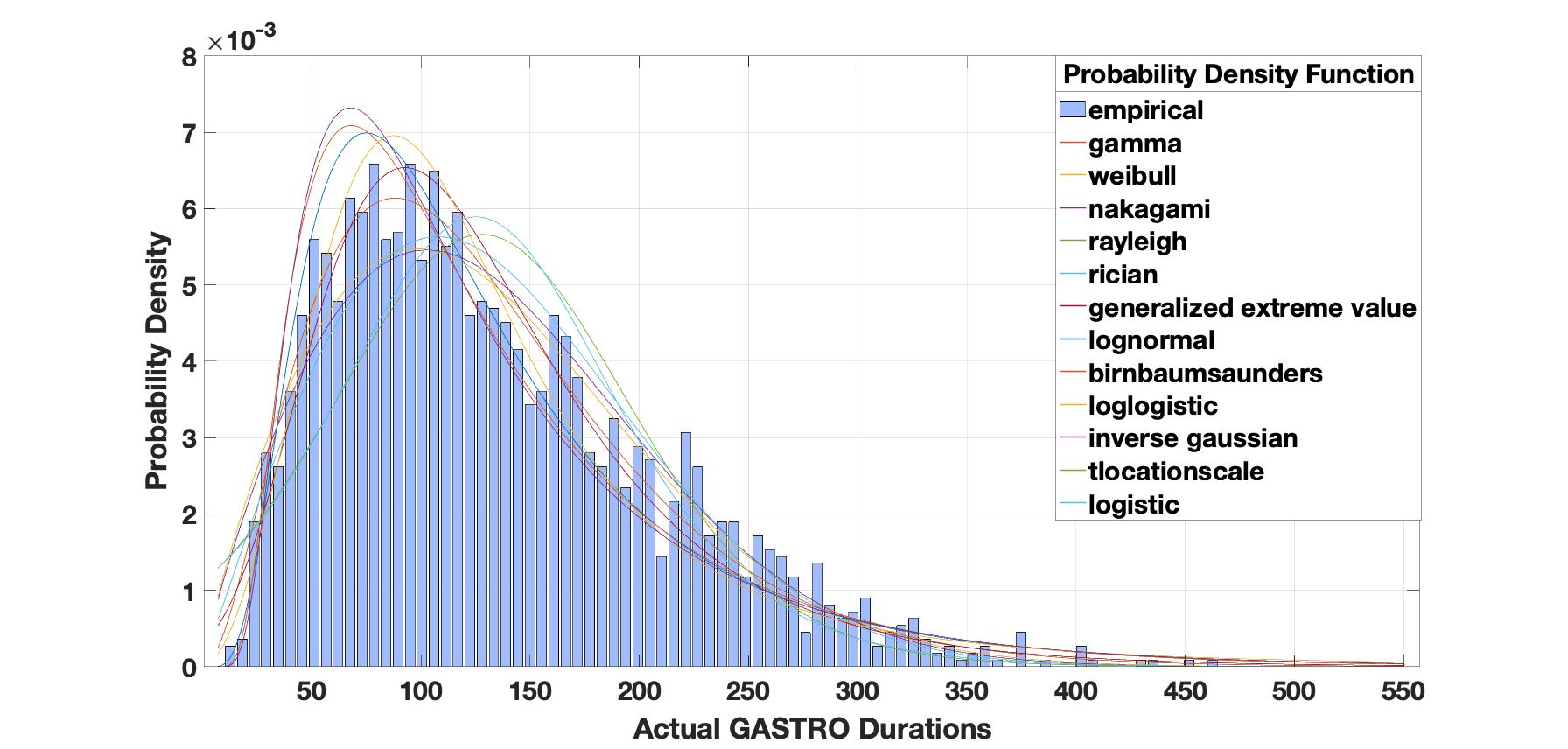}
        \caption{GASTRO}
    \end{subfigure}%
    
 \begin{subfigure}[b]{0.5\textwidth}
          \centering
        \includegraphics[width=\textwidth]{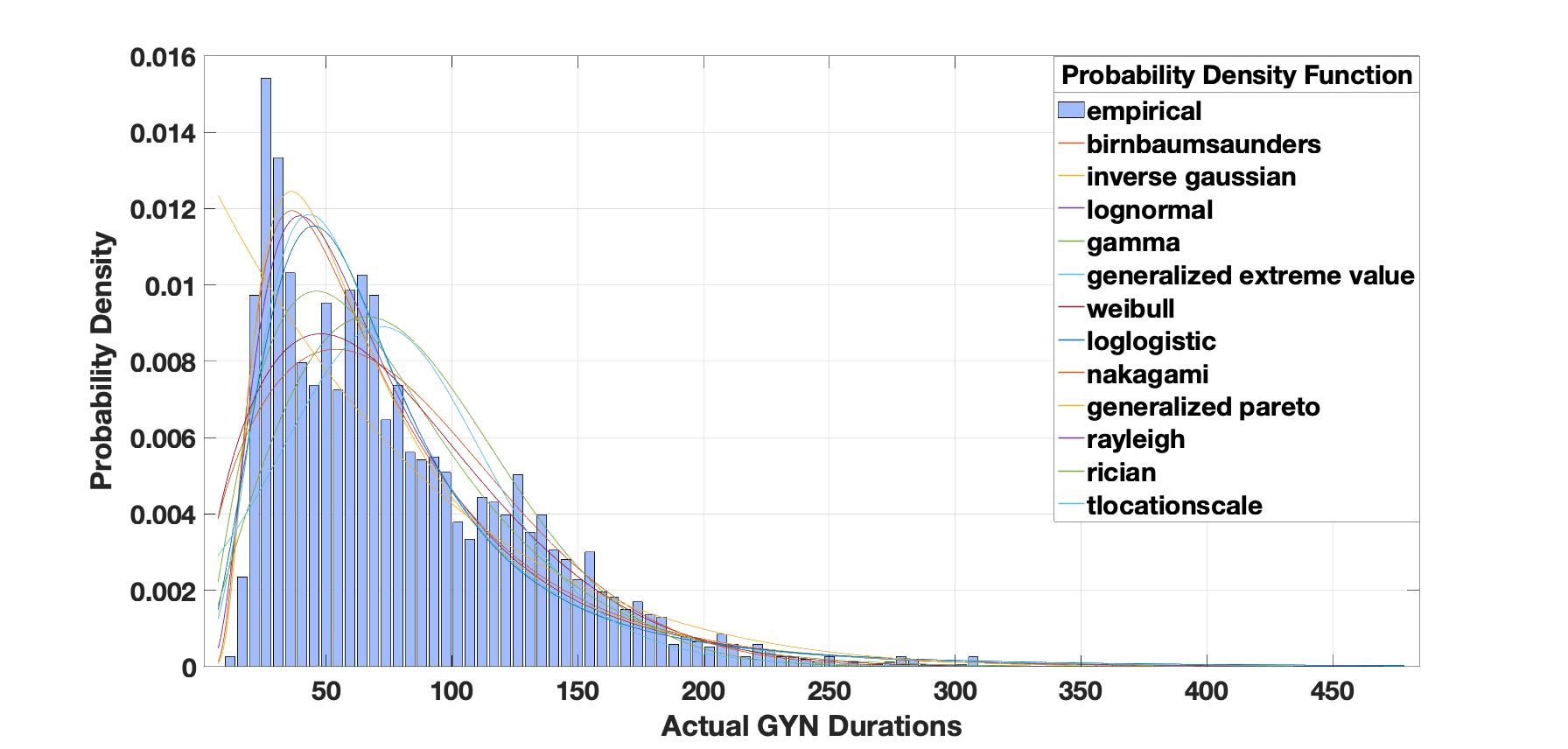}
        \caption{GYN}
    \end{subfigure}%
 \begin{subfigure}[b]{0.5\textwidth}
          \centering
        \includegraphics[width=\textwidth]{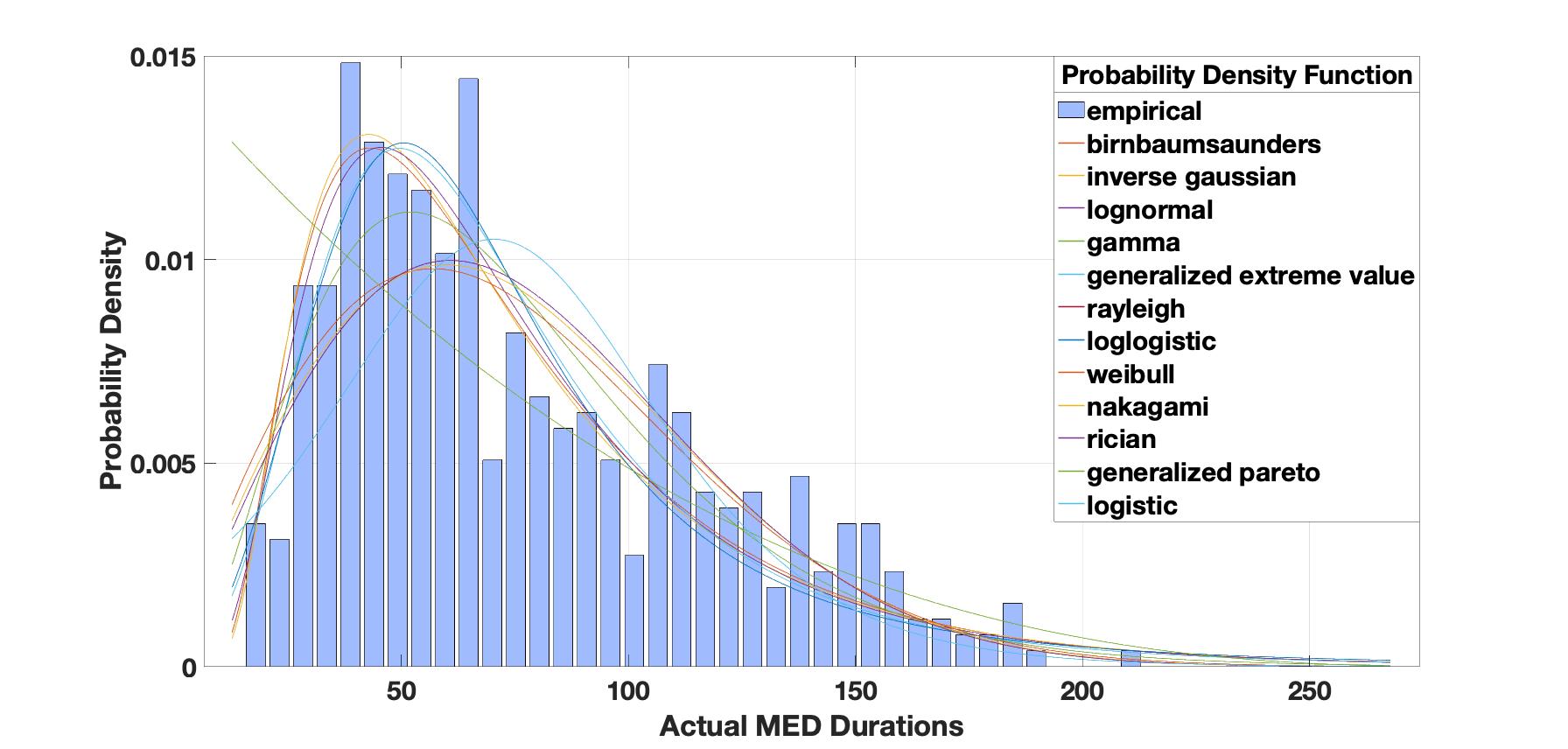}
        \caption{MED}
    \end{subfigure}%

\begin{subfigure}[b]{0.5\textwidth}
          \centering
        \includegraphics[width=\textwidth]{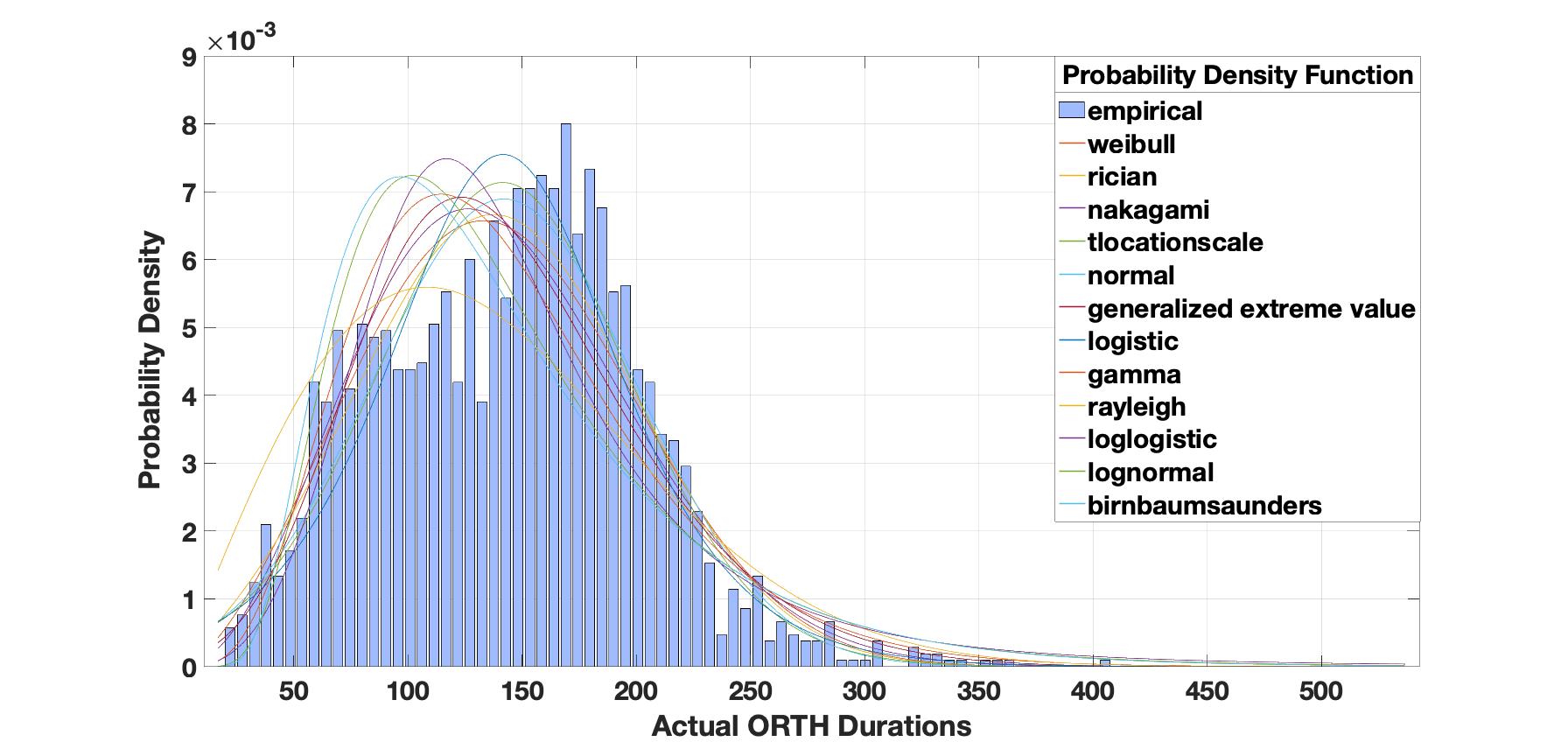}
        \caption{ORTH}
    \end{subfigure}%
\begin{subfigure}[b]{0.5\textwidth}
          \centering
        \includegraphics[width=\textwidth]{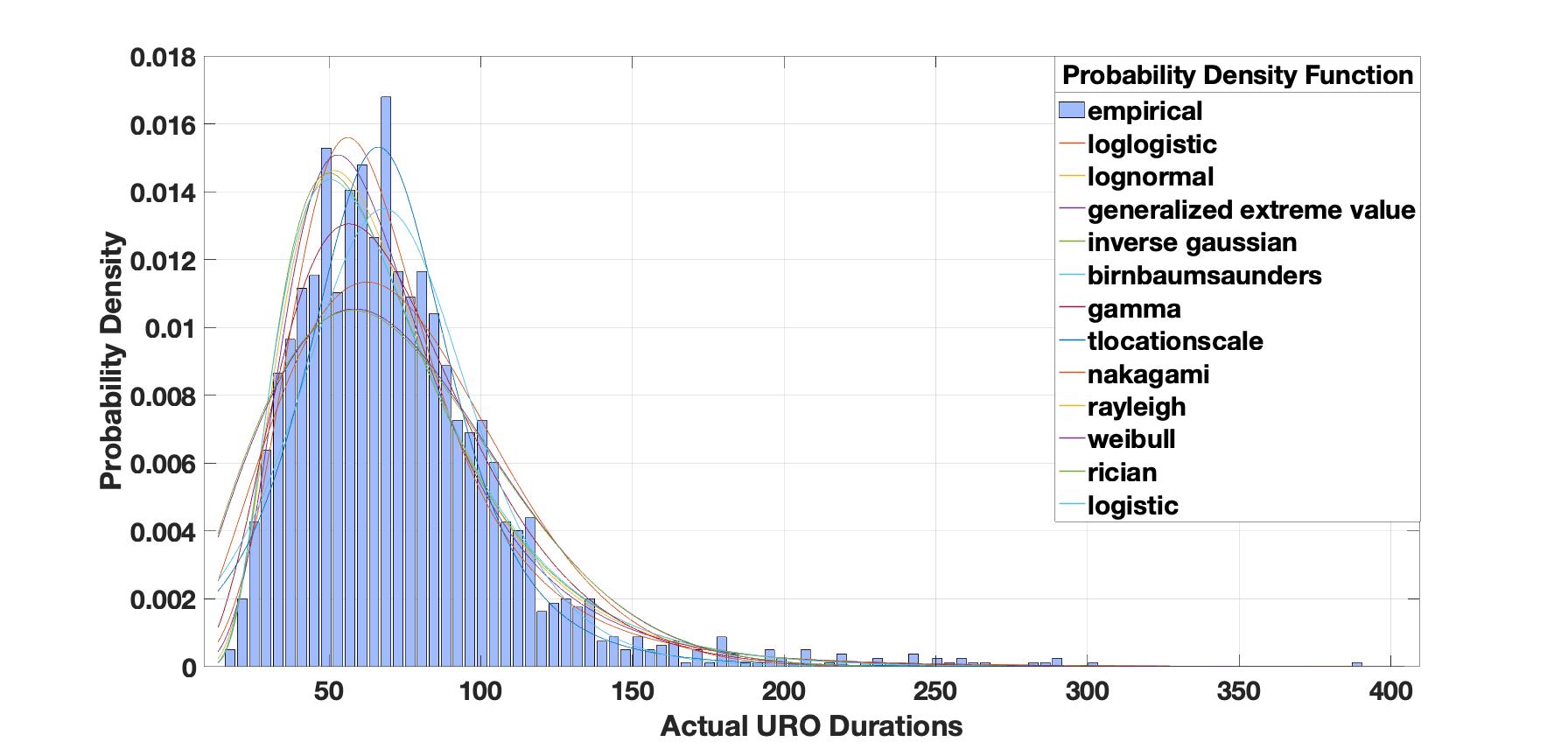}
        \caption{URO}
    \end{subfigure}%
\caption{\textcolor{black}{The empirical and fitted probability distributions for actual surgery duration.}}\label{Fig1:duration}
\end{figure}

Scheduling surgeries in ORs is a complex task primarily due to the significant variability of surgery duration and limited OR capacity. To measure this variability, we use \textcolor{black}{three-years worth of actual surgery durations with respect to six different surgical specialties, namely} General Cardiology (CARD), Gastroenterology (GASTRO), Gynecology (GYN), Medicine (MED), Orthopedics (ORTH), and Urology (URO) (from \cite{mannino2012pattern} and \cite{manninoData}). Figure~\ref{Fig1:duration} presents the empirical and fitted distributions of actual surgery durations data for each specialty. \textcolor{black}{This figure clearly illustrates significant variability in durations within and across surgery types.} Furthermore, there is a wide range of possible probability distributions for modeling the variability (uncertainty) in surgery durations of each type, suggesting distributional ambiguity (i.e., uncertainty of probability distribution).  Such uncertainty and ambiguity are hard to predict and model in advance when elective surgery is scheduled and the OR schedule is constructed.

Ignoring uncertainty and ambiguity may lead to devastating consequences, most notably unpredictable OR utilization, overtime, idle time, and can lead to surgery cancellation and thus sub-optimal quality of care \citep{carello2014cardinality, denton2010optimal, hof2017case,  may2011surgical, shehadeh2019analysis,  wang2014column, wang2019distributionally, xiaoreserved}. To model uncertainty, most of the existing OR scheduling approaches assume that the exact distribution of surgery duration is known (often lognormal) and employ stochastic programming (SP) with  sample average approximation (SAA)   (see, e.g., \cite{dentonopt}, \cite{lamiri2008stochastic}, \cite{shylo2013stochastic}, \cite{batun2011operating}). The SAA approaches  assume that the hospital manager is risk-neutral and evaluates the overtime and idle time costs via the sample average, in which the approximation accuracy improves with the increase of sample size, but the computation becomes cumbersome as well \citep{birge2011introduction, shapiro2009lectures}.

\textcolor{black}{SP remains the state-of-the-art approach to model uncertainty in healthcare scheduling and other application domains. However, the applicability of the SP approach is limited to the case in which we know the distribution of surgery durations of each type. In practice, however, it is unlikely that decision-makers have access to 
a sufficient amount of high-quality data to estimate their probability distributions accurately \citep{viapiano2000operating, keyvanshokooh2020coordinated, Survery2020, wang2019distributionally}}. Even when hospitals collect data on surgery durations, this data may not be readily available due to privacy issues.

In view of possible distributional ambiguity, if we solve a model with a data sample from a biased distribution (as in the SP approach), then the resulting (biased) optimal scheduling decisions may have a disappointing out-of-sample operational performance (e.g., excessive overtime) under the true probability distribution when the schedule is implemented in practice \citep{esfahani2018data}.

\textcolor{black}{Alternatively, we can design a so-called \textit{ambiguity set} (i.e., a family) of all distributions that possess certain partial information about the durations or are close enough to an empirical distribution. Then, we can formulate a distributionally robust optimization (DRO) problem that determines the optimal elective surgery assigning decisions that minimize the sum of patient-related costs plus the maximum expected cost associated with overtime and idle time of ORs over all distributions defined in the ambiguity set. Note that in DRO,  the optimization is based on the worst-case distribution within the ambiguity set, i.e., the distribution is a decision variable \citep{rahimian2019distributionally}.}

\textcolor{black}{DRO has received significant attention recently in many application domains due to the following three primary benefits. First, DRO alleviates the unrealistic assumption of the decision-maker's complete knowledge of distributions. Second, DRO approaches acknowledge the presence of distributional uncertainty. Thus, depending on the ambiguity set used, DRO solutions can anticipate the possibility of disappointing out-of-sample consequences. Third, several studies have proposed techniques to derive DRO models for real-life optimization problems that are more computationally tractable than their SP counterparts, and we have seen many successful applications of DRO to OR and surgery scheduling in recent years (see, e.g., \cite{bansal2021distributionally, keyvanshokooh2020coordinated, jian2017integer, mak2014appointment, shehadeh2020distributionally, ShehPad2020, wang2019distributionally} to name a few).}

\textcolor{black}{The computational tractability of a DRO model depends on the ambiguity set used \citep{rahimi2020comprehensive}. There are various techniques to construct the ambiguity set, and it could be based on moment information \citep{delage2010distributionally} or statistical measures such as the Wasserstein distance \citep{esfahani2018data}.} Most of the early DRO approaches employ moment ambiguity sets, consisting of all distributions sharing certain moments, e.g., the first moment and support. Mean-support ambiguity sets often lead to computationally tractable reformulations (see, e.g., \cite{jian2017integer},  \cite{shehadeh2020distributionally}, and \cite{wang2019distributionally}). However, asymptotic properties cannot be guaranteed because moment-based ambiguity sets only incorporate descriptive statistics (i.e., moments information). \textcolor{black}{New studies have shifted toward distance-based DRO approaches that define the ambiguity set using a distance metric  (e.g., Wasserstein distance) to describe the deviation from a reference (often empirical) distribution.} The main advantage of these approaches is that they enable observed (possibly small-size) data to be incorporated directly and effectively in the optimization problem.

\textcolor{black}{Despite the potential advantages of DRO, there are no moment-based, distance-based, or any other DRO models for elective surgery scheduling in flexible ORs.} This inspires us to study a data-driven DRO approach that can use some partial information on the true underlying distribution of surgery duration from a small set of historical realizations and investigate the properties and performance of such an approach for this specific problem. In particular, we aim to investigate the value of the distance-based DRO approach in modeling uncertainty in the specific OR scheduling problem that we consider in this paper.

\subsection{\textbf{Contributions}}

\noindent In this paper, we address the distributional ambiguity of surgery durations using a data-driven distance-based DRO approach. First, we construct an ambiguity set of probability distributions, incorporating 1-Wasserstein distance and the empirical support set of surgery duration. We then define a data-driven \textit{\underline{d}istributionally robust \underline{s}urgery  \underline{a}ssignment} (DSA) problem, which aim to determine optimal elective surgery assignment decisions to available surgical blocks over a planning horizon to minimize the sum of  patient-related penalty costs (i.e., cost of performing or delaying elective surgery) plus the worst-case (i.e., maximum) expected cost associated with overtime and idle time of the ORs. We take the expectation overall distributions residing in  the data-driven 1-Wasserstein ambiguity set. We summarize our main contributions as follows.
\begin{itemize}\itemsep0em
\item \textcolor{black}{To the best of our knowledge, and according to our literature review in Section~\ref{sec2:Literature}, our paper is the first to propose a data-driven DRO model (denoted as W-DRO) with a 1-Wasserstein ambiguity set for elective surgery scheduling in flexible ORs. In addition, we derive a DRO model (denoted as M-DRO) based on a mean-range ambiguity set and compare the operational performance of the two DRO models, demonstrating where significant performance improvements could be gained with the proposed W-DRO model.}

\item Using DSA structural properties, we derive an equivalent mixed-integer linear programming (MILP) reformulation of \textcolor{black}{our proposed min-max DRO model} that can be implemented and efficiently solved using off-the-shelf optimization software,  not requiring customized algorithmic development or parameter tuning, thereby enabling practitioners to use our  MILP model.  As argued by \cite{shehadeh2019analysis}, such  implementability \textit{is necessary for an optimization-based decision support tool to gain wide adoption in OR and other healthcare systems that do not have ongoing access to support staff with optimization expertise. }

\item To demonstrate the value of a DRO approach for DSA, we conduct an extensive numerical experiment using real-world surgery data. Specifically, our results demonstrate that our proposed W-DRO model (1) yields robust decisions with both small and large data size; (2) enjoys both asymptotic consistency and finite-data guarantees, (3) have a good out-of-sample performance under perfect and misspecified \textcolor{black}{distributional information}, even when a small data set is used in the optimization; (4) is less conservative than the mean-support based DRO model; and (5) is computationally efficient with solution times sufficient for real-world implementation. \textcolor{black}{In addition, we compare the flexible and dedicated ORs policies and show the pros and cons of each in terms of overtime, OR utilization, and access (measured by the number of scheduled surgeries under each policy).}

\item Finally,  we extend our model to include the decisions of how many OR blocks to open. We drive an equivalent MILP of this extension and conduct an experiment demonstrating the superior performance of this extension over the classical SP approach.

\end{itemize}

\subsection{\textbf{Structure of the paper}}

\color{black}
\noindent The remainder of the paper is structured as follows. In Section~\ref{sec2:Literature}, we review the relevant literature.  In Section~\ref{sec3:Defintions}, we detail our problem setting. In Section~\ref{sec:SPmodel}, we present our SP model. In Section~\ref{W_DASBA}, we present and analyze our proposed DRO model. In Section~\ref{Appx:Block}, we extend our proposed DRO model to include the decisions of how many OR blocks to open. In Section~\ref{sec:computational},  we present our numerical experiments and corresponding insights.   Finally, we draw conclusions and discuss future directions in Section~\ref{Sec:Conclusion}.

\color{black}

\section{Relevant Literature}\label{sec2:Literature}

\noindent In this section, we review the literature most relevant to our problem. Specifically, we focus on studies that proposed stochastic optimization approaches for OR planning or surgery scheduling, with two competing classes: elective and emergency.  For comprehensive surveys of OR and surgery planning, we refer to \cite{cardoen2010operating}, \cite{gartner2019flexible}, \cite{hof2017case}, \cite{may2011surgical}, \cite{rahimi2020comprehensive}, \cite{samudra2016scheduling}, \cite{Survery2020}, and \cite{zhu2019operating}. We refer to \cite{pinedo2016scheduling} for a detailed survey of a wide range of scheduling problems, including their theory, algorithms, and applications. We also recognize that there is prior literature that employs simulation to study elective surgery duration and emergency arrivals (see, e.g., \cite{wang2016discrete} and references therein), similar objectives to us in deterministic models (see, e.g., \cite{fei2009solving} and references therein). Our problem is somewhat similar to Newsvendor problem, in which one must find a trade-off between reserving too much versus too little OR time for emergency surgeries (see, e.g., \cite{choi2012handbook} and \cite{nagarajan2014prospect}). In addition to solving the optimization model directly, some studies propose algorithmic and decomposition approaches, including logic-based Benders decomposition approaches (see, e.g., the recent work of \cite{guo2021logic} and references therein). We refer to the surveys mentioned above and references therein for other tactical and strategic decisions. For example, \cite{anjomshoa2018exact} propose an approach for tactical planning and patient selection for elective surgeries motivated by a case study in a local hospital in Melbourne.

Most of the existing stochastic optimization approaches assume that elective and emergency surgeries have dedicated ORs and focus on elective surgery scheduling  (see, e.g., \cite{anjomshoa2018exact, batun2011operating, dentonopt, deng2019chance, guo2021logic, gul2011bi, neyshabouri2017two, shehadeh2019analysis, zhu2019operating}, and the references therein). Dedicating ORs to emergency cases has the apparent advantage of having OR space readily available for emergency cases, which can be handled promptly without impacting the scheduled elective surgeries.  However, the drawback is the low utilization of costly OR resources when emergency cases do not materialize, or emergency arrivals are lower than expected  \citep{xiaoreserved}.  Herein, we focus on flexible OR where emergency patients are accommodated in the existing elective surgery schedule. However, our DRO model can be used for elective surgery scheduling in ORs where emergency cases have dedicated blocks.  \textcolor{black}{In addition, we focus on two-stage models with recourse. Recourse models result when some decisions must be fixed before information relevant to uncertainty is available (e.g., in our problem surgery assignment to OR), while some decisions can be delayed until this information is available (e.g., actual schedule on the day of service and the associated overtime and idle time of the ORs). Two-stage models with recourse have been widely employed to model and solve stochastic (healthcare) scheduling problems \citep{ahmadi2017outpatient, Survery2020}.}

\cite{gerchak1996reservation}, \cite{lamiri2008column}, \cite{lamiri2008stochastic} proposed SP models for the elective surgery planning problem in flexible OR shared between elective and emergency surgeries. The planning problem in \cite{lamiri2008column} is closely related to our work, and consists of assigning elective patients to different ORs over a planning horizon in order to minimize elective-patient-related costs (i.e., cost of performing or delaying surgery to the next planning horizon) and the expected operating rooms' utilization (overtime and idle time) costs. The expectation in \cite{lamiri2008column}'s model is taken with respect to the assumed known distribution of surgery duration.  To solve their SP, \cite{lamiri2008column} proposed a column generation approach, which can produce near-optimal solutions in a short computation time for problems with 12 operating rooms. In contrast, \cite{lamiri2008stochastic} employed Monte Carlo Optimization to solve a sample-average approximation of their SP.

These SP models assume that the hospital manager is risk-neutral and can fully characterize the probability distribution of surgery duration. In practice, it is implausible that hospital managers and other decision-makers have sufficient high-quality data to infer the distribution of random parameters, especially in healthcare \citep{macario2010possible, viapiano2000operating, wang2019distributionally}. Moreover, the distribution of random parameters is often ambiguous \citep{mak2014appointment, qi2017mitigating, wang2019distributionally}. As pointed out by \cite{wang2019distributionally}, various decisions theory, empirical, and neural studies show that decision makers tend to be ambiguity-averse \citep{halevy2007ellsberg, hsu2005neural}.  This is especially the case in healthcare and OR practice due to the high cost of care and indirect costs resulting from adverse outcomes such as overtime, leading to surgery cancellation and OR staff fatigue and, consequently, sub-optimal care \citep{carello2014cardinality, denton2010optimal,hof2017case,  may2011surgical, shehadeh2019analysis, wang2019distributionally, wang2014column, xiaoreserved}.

Various robust approaches have been proposed to model the risk-averse nature of decision-makers and uncertain parameters based on partial information of the distribution. The classical robust optimization (RO) approach model the uncertain parameters based on an uncertainty set of possible outcomes with some structure (e.g., ellipsoid or polyhedron, see,  \cite{bertsimas2004price},  \cite{ben2015deriving} \cite{soyster1973convex}, \cite{rowse2015robust}, \cite{neyshabouri2017two}). In RO, optimization is based on the worst-case scenario within the uncertainty set. By optimizing based on the worst-case scenario, classical RO approaches may yield over-conservative solutions and sub-optimal decisions for the other more likely to be observed scenarios \citep{chen2019robust, delage2018value, wang2019distributionally}.

In contrast to RO, the DRO approach assumes that the distribution of random parameters resides in a so-called ambiguity set, i.e., a family of all possible distributions that share some common properties of random parameters. Accordingly, optimization is based on the distributions residing within the ambiguity set, i.e., the probability distribution of a random parameter is a decision variable in DRO.   Most early DRO studies employ the moment ambiguity set, which consists of all distributions sharing particular moments (e.g., the first and second moments).  Notably,  \cite{wang2019distributionally} is the first (and only study) to employ DRO for elective surgery scheduling and surgery block allocation. Using an ambiguity set based on mean, mean absolute deviation, and support of surgery duration,  \cite{wang2019distributionally} propose a DRO model that determines the number of ORs to open and assign the surgeries in a daily listing to the open ORs.  \cite{wang2019distributionally}' DRO model aims to minimize OR opening costs and expected OR overtime cost.  \cite{wang2019distributionally} did not consider emergency surgery  and ignore idle time costs, which is an important metric of OR efficiency and utilization \citep{girotto2010optimizing, jebali2015stochastic, liu2019integrated, ShehPad2020}.

  Moment-support-based DRO approaches often lead to tractable reformulations, but they known to have weak convergence properties. New studies have shifted toward distance-based DRO approaches that construct ambiguity sets in the vicinity of a reference (empirical) distribution.  The distance-based DRO models' main advantage is that they enable possibly small-size data to be incorporated directly and effectively in the optimization problem. In this paper, we construct a 1-Wasserstein distance-based ambiguity set of all probability distributions of  surgery durations. Various authors have shown that 1-Wasserstein ambiguity admits tractable reformulation in many real-world applications (see, e.g., \cite{jiang2019data} \cite{saif2020data}). Then, we propose a DRO model for elective surgery planning in flexible ORs. The objective is to minimize the sum of patient-related penalty costs (i.e., cost of performing or delaying elective surgery) and the expectation of OR overtime and idle time cost.  We take the expectation over all distributions residing in the data-driven 1-Wasserstein ambiguity set. We derive an equivalent MILP reformulation of our min-max model. Our MILP can be implemented and solved directly with standard optimization software packages, not requiring customized algorithmic development or tuning, thereby enabling practitioners to use our  MILP.

   We use our MILP and real-world surgery data to conduct extensive computational experiments, demonstrating our approach's advantages.  Finally, we extend our model to include the decision of how many OR blocks to open, derive an equivalent MILP of this extension, and present an experiment demonstrating the performance of this extension. To the best of our knowledge, and according to recent survey papers on OR and surgery scheduling \citep{zhu2019operating}, our paper is the first to propose and analyze DRO approaches for this specific surgery and OR planning problem.

\color{black}

\section{Problem Description}\label{sec3:Defintions}
\color{black}

\noindent \textcolor{black}{Let us now start to introduce our problem. We consider a hospital with $B$ available surgical blocks (or ORs) over an arbitrary planning horizon of $T$ days (e.g., a week).  Each surgical block $b \in [B]:=\{1\ldots,B\}$ is assigned to a particular OR, has a pre-allocated length of time $\mathcal{L}_b$, and is dedicated to only one type of surgical specialty.  For any $b \in [B]$, the time length $\mathcal{L}_b$ is typically long enough so that  multiple surgeries can be performed during that time. Note that there can be multiple blocks for the same specialty during a cycle of the OR schedule. The OR capacity is shared among two competing surgery classes: a known number $I$ of elective surgeries that are to be planned in advance, and random emergency surgeries that must be performed on the day of arrival.}

\textcolor{black}{Each elective surgery $i \in [I]$ has a surgery type and can be assigned to any of the blocks dedicated to the corresponding surgery  type during the planning horizon. Associated with each elective case, there is a cost for performing (scheduling) and for rejecting the surgery (i.e., cost of delaying surgery to the next planning horizon). Let's assume that patients who are not assigned to any surgery block during the current planning horizon are assigned to a dummy block $b^\prime \notin [B]$. Furthermore, let $c_{i,b}$ and $c_{i,b'}$ represent the costs of performing and postponing surgery $i \in [I]$, respectively. Then, it is reasonable to assume that $c_{i,b'}>c_{i,b},$ for all  $b \in [B]$ (a common assumption in the elective surgery scheduling literature, see, e.g.,  \cite{jebali2015stochastic}, \cite{min2010scheduling}, \cite{ ShehPad2020}, and \cite{zhang2019two}). Moreover, the cost of assigning a surgery to a block is block-dependent to take into account the surgery's waiting time on the list and clinical priority, among other potential factors \citep{jebali2015stochastic, lamiri2008column,lamiri2008stochastic}.} However, determining this cost is out of the scope of this paper.

Elective surgery duration ($d_i$, $i \in [I]$)  are random and depend on surgery type. As defined in \cite{lamiri2008column} and  \cite{lamiri2008stochastic}, the total random surgery duration or capacity, $e_b$, needed for emergency cases in each surgical block $b \in [B]$ is random. The probability distributions of these random parameters are ambiguous (unknown).

Given a waiting list of $I$ elective surgeries and their types, we are interested in determining a plan that specifies the number (or subset) of elective surgeries to schedule in each surgery block (equivalently, surgery assignments to the available surgery blocks).  The plan's quality is a function of costs related to performing (scheduling) or not performing (rejecting) elective surgeries and costs related to OR overtime and total idle time. Overtime occurs when surgeries assigned to block $b$ are not completed within [0,  $\mathcal{L}_b$]. Idle time occurs when surgeries assigned to block $b$ are completed before $\mathcal{L}_b$. \textcolor{black}{Note that overtime and idle time are random (second-stage) performance metrics.  That is, they are observed after the realization of random parameters while the planned schedule is being executed. In contrast, costs related to scheduling or not scheduling an elective surgery are fixed planning (first-stage) costs incurred when the OR schedule is constructed}. We make the following assumption based on prior studies: 
  \begin{itemize}\itemsep0em  
\item[A1.] The set of elective patients, their surgery types, and priorities are known with certainty (a standard assumption in the elective surgery scheduling literature, see, e.g., \cite{jebali2015stochastic}, \cite{ lamiri2008column}, \cite{lamiri2008stochastic}, \cite{ min2010scheduling},  and \cite{neyshabouri2017two}). 
\item[A2.] The planning horizon $T$ is an integer multiple of the surgery schedule cycle length.
\item[A3.] Surgical blocks and their assignments to operating rooms during a cycle of surgery schedule are previously determined \citep{lamiri2008column,lamiri2008stochastic,  neyshabouri2017two}.
\item[A4.] For ease of modeling,  we assume that patients who are not assigned to any surgery block during the current planning horizon are assigned to a dummy block $b^\prime$ \citep{ lamiri2008column,lamiri2008stochastic}.
\item[A5.] Patients of the same type have identical probability distributions of surgery duration  (a common assumption in prior studies).
\item[A6.] We assume that we know the lower ($\dLb, \eLb$) and upper ($\dUb$, $\eUb$) bounds on $\db$ and $\eb$ (this is a mild and realistic assumption that holds in healthcare scheduling and mimic the OR practice). Mathematically, we consider support $\calS=\calS^{\tiny d} \times \calS^{\tiny e}$, where $\calS^{\tiny d}$ and $\calS^{\tiny e}$ are respectively the supports of random parameters $d$ and $e$, and defined as follows:
\begin{align}\label{eq:support}
&  \calS^{\tiny d}:=\left\{ \db  \geq 0: \begin{array}{l} \dL_i
\leq d_i \leq \dU_i,  \forall i \in [I] \end{array}  \right\},\\
&  \calS^{\tiny e}:=\left\{\eb \geq 0: \begin{array}{l} \eL_b \leq e_b \leq \eU_b,  \forall b \in [B]  \end{array}  \right\}.
\end{align}
\end{itemize}

 For modeling convenience, prior literature suggests that one can avoid excessive overtime by assigning a high weight ($c_b^{\ov}$) on overtime in the objective and/or setting $\mathcal{L}_b$=(planned length of $b-\Delta$) with $\Delta>0$. Our model can accommodate such constraints. 
 
 \vspace{1mm}

\noindent \textbf{Additional Notation:} For $a,b \in \mathbb{Z}$, we define $[a]:=\lbrace1,2,\ldots, a \rbrace$ and $[a,b]_\mathbb{Z}:=\lbrace c \in \mathbb{Z}: a \leq c \leq b \rbrace$, i.e., $[a,b]$ represent the set of running integer indices $\{a, a+1, a+2, \ldots, b \}$. The abbreviations ``w.l.o.g.'' and ``w.l.o.o.'' respectively represent ``without loss of generality'' and ``without loss of optimality.''  We use boldface notation to denote vectors, e,g., $\pmb{d}:=[d_1, d_2, \ldots, d_N]^\top$. A complete listing of the parameters and decision variables can be found in Table~\ref{table:notation}.

\begin{table}[t]  
\small
\center
   \renewcommand{\arraystretch}{0.9}
  \caption{Notation.} 
\begin{tabular}{ll}
\hline
\multicolumn{2}{l}{\textbf{Parameters and sets}} \\
$I$ & number, or set, of surgery  \\
$B$ & number, or set, of OR blocks \\
$c_{i,b}$ & cost of assigning surgery $i$ to block $b$\\
$c_b^{\ov}$ & unit overtime cost for each block $b$\\
$c_b^{\g}$ & unit idle cost for each block $b$\\
$d_i$ & duration of surgery $i$\\
$\dL/\dU_i$ & lower/upper bounds of $d_i$   \\
$e$ & capacity needed for emergency cases \\
$\eL_i/\eU_i$ & lower/upper bounds of $e_i$   \\
$\mathcal{L}_b$ & capacity or planned length of block $b$\\
 \multicolumn{2}{l}{\textbf{First-stage decision variables } } \\
$y_{i,b}$ &   $\left\{\begin{array}{ll}
1, & \mbox{surgery }  i  \text{ is assigned to OR block } b     \\
0, & \mbox{otherwise.}
\end{array}\right.$ \\
 \multicolumn{2}{l}{\textbf{Second-stage decision variables } } \\
 $o_b$ & continuous decision variable capturing overtime in block $b$ \\
  $g_b$&  continuous decision variable capturing idle time in block $b$ \\
\hline
\end{tabular}\label{table:notation}
\end{table} 


\section{\textbf{Stochastic Programming Model}}\label{sec:SPmodel}
\noindent \textcolor{black}{In this section, we present a two-stage SP formulation of the problem that assumes that the probability distributions of surgery durations are known.} \textcolor{black}{First, let us introduce the variables and constraints defining the first-stage of this SP model. For each $i \in [I]$ and $b \in [B]\cup \{b'\}$, we define a binary decision variable $y_{i,b}$ that equals 1 if we assign elective surgery $i$ to block $b$, and is zero otherwise.}  We define feasible region $\calY$ of $\yb$ in \eqref{eq:RegionY}  such that each elective surgery is assigned to one surgical block with the obvious convention that a surgery assigned to the dummy block $b'$ (i.e., $y_{i,b'}=1$) is postponed to the next planning period (i.e., rejected in the current planning period). 
\begin{align}
\calY&:=\left\{ \yb:  \begin{array}{l} \sum \limits_{b \in [B] \cup \lbrace b^\prime \rbrace} y_{i,b} =1, \  \forall i \in [I],  \ y_{i,b} \in \lbrace 0, 1 \rbrace, \  \forall i \in [I], b \in  [B] \cup \lbrace b^\prime \rbrace \end{array} \right\}.  \label{eq:RegionY}
\end{align} 

\textcolor{black}{Let us now introduce the parameters and variables defining our second-stage problem. For each $i \in [I]$}, we let $d_i$ represents the duration of elective surgery $i$. For all $b \in [B]$, we define $e_b$ as the  duration of emergency surgery in block $b$ (i.e.,  capacity needed for emergency cases as defined in  \cite{lamiri2008column, lamiri2009optimization}). \textcolor{black}{For all $b \in [B]$, we define nonnegative continuous decision variables $o_b$ and $g_b$ respectively to represent block $b$ overtime and idle time.} For all $b \in [B]$, we define $c_b^{\ov}$ and $c_b^{\g}$ as the nonnegative unit penalty costs of overtime and  idle time, respectively. Finally, we define $\Prob_\xi$ as the the probability distribution of $\xi=[\db, \eb]^\top$, and $\E_{\Prob_\xi}$ as the expectation under distribution $\Prob_\xi$.

\textcolor{black}{The SP model can now be stated as follows  \citep{lamiri2008stochastic, lamiri2009optimization}:}
\begin{align} 
(\text{SP}) \ \ \ \hat{Z}&= \min \limits_{\yb \in \calY} \Bigg \{ \sum_{i\in I} \sum_{b \in B\cup\{b'\}} c_{i,b} y_{i,b} +\E_{\Prob_\xi} [f(\yb,\xi)] \Bigg \},   \label{SP}
\end{align}

where for a feasible $\yb \in \calY$ and a joint realization of uncertain parameters $\xi=[\db, \eb]^\top$, we compute the random cost associated with  overtime and idle time using the following linear program (LP).
\allowdisplaybreaks
\begin{subequations}\label{2ndstage}
\begin{align}
  f(\yb, \xi):=   \min_{\ob,\gb} &\sum_{b \in B} \big( c_b^{\ov}o_b+ c_b^{\g}o_g \big) \label{2nd_Obj}\\
  \text{s.t.} & \ o_b-g_b=\sum_{i \in I} d_iy_{i,b}+e_b-\mathcal{L}_b, && \forall b \in [B],\label{2nd_Const1}\\ 
  			& \ (o_b,g_b) \geq 0, && \forall b \in [B]. \label{2nd_Const2}
\end{align} 
\end{subequations}

Formulation \eqref{SP} aim to find first-stage scheduling decisions $\yb \in \calY$ that minimize the first-stage cost (i.e., costs related to performing/scheduling or not performing/rejecting elective surgeries) plus the expected cost associated with overtime and idle time.  Constraint \eqref{2nd_Const1} yield either block overtime or idle time based on the durations of surgeries preformed in this block.  Finally, constraints \eqref{2nd_Const2} specify feasible ranges of the second-stage decision variables. \textcolor{black}{It is easy to verify that the second-stage (recourse) formulation \eqref{2ndstage} is feasible for any feasible first-stage decisions. Thus, we have a relatively complete recourse.}

\section{\textbf{DRO model over Wasserstein ambiguity set (W-DRO)}}\label{W_DASBA}
\noindent \textcolor{black}{In this section, we present our proposed DRO models for the DSA  that does not assume that the probability distributions of surgery durations are known. Specifically, we consider the case that the joint distribution $\Prob_\xi$ of $\xi:=[\pmb{d}, \pmb{e}]^\top$ may be observed via a possibly small finite set $\{\hat{\xi}^1, \ldots, \hat{\xi}^N\}$ of $N$ samples, which may come from the limited historical realizations or a reference empirical distribution.} Accordingly, we construct an ambiguity set based on 1-Wasserstein distance  (i.e., we use $\ell_1$--norm in the definition of Wasserstein metric), which allows us to derive a tractable model. Then, we formulate our DRO model (denoted as W-DRO) using this  ambiguity set.

\textcolor{black}{First, let us define the 1-Wasserstein distance}. Suppose that random vectors $\xi_1$ and $\xi_2$ follow $\mathbb{Q}_1$ and $\mathbb{Q}_2$, respectively, where distributions $\mathbb{Q}_1$ and $\mathbb{Q}_2$ are defined over the common support $\calS$. The 1-Wasserstein distance $\text{dist}(\mathbb{Q}_1, \mathbb{Q}_2)$ between $\mathbb{Q}_1$ and $\mathbb{Q}_2$ is the cost of an optimal transportation plan for moving the probability mas in one so it becomes identical to the other, where the cost of moving from $\xi_1$ to $\xi_2$ equals to the norm $||\xi_1-\xi_2||$. Mathematically, 
\begin{equation}\label{W_distance}
\text{dist}(\mathbb{Q}_1, \mathbb{Q}_2):= \inf_{\Pi \in \calP(\mathbb{Q}_1, \mathbb{Q}_2)}  \Bigg \{\int_{\calS} ||\xi_1-\xi_2|| \ \Pi(\text{d}\xi_1, \text{d}\xi_2) \Bigg | \begin{array}{ll} & \Pi \text{ is a joint distribution of }\xi_1 \text{ and } \xi_2\\
& \text{with marginals } \mathbb{Q}_1 \text{ and }  \mathbb{Q}_1, \text{ respectively} 
 \end{array} \Bigg \}
\end{equation}

where $\calP(\mathbb{Q}_1, \mathbb{Q}_2) $ is the set of all joints distributions of ($\xi_1$, $\xi_2$) supported on $\calS$ with marginals $\mathbb{Q}_1$ and $\mathbb{Q}_2$.  Accordingly, we construct the following $1$-Wasserstein ambiguity set:
\begin{align}\label{W_ambiguity}
\calF (\hat{\mathbb{P}}_\xi^N, \epsilon)= \left\{ \mathbb{Q}_\xi \in \calP(\calS) \middle| \begin{array}{l} \text{dist}(\mathbb{Q}_\xi, \hat{\mathbb{P}}_\xi^N) \leq \epsilon\end{array} \right\},
\end{align}

where $ \calP(\calS)$ is the set of all \textcolor{black}{joint} probability distributions supported on $\calS$, $\hat{\mathbb{P}}_\xi^N=\frac{1}{N} \sum_{n=1}^N \delta_{\hat \xi^n}$ is the empirical distribution of $\xi$ based on $N$ i.i.d samples, and $\epsilon >0$ is the radius of the ambiguity set. The set $\calF (\hat{\mathbb{P}}_\xi^N, \epsilon)$ represent a Wasserstein ball of radius $\epsilon$ centered at the empirical distribution $\hat{\mathbb{P}}_\xi^N$. \textcolor{black}{Note that we do not make any assumptions on elective and emergency surgeries durations (i.e., $\xi$). Thus, they can be dependent, correlated, or independent.}

 Using the ambiguity set $\calF(\hat{\mathbb{P}}_\xi^N, \epsilon)$, we formulate our W-DRO model for DSA as follows:
\begin{align}
(\text{W-DRO}) \ \ \ \hat{Z}(N, \epsilon)&= \min \limits_{\yb \in \calY} \Bigg \{ \sum_{i\in I} \sum_{b \in B\cup\{b'\}} c_{i,b} y_{i,b} + \sup_{\mathbb{P}_{\xi} \in \calF (\hat{\mathbb{P}}_\xi^N, \epsilon) }\mathbb{E}_{\mathbb{P}_\xi} [f(\yb,\xi) ] \Bigg \}.  \label{DASBA_Obj}
\end{align} \label{DASBA}

\textcolor{black}{Formulation \eqref{DASBA} finds first-stage decisions $\yb \in \calY$ that minimize the first-stage cost plus the maximum expectation of the second-stage cost,  where the expectation is taken over all distributions residing in $\calF (\hat{\mathbb{P}}_\xi^N, \epsilon)$.}

 We often seek asymptotic consistency in data-driven applications such as W-DRO. In particular, one expects that as the sample size $N$ increases to infinity, the optimal value of problem \eqref{DASBA} $\hat{Z}(N, \epsilon)$ converges to $Z^*$ (the optimal value of the  SP model in \eqref{SP} with perfect knowledge of $\Prob_\xi$).  Accordingly, as $N$ increases, an optimal solution $\yb$ to  W-DRO tends to the optimal solution of problem \eqref{SP}. Additionally, if $\hat{Z}(N, \epsilon)>Z^*$ almost surely, then W-DRO provides a safe upper bound guarantee on the expected total cost with any finite data size $N$. Assumption A6 indicates that the support set $\calS$ is non-empty, convex, and compact. As such, we can make use of existing theory in establishing the asymptomatic consistency and finite sample guarantee of our W-DRO model. In particular, given assumption A6, lemma 1 in \cite{jiang2019data} and Theorem 2 in \cite{fournier2015rate} assure that the $\ell_p$-Wasserstein ambiguity set incorporates the true distribution $\Prob_\xi$ of $\xi$ with high confidence. Theorem 1 in \cite{jiang2019data} and Theorem 3.6 in \cite{esfahani2018data} assure asymptotic consistency, and Theorem 2 in \cite{jiang2019data} and Theorem 3.5 in \cite{esfahani2018data} assure finite-data guarantee of W-DRO. We refer to these papers for detailed proofs and discussion on these results. For completeness, below, we provided the adapted Lemma 1, Theorem 1, and Theorem 2  and provide customized proofs to fit our problem in \ref{Appix:Proof_Lemma1}--\ref{Appix:Proof_Thrm2}.

\color{black}

\vspace{2mm}

\noindent  \textbf{Lemma 1.} (Adapted from \cite{jiang2019data} and Theorem 2 in \cite{fournier2015rate}). \textit{Suppose that Assumption A6 holds. Then there exist none-negative constants $c_1$ and $c_2$ such that, for all $N\geq 1$ and $\beta \in (0, \min \{ 1, c_1\} )$,}
\begin{equation*}
\Prob^N \big \{ \text{dist} (\Prob_\xi, \hat{\Prob}_\xi^N) \leq \epsilon_N(\beta) \} \geq 1-\beta,
\end{equation*}
\textit{where $\Prob_\xi^N$ represents the product measure of $N$ copies of $\Prob_\xi$ and $\epsilon_N(\beta)= \Big [ \frac{\log(c_1\beta^{-1})}{c_2N}\Big]^{\frac{1}{\max \{ 3p, n\}}}$ (see \ref{Appix:Proof_Lemma1} for a proof).}

As detailed in \cite{jiang2019data}, Lemma 1 assures that the $p$-Wasserstein ambiguity ($1$-Wasserstein ambiguity in our problem) contains the true distribution $\Prob_\xi$ with high confidence. 

\vspace{2mm}

\noindent \textbf{Theorem 1.} (Asymptotic consistency, adapted from \cite{jiang2019data} and Theorem 3.6 of \cite{esfahani2018data}). \textit{Suppose that Assumption A6 holds. Consider a sequence of confidence levels $\{ \beta_N\}_{N \in \mathbb{R}}$ such that $\sum_{N=1}^\infty \beta_N <\infty $ and $\lim_{N\rightarrow \infty} \epsilon_N(\beta_N)=0$, and let $(\hat{y} (N, \epsilon_N(\beta_N))$ represents an optimal solution to W-DRO with the ambiguity set $\calF (\hat{\Prob}_\xi^N, \epsilon_N(\beta_N))$. Then, $\Prob_\xi^\infty$- almost surely we have $\hat{Z}(N, \epsilon_N(\beta_N))\rightarrow Z^*$ as $N\rightarrow \infty$. In addition, any accumulation points of  $\{\hat{y} (N, \epsilon_N(\beta_N) \}_{N \in \mathbb{N}}$ is an optimal solution of \eqref{SP} $\Prob_\xi^\infty$- almost surely (see \ref{Appx:Proof_Thrm1} for a proof).}

\vspace{2mm}

\noindent \textbf{Theorem 2.} (Finite-data guarantee, adapted from \cite{jiang2019data} and Theorem 3.5 in \cite{esfahani2018data}).  \textit{For any $\beta \in (0, 1)$, let $\hat{y}(N, \epsilon_N(\beta_N)),$ represent an optimal solution of W-DRO with ambiguity set $\calF (\hat{\Prob}_\xi^N, \epsilon_N(\beta_N))$. Then,} 
$$\Prob_{\xi}^N \big \{  \E_{\Prob_\mathbf{\xi}} [f (  \hat{y}(N, \epsilon_N(\beta_N),\xi )] \leq \hat{Z} (N, \epsilon_N(\beta_N))\big\} \geq 1- \beta.$$ \textit{(see \ref{Appix:Proof_Thrm2} for a proof).}

\color{black}
\subsection{\textbf{Reformulation of the W-DRO model}}\label{sec:reformualtion}
\noindent  \textcolor{black}{Recall that $f(\cdot)$ is defined by a minimization problem. Thus, in \eqref{DASBA}, we have an inner max-min problem. As such, it is not straightforward to solve the W-DRO model in \eqref{DASBA}.  In this section, we derive an equivalent reformulation of the W-DRO model that is solvable.} First, in Proposition~\ref{Prop1}, we present an equivalent dual formulation of the inner maximization problem in \eqref{DASBA}  (see \ref{Proof_Prop1} for a detailed proof).
\begin{prop}\label{Prop1} 
For fixed $\yb \in \calY$, problem $\sup_{\Prob_{\xi}\in \calF (\hat{\mathbb{P}}_\xi^N, \epsilon) }\mathbb{E}_{\mathbb{P}_\xi} [f(\yb, \xi) ] $ in \eqref{DASBA}  is equivalent to 
\begin{subequations}\label{DualOfInner}
\begin{align}
\inf_{\rho} &  \Big\{ \epsilon \rho + \frac{1}{N} \sum_{n=1}^N \sup_{\xi \in \calS}  \{ f(\yb, \xi)-\rho || \xi -\hat{\xi}^n || \}  \Big\} \label{DualOfInner:Obj1}\\
\text{s.t.}  &  \  \rho \geq 0. \label{DualOfInner:C2}
\end{align}
\end{subequations}
\end{prop}

Again, formulation \eqref{DualOfInner} is not directly solvable using standard techniques. Fortunately, given that the supports of $\pmb{d}$ and $\pmb{e}$ are rectangular and finite (Assumption A6), we next show that we can reformulate each of the inner maximization problems in \eqref{DualOfInner} as a linear program (LP) for fixed $\rho \geq 0$ and $\yb \in \calY$.  First, we observe that for each $\yb \in \calY$, the recourse problem can be decomposed by each block, i.e.,  $f(\yb, \xi)=\sum_{b \in B} f_b(y,\xi)$, where for each $b \in B$:
 \begin{subequations}\label{2nd_block}
\begin{align}
 f_b(y, \xi):= \min_{o_b, g_b} &  \ \Big( c_b^{\ov}o_b+ c_b^{\g}o_g \Big) \label{2nd_block_Obj}\\
  \text{s.t.} & \ o_b-g_b=\sum_{i \in I} d_iy_{i,b}+e_b-\mathcal{L}_b, \label{2nd_block_Const1}\\ 
  					& \ (o_b,g_b) \geq 0. &&  \label{2nd_block_Const2}
\end{align} 
\end{subequations}

Let $\beta_b$ be the dual variables associated with constraint \eqref{2nd_block_Const2}. The dual of the LP in \eqref{2nd_block} is as follows
\allowdisplaybreaks
\begin{subequations}\label{Dual_2nd}
\begin{align}
 f_b(\yb, \xi) := \max_{\beta}  & \Big ( \sum_{i \in I} d_iy_{i,b}+e_b-\mathcal{L}_b \Big ) \beta_b \label{Dual_Obj}\\
\text{s.t. } & \mathcal{B}:=\{-c_b^{\g} \leq \beta_b \leq c_b^{\ov}  \}. 
\end{align}
\end{subequations}

 Next,  for fixed $\yb \in \calY$ and $\rho \geq 0$,  we define $Q^n(\rho, \yb)= \sup_{\xi \in \calS}  \{ f(\yb, \xi)-\rho | \xi -\hat{\xi}^n |  \}$, for all $n=1,\ldots,N$. Given the dual formulation of $f_b(y,\xi)$ with feasible region $ \mathcal{B}:=\{-c_b^{\g} \leq \beta_b \leq c_b^{\ov}  \} $, we write the problem of computing  $Q^n_b(\rho, \yb)$ as follows for each $n \in [N]$ and $b \in [B] $.
\allowdisplaybreaks
\begin{subequations}\label{W_Inner}
\begin{align}
Q^n_b(\rho, \yb)&=\sup \limits_{(\pmb{\du, \e})\in \calS}  \sup_{\beta \in \mathcal{B}}  \Bigg \{ \Big ( \sum_{i \in I} d_iy_{i,b}+e_b-\mathcal{L}_b \Big )\beta_b - \rho \sum_{i \in I} |d_i-\dhat_i^n|y_{i,b}- \rho |e_b-\ehat_b^n|\Bigg\}\\
&\equiv \max \limits_{(\pmb{\du, \e})\in \calS}  \max_{\beta \in \mathcal{B}} \Bigg \{ \sum_{i \in I} \big[d_i \beta_b-\rho|d_i-\dhat_i^n| \big] y_{i,b} +  \big[ e_b \beta_b- \rho |e_b-\ehat_b^n|\big] - \mathcal{L}_b \beta_b \Bigg\}.
\end{align}
\end{subequations}

It is easy to verify that for fixed  $\yb \in \calY$, $\du \in [\dL, \dU]$, $\e \in [\eL,\eU]$, function  $\max_{\beta \in \mathcal{B}} \{ \big (\sum_{i \in I} d_iy_{i,b}+e_b-\mathcal{L}_b \big )\beta_b \} $ in \eqref{W_Inner} is convex in variable $\beta_b$. Hence, problem $\max_{\beta_b \in \mathcal{B}} \{ \cdot \} $ is a convex maximization problem. It follows from the fundamental convex analysis  that there exists an optimal  solution $\beta_b^*$ to problem  $\max_{\beta \in \mathcal{B}} \{ \cdot \} $ at one of the extreme points $\hat{\beta}_b$ of polyhedron $ \mathcal{B}:=\{-c_b^{\g} \leq \beta_b \leq c_b^{\ov}  \} $. In any extreme point, constraint $-c_b^{\g} \leq \hat{\beta}_b \leq c_b^{\ov} $ is binding at either the lower bound or upper bound. Thus, given that  $\du \in [\dL, \dU]$ and $\e \in [\eL,\eU]$, it is easy to verify that  the optimal objective value of problem \eqref{W_Inner} is equivalent to 
\begin{align}
\eta_b^{n*}= \max \left\{\begin{array}{l}  \sum\limits_{i \in I} \big[\dU_i c_b^{\ov}-\rho|\dU_i-\dhat_i^n| \big] y_{i,b} +  \big[ \eU_b c_b^{\ov}- \rho |\eU_b-\ehat_b^n|\big]  - \mathcal{L}_b  c_b^{\ov}  \\
\sum\limits_{i \in I}  \big[\dU_i c_b^{\ov}-\rho|\dU_i-\dhat_i^n| \big] y_{i,b} +  \big[ \eL_b c_b^{\ov}- \rho |\eL_b-\ehat_b^n|\big]  - \mathcal{L}_b c_b^{\ov}  \\
\sum\limits_{i \in I} \big[\dU_i  (-c_b^{\g})-\rho|\dU_i-\dhat_i^n| \big] y_{i,b} +  \big[ \eU_b  (-c_b^{\g})- \rho |\eU_b-\ehat_b^n|\big]  - \mathcal{L}_b   (-c_b^{\g}) \\
\sum\limits_{i \in I}  \big[\dU_i (-c_b^{\g})-\rho|\dU_i-\dhat_i^n| \big] y_{i,b} +  \big[ \eL_b  (-c_b^{\g})- \rho |\eL_b-\ehat_b^n|\big]  - \mathcal{L}_b (-c_b^{\g}) \\
\sum\limits_{i \in I} \big[\dL_i c_b^{\ov}-\rho|\dL_i-\dhat_i^n| \big] y_{i,b} +  \big[ \eU_b c_b^{\ov}- \rho |\eU_b-\ehat_b^n|\big]  - \mathcal{L}_b  c_b^{\ov}  \\
\sum\limits_{i \in I}  \big[\dL_i c_b^{\ov}-\rho|\dL_i-\dhat_i^n| \big] y_{i,b} +  \big[ \eL_b c_b^{\ov}- \rho |\eL_b-\ehat_b^n|\big]  - \mathcal{L}_b c_b^{\ov}  \\
\sum\limits_{i \in I} \big[\dL_i  (-c_b^{\g})-\rho|\dL_i-\dhat_i^n| \big] y_{i,b} +  \big[ \eU_b  (-c_b^{\g})- \rho |\eU_b-\ehat_b^n|\big]  - \mathcal{L}_b   (-c_b^{\g}) \\
\sum\limits_{i \in I}  \big[\dL_i (-c_b^{\g})-\rho|\dL_i-\dhat_i^n| \big] y_{i,b} +  \big[ \eL_b  (-c_b^{\g})- \rho |\eL_b-\ehat_b^n|\big]  - \mathcal{L}_b (-c_b^{\g}) \\
\sum\limits_{i \in I} \dhat_i^n c_b^{\ov} y_{i,b} +  \ehat_b^n c_b^{\ov} - \mathcal{L}_b  c_b^{\ov}  \\
\sum\limits_{i \in I} \dhat_i^n  (-c_b^{\g})y_{i,b} +  \ehat_b^n  (-c_b^{\g}) - \mathcal{L}_b   (-c_b^{\g})
 \end{array} \right\}.
\end{align}

Accordingly, we can formulate problem \eqref{W_Inner} as the following linear program for a fixed  $\yb \in \calY$
\begin{subequations}\label{W_LB_Inner}
\begin{align}
Q^n_b(\rho, \yb)= \min & \ \eta_b^n\\
& \eta_b^n \geq \sum_{i \in I} \big[\dU_i c_b^{\ov}-\rho|\dU_i-\dhat_i^n| \big] y_{i,b} +  \big[ \eU_b c_b^{\ov}- \rho |\eU_b-\ehat_b^n|\big]  - \mathcal{L}_b  c_b^{\ov},  \label{W_LB_Inner_C1}\\
& \eta_b^n \geq  \sum_{i \in I}  \big[\dU_i c_b^{\ov}-\rho|\dU_i-\dhat_i^n| \big] y_{i,b} +  \big[ \eL_b c_b^{\ov}- \rho |\eL_b-\ehat_b^n|\big]  - \mathcal{L}_b c_b^{\ov},  \label{W_LB_Inner_C2}  \\
& \eta_b^n \geq  \sum_{i \in I} \big[\dU_i  (-c_b^{\g})-\rho|\dU_i-\dhat_i^n| \big] y_{i,b} +  \big[ \eU_b  (-c_b^{\g})- \rho |\eU_b-\ehat_b^n|\big]  - \mathcal{L}_b   (-c_b^{\g}),  \label{W_LB_Inner_C3}\\
& \eta_b^n \geq  \sum_{i \in I}  \big[\dU_i (-c_b^{\g})-\rho|\dU_i-\dhat_i^n| \big] y_{i,b} +  \big[ \eL_b  (-c_b^{\g})- \rho |\eL_b-\ehat_b^n|\big]  - \mathcal{L}_b (-c_b^{\g}), \label{W_LB_Inner_C4}\\
& \eta_b^n \geq  \sum_{i \in I} \big[\dL_i c_b^{\ov}-\rho|\dL_i-\dhat_i^n| \big] y_{i,b} +  \big[ \eU_b c_b^{\ov}- \rho |\eU_b-\ehat_b^n|\big]  - \mathcal{L}_b  c_b^{\ov},  \label{W_LB_Inner_C5}\\
& \eta_b^n \geq  \sum_{i \in I}  \big[\dL_i c_b^{\ov}-\rho|\dL_i-\dhat_i^n| \big] y_{i,b} +  \big[ \eL_b c_b^{\ov}- \rho |\eL_b-\ehat_b^n|\big]  - \mathcal{L}_b c_b^{\ov},  \label{W_LB_Inner_C6}\\
& \eta_b^n \geq  \sum_{i \in I} \big[\dL_i  (-c_b^{\g})-\rho|\dL_i-\dhat_i^n| \big] y_{i,b} +  \big[ \eU_b  (-c_b^{\g})- \rho |\eU_b-\ehat_b^n|\big]  - \mathcal{L}_b   (-c_b^{\g}), \label{W_LB_Inner_C7}\\
& \eta_b^n \geq  \sum_{i \in I}  \big[\dL_i (-c_b^{\g})-\rho|\dL_i-\dhat_i^n| \big] y_{i,b} +  \big[ \eL_b  (-c_b^{\g})- \rho |\eL_b-\ehat_b^n|\big]  - \mathcal{L}_b (-c_b^{\g}), \label{W_LB_Inner_C8}\\
&\eta_b^n  \geq \sum_{i \in I} \dhat_i^n c_b^{\ov} y_{i,b} +  \ehat_b^n c_b^{\ov} - \mathcal{L}_b  c_b^{\ov},  \label{W_LB_Inner_C9}\\
& \eta_b^n \geq \sum_{i \in I} \dhat_i^n  (-c_b^{\g})y_{i,b} +  \ehat_b^n  (-c_b^{\g}) - \mathcal{L}_b   (-c_b^{\g}).
 \label{W_LB_Inner_C10}
\end{align}
\end{subequations}
Summing over $b \in [B]$ and $n \in [N]$ and combining $Q^n_b(\rho, \yb)$ in the from of \eqref{W_LB_Inner} with the outer minimization problem in \eqref{DualOfInner}, we derive the following equivalent reformulation of \eqref{DualOfInner}.
\begin{subequations}\label{DualOfInner2}
\begin{align}
\inf_{\rho \geq 0, \pmb{\eta}} & \ \   \Big\{ \epsilon^p \rho + \frac{1}{N} \sum_{n=1}^N  \eta_b^n  \Big\} \label{DualOfInner:Obj1_2}\\
\text{s.t.}  &  \ \  \eqref{W_LB_Inner_C1}-\eqref{W_LB_Inner_C10}. \label{DualOfInner:C2_2}
\end{align}
\end{subequations}
Combining the inner problem in the form of \eqref{DualOfInner2} with the outer minimization problem in  \eqref{DASBA_Obj}, we derive the following equivalent  mixed-integer non-linear programming (MINLP) reformulation of the DSA model in \eqref{DASBA_Obj}.
\begin{subequations}  \label{DASBA_MINLP}
\begin{align}
\hat{Z}(N, \epsilon)= \min  &  \Bigg \{ \sum_{i\in I} \sum_{b \in B\cup\{b'\}} c_{i,b} y_{i,b} +  \epsilon \rho + \frac{1}{N} \sum_{n=1}^N \eta^n_b \Bigg\} \label{DASBA_Obj_MINLP}\\
 \text{ s.t. }& \ \yb \in \calY, \ \rho \geq 0, \ \eqref{W_LB_Inner_C1}-\eqref{W_LB_Inner_C8}.
\end{align}
\end{subequations}

Note that constraints \eqref{W_LB_Inner_C1}-\eqref{W_LB_Inner_C8} contains the interaction terms $\rho y_{i,b}$ with binary variables $y_{i,b}$  and the non-negative continuous variable $\rho$, for all $i\in[I]$ and $b \in [B]$. To linearize, we define variables $\pi_{i,b}=\rho y_{i,b}$, for all $i\in[I]$ and $b \in [B]$, and introduce the following inequalities for variables $\pi_{i,b}$
\begin{subequations}
\begin{align}
&\pi_{i,b} \geq 0, \quad  \ \quad \pi_{i,b} \geq \bar \rho (y_{i,b}-1)+\rho, \label{Mac1}\\
&\pi_{i,b} \leq \bar \rho y_{i,b}, \quad \pi_{i,b} \leq \rho.  \label{Mac2}
\end{align}
\end{subequations}

Accordingly, formulation \eqref{DASBA_MINLP} is equivalent to the following MILP:
\begin{subequations}  \label{DASBA_MILP}
\begin{align}
\hat{Z}(N, \epsilon)= \min  &  \Bigg \{ \sum_{i\in I} \sum_{b \in B\cup\{b'\}} c_{i,b} y_{i,b} +  \epsilon \rho + \frac{1}{N} \sum_{n=1}^N \eta^n_b \Bigg\} \label{DASBA_Obj_MILP}\\
 \text{ s.t. }& \ \yb \in \calY, \ \rho \geq 0, \ \eqref{W_LB_Inner_C1}-\eqref{W_LB_Inner_C8}, \eqref{Mac1}-\eqref{Mac2}.
\end{align}
\end{subequations}

\section{Extension: Block Allocation and Surgery Scheduling in Flexible ORs}\label{Appx:Block}
\color{black}
\noindent In this section, we extend our W-DRO model to determine which OR to open and then assign surgery to open OR. Specifically, we consider a Surgery Block Allocation (SBA) problem, which, as defined in \cite{wang2019distributionally}, aims to determine the number of surgical blocks (i.e., ORs) to open and the number of elective surgeries to schedule in each open block. The fixed cost of opening OR or block $b \in [B]$ is $c_b^{\f}$, which incorporates staffing and equipment costs. The objective is to minimize the total cost consisting of the fixed cost of 
opening ORs or surgical blocks,  fixed cost related to scheduling or rejecting elective surgery, and the expected cost associated with OR overtime and idle time. Next, we propose a distributionally robust SBA (DSBA) model with 1-Wasserstein ambiguity (denoted as W-DSBA) for this problem.

Let us introduce additional variables defining our W-DSBA model.  For all $b \in [B]$, we define a binary decision variable $x_b$, which equal 1 if block $b$ is open and 0 otherwise.  We define feasible region $\calY'$ of $(\xb, \yb)$ in \eqref{eq:RegionY2}  such that each elective surgery is assigned to one open surgical block or postponed to the next planning period. 
\color{black}
\begin{align}
\calY'&:=\left\{(\xb, \yb):  \begin{array}{l} \ \ \ y_{i,b} \leq x_b, \ \forall b \in [B],  \\ \sum \limits_{b \in [B] \cup \lbrace b^\prime \rbrace} y_{i,b} =1, \  \forall i \in [I], \\
 \ x_b \in \{0, 1\}, \ \forall  b \in  [B], \\
  \ y_{i,b} \in \lbrace 0, 1 \rbrace, \  \forall i \in [I], b \in  [B] \cup \lbrace b^\prime \rbrace \end{array} \right\}. \label{eq:RegionY2}
\end{align} 

 Using the ambiguity set $\calF_\xi(\hat{\mathbb{P}}_\xi^N, \epsilon)$ defined in \eqref{W_ambiguity}, we formulate the following DRO model for the DSBA problem.
\begin{subequations}  \label{DASBA2}
\begin{align}
(\text{W-DSBA}) \ \hat{Z}'(N, \epsilon)&= \min \limits_{(\xb, \yb) \in \calY'} \Bigg \{ \sum_{b \in B} c_b^{\f}x_b+ \sum_{i\in I} \sum_{b \in B\cup\{b'\}} c_{i,b} y_{i,b} + \sup_{\Prob_{\xi} \in \calF (\hat{\mathbb{P}}_\xi^N, \epsilon) }\mathbb{E}_{\mathbb{P}_\xi} [f(\xb, \yb,\xi)] \Bigg \}, \label{DASBA_Obj2}
\end{align}
\end{subequations}
where for a feasible $(\xb, \yb) \in \calY'$ and a realization of $\xi$:
\begin{subequations}\label{2ndstage2}
\begin{align}
  f(\xb, \yb,\xi):=   \min_{\ob, \gb} &\sum_{b \in B} \big( c_b^{\ov}o_b+ c_b^{\g}o_g \big) \label{2nd_Obj2}\\
  \text{s.t.} & \ o_b-g_b=\sum_{i \in I} d_iy_{i,b}+e_bx_b-\mathcal{L}_b x_b,&& \forall b \in [B]\label{2nd_Const12}\\ 
  			& \ (o_b,g_b) \geq 0, && \forall b \in [B]. \label{2nd_Const22}
\end{align} 
\end{subequations}

Formulation \eqref{DASBA2} finds first-stage decisions $(\xb, \yb) \in \calY'$ that minimize the sum of (1) fixed cost of opening ORs (first term); (2) fixed cost of scheduling or postponing elective surgeries (second term); and (3) maximum expected cost of OR overtime and idle time over all distributions defined in the ambiguity set $\calF (\hat{\mathbb{P}}_\xi^N, \epsilon) $. \textcolor{black}{We remark that formulation \eqref{DASBA2} is the first DRO formulation with Wasserstein ambiguity for the DSBA problem in flexible ORs. \cite{wang2019distributionally} proposed a moment-based DRO model for surgery block allocation in dedicated OR. The objective is to determine the ORs to open and assign the surgeries in a daily listing to the ORs, minimizing the weighted sum of OR opening costs and expected overtime costs. However, \cite{wang2019distributionally}'s model cannot be used for flexible ORs because it does not account for emergency surgeries.}

 Using the same reformulation techniques in Section~\ref{sec:reformualtion}, we derive the following equivalent reformulation of the W-DSBA model in \eqref{DASBA2}
\begin{subequations}  \label{DASBA_MINLP2}
\begin{align}
\hat{Z}'(N, \epsilon)= \min  &  \Bigg \{ \sum_{b \in V} c_b^{\f}x_b+\sum_{i\in I} \sum_{b \in B\cup\{b'\}} c_{i,b} y_{i,b} +  \epsilon \rho + \frac{1}{N} \sum_{n=1}^N \eta^n_b \Bigg\} \label{DASBA_Obj_MINLP2}\\
 \text{ s.t. }& \ (\xb,\yb) \in \calY', \ \rho \geq 0, \\
& \eta_b^n \geq \sum_{i \in I} \big[\dU_i c_b^{\ov}-\rho|\dU_i-\dhat_i^n| \big] y_{i,b} +  \big[ \eU_bx_b c_b^{\ov}- \rho x_b |\eU_b-\ehat_b^n|\big]  - \mathcal{L}_b x_b c_b^{\ov},  \label{W_LB_Inner_C12}\\
& \eta_b^n \geq  \sum_{i \in I}  \big[\dU_i c_b^{\ov}-\rho|\dU_i-\dhat_i^n| \big] y_{i,b} +  \big[ \eL_b x_b  c_b^{\ov}- \rho x_b  |\eL_b-\ehat_b^n|\big]  - \mathcal{L}_bx_b  c_b^{\ov},  \label{W_LB_Inner_C22}  \\
& \eta_b^n \geq  \sum_{i \in I} \big[\dU_i  (-c_b^{\g})-\rho|\dU_i-\dhat_i^n| \big] y_{i,b} +  \big[ \eU_b  x_b (-c_b^{\g})- \rho x_b  |\eU_b-\ehat_b^n|\big]  - \mathcal{L}_bx_b    (-c_b^{\g}) , \label{W_LB_Inner_C32}\\
& \eta_b^n \geq  \sum_{i \in I}  \big[\dU_i (-c_b^{\g})-\rho|\dU_i-\dhat_i^n| \big] y_{i,b} +  \big[ \eL_b x_b  (-c_b^{\g})- \rho x_b |\eL_b-\ehat_b^n|\big]  - \mathcal{L}_b x_b  (-c_b^{\g}), \label{W_LB_Inner_C42}\\
& \eta_b^n \geq  \sum_{i \in I} \big[\dL_i c_b^{\ov}-\rho|\dL_i-\dhat_i^n| \big] y_{i,b} +  \big[ \eU_b x_b c_b^{\ov}- \rho x_b  |\eU_b-\ehat_b^n|\big]  - \mathcal{L}_b  x_b c_b^{\ov} , \label{W_LB_Inner_C52}\\
& \eta_b^n \geq  \sum_{i \in I}  \big[\dL_i c_b^{\ov}-\rho|\dL_i-\dhat_i^n| \big] y_{i,b} +  \big[ \eL_b x_b c_b^{\ov}- \rho x_b  |\eL_b-\ehat_b^n|\big]  - \mathcal{L}_b x_b c_b^{\ov},  \label{W_LB_Inner_C62}\\
& \eta_b^n \geq  \sum_{i \in I} \big[\dL_i  (-c_b^{\g})-\rho|\dL_i-\dhat_i^n| \big] y_{i,b} +  \big[ \eU_b x_b   (-c_b^{\g})- \rho x_b  |\eU_b-\ehat_b^n|\big]  - \mathcal{L}_b  x_b  (-c_b^{\g}), \label{W_LB_Inner_C72}\\
& \eta_b^n \geq  \sum_{i \in I}  \big[\dL_i (-c_b^{\g})-\rho|\dL_i-\dhat_i^n| \big] y_{i,b} +  \big[ \eL_b x_b  (-c_b^{\g})- \rho x_b |\eL_b-\ehat_b^n|\big]  - \mathcal{L}_b  x_b (-c_b^{\g}), \label{W_LB_Inner_C82}\\
&\eta_b^n  \geq \sum_{i \in I} \dhat_i^n c_b^{\ov} y_{i,b} +  \ehat_b^n x_b  c_b^{\ov} - \mathcal{L}_b  x_b c_b^{\ov}, \label{W_LB_Inner_C92}\\
& \eta_b^n \geq \sum_{i \in I} \dhat_i^n  (-c_b^{\g})y_{i,b} +  \ehat_b^n x_b  (-c_b^{\g}) - \mathcal{L}_bx_b    (-c_b^{\g}).
 \label{W_LB_Inner_C102}
\end{align}
\end{subequations}

\noindent Note that constraints \eqref{W_LB_Inner_C12}-\eqref{W_LB_Inner_C82} contains the interaction terms $\rho y_{i,b}$ and $\rho x_b$ with binary variables $y_{i,b}$ and $x_b$  and the nonnegative continuous variable $\rho$, for all $i\in[I]$ and $b \in [B]$.  To linearize, we use variables $\pi_{i,b}=\rho y_{i,b}$, for all $i\in[I]$ and $b \in [B]$, and introduce inequalities \eqref{Mac1}--\eqref{Mac2} for variables $\pi_{i,b}$. We also define variables $\tau_b=\rho x_b$, and introduce the following inequalities for variables $\tau_b$, for all $b \in [B]$.
\begin{subequations}
\begin{align}
&\tau_b \geq 0, \quad  \ \quad \tau_b \geq \bar \rho (x_{b}-1)+\rho, \label{Mac12}\\
&\tau_b \leq \bar \rho x_{b}, \quad \tau_p \leq \rho.  \label{Mac22}
\end{align}
\end{subequations}
Accordingly,  formulation \eqref{DASBA_MINLP2} is equivalent to the following MILP:
\begin{subequations}  \label{DASBA_MILP2}
\begin{align}
\hat{Z}(N, \epsilon)= \min  &  \Bigg \{\sum_{b \in V} c_b^{\f}x_b+ \sum_{i\in I} \sum_{b \in B\cup\{b'\}} c_{i,b} y_{i,b} + \epsilon \rho + \frac{1}{N} \sum_{n=1}^N \eta^n_b \Bigg\} \label{DASBA_Obj_MILP2}\\
 \text{ s.t. }& \ (y,x) \in \calY, \ \rho \geq 0, \ \eqref{W_LB_Inner_C12}-\eqref{W_LB_Inner_C102}, \eqref{Mac1}-\eqref{Mac2}, \eqref{Mac12}-\eqref{Mac22}.
\end{align}
\end{subequations}

\section{Computational Experiments}\label{sec:computational}

\noindent In this section, we use publicly available actual surgery duration data to construct several instances of DSA and compare the performance of three proposed models: our 1-Wasserstein-based DRO model (W-DRO), moment-based DRO model (M-DRO), and a sample average approximation (SAA) of the SP model.   In M-DRO, we construct an ambiguity set  based on the mean and support of surgery duration (see \ref{Appx:MoomentModel} for the formulation). The SAA model solves model \eqref{SP} with $\mathbb{P}_\xi$ replaced by an empirical distribution based on $N$ samples of random parameters (see \ref{Appex:SAA} for the formulation). In Section \ref{subsection:Expt_Setup}, we describe the set of problem instances that we constructed and discuss \textcolor{black}{other aspects of our experimental setup}.  In Section \ref{subsection:Wass_Eps}, we examine the effect of the radius $\epsilon$ on the performance of the W-DRO model. In Section \ref{subsection:Out-of-sample}, we compare the out-of-sample performance of the models under unseen data 
from the in-sample distribution (i.e., the distribution we used in the optimization step). In Section~\ref{subsection:Out-of-sample2}, we compare the out-of-sample performance of the models \textcolor{black}under unseen data from a {distribution different than the one we used in the optimization.} In Section~\ref{sec:CPU},  we analyze the W-DRO and SAA models solution times.  In  Section~\ref{sec:Dedicate}, we compare the flexible OR and dedicated OR policies. Finally, in Section~\ref{sec:blockallocations}, we present an experiment \textcolor{black}{comparing the performance the proposed W-DASBA model for the surgical block allocation problem in Section~\ref{Appx:Block} and its SP counterpart}.

\subsection{Description of the experiments} \label{subsection:Expt_Setup}

\noindent  Our computational study is based on anonymized real-world surgery data presented by \cite{mannino2012pattern} and \cite{manninoData}. This data set involves three-years (10,390 observations) worth of daily surgery data that belong to six different surgical specialties, namely General Cardiology (CARD), Gastroenterology (GASTRO), Gynecology (GYN), Medicine (MED), Orthopedics (ORTH), and Urology (URO).  \textcolor{black}{Table~\ref{table:SurgeryStats} presents the mean and standard deviation of elective surgery duration based on surgery type. These values where computed from publicly available data that is referenced in~\cite{manninoData} and \cite{mannino2012pattern}.}

 We assume that we have 10 available ORs and $B=$32 surgical blocks (the largest benchmark surgical suite and block schedule in the literature; see \cite{min2010scheduling} and \cite{ShehPad2020}). Table~\ref{table:SchedulingProblems} presents the weekly assignments of surgical blocks to ORs. Each block $b \in [B]$ is $\mathcal{L}_b=$8-hours long, and an elective surgery can be assigned to any of the blocks allocated to the corresponding surgery type during the planning horizon.  \textcolor{black}{It is worth to mention that the distribution of surgery blocks through OR rooms and the week is often influenced by surgeons' and surgical teams' schedules, and the setup of the operating rooms, among other factors.}

 \begin{table}[t!]
  \footnotesize
 \center 
   \renewcommand{\arraystretch}{0.4}
  \caption{Statistics for surgery duration (in minutes) based on surgery type. Notation: percent\%  is the percentage of patients needing a specific type of surgery, $\mu^{\tiny d}$ and  $\mu^{\tiny l}$ are respectively the empirical mean and standard deviations of $d$.}
\begin{tabular}{llllllll}
\hline
Surgery type	& percent\% &  $\mu^{\tiny d} $	&		$\sigma^{\tiny d}$ \\
\hline	
CARD & 14.01 & 99& 53  \\						
GASTRO & 17.79 & 132 & 76  \\
GYN & 27.81 & 78 & 52 \\
MED& 4.41 & 75 & 72\\
ORTH& 17.81 & 142 & 58\\
URO& 17.98 & 72 & 38\\
\hline	
\end{tabular} 
\label{table:SurgeryStats}
\end{table}

\begin{table}[t!]
 \footnotesize
 \center 
 \footnotesize
   \renewcommand{\arraystretch}{0.4}
  \caption{Block schedule.}
\begin{tabular}{cllllllllllll}
\hline
& & \\
\textbf{OR room} & \textbf{Monday}  & \textbf{Tuesday} &\textbf{ Wednesday}& \textbf{Thursday}& \textbf{Friday} \\
& &  \\
\hline
1 & GASTRO& GASTRO& GASTRO & &  \\
2& & &  GASTRO& GASTRO&GASTRO\\
3& CARD&  & CARD & & CARD  \\
4& ORTH & ORTH & & ORTH & ORTH\\
5&  & ORTH & MED & &  \\
6 &  GYN & GYN & GYN & GYN \\
7& & GYN & GYN & GYN &   GYN \\
8 & URO & URO & & URO & URO \\
9 &  CARD & & URO & & CARD \\
10 & URO &  & ORTH & &  \\
\hline																									
\end{tabular} 
\label{table:SchedulingProblems}
\end{table}


We consider two different cost structures for the objective function:  (1) Cost1: $c^{\mbox{\tiny o}}=\$26/\text{min}$, $c^{\mbox{\tiny  g}}=c^{\mbox{\tiny o}}/1.5$, and (2) Cost2: $c^{\mbox{\tiny o}}=\$26/\text{min},$ $c^{\mbox{\tiny g}}=0$. An overtime cost of $\$26$ per minute is based on the work of \cite{min2010scheduling}, \cite{ShehPad2020}, and\cite{stodd1998operating}. We fix the ratio $c^{\mbox{\tiny o}}/c^{\mbox{\tiny  g}}$ to 1.5 as in prior surgery scheduling studies \citep{jebali2015stochastic, liu2019integrated, shehadeh2019analysis}.  

We follow the same procedure in the literature to generate patient costs $c_{i,b}$ and $c_{i,b^\prime}$ (with surgery of the same type having common values of these parameters). We compute the lower and upper bound of surgery duration from the data . \cite{lamiri2008column} and \cite{lamiri2008stochastic} assumed that the daily capacity needed for emergency cases to be exponentially distributed with a mean of 3 and 2 hours, respectively, and the mean of the elective case to be ~2 hours. This indicates that, on average, the capacity required for emergency surgery in each block is approximately equivalent to the capacity needed for 1 scheduled elective surgery. We apply this logic to construct the parameters needed for emergency surgeries of each type.

 We implemented and solved the three models using the AMPL modeling language \textcolor{black}{and use CPLEX (version 12.6.2)  as the solver with its default settings.}    We performed all experiments on a MacBook Pro with an Intel Core i7 processor, 2.6 GHz CPU, and 16 GB (2667 MHz DDR4) of memory.

\subsection{Effect of $\epsilon$ in W-DRO model} \label{subsection:Wass_Eps}

\noindent The Wasserstein ball's radius $\epsilon$ is an input parameter  to the W-DRO model, and a larger radius implies that we seek a more distributionally robust solutions. In this section, we investigate the impact of $\epsilon$ on the out-of-sample performance of the W-DRO's optimal solution  $\pmb{\yhat} (\epsilon, N)$  (i.e.,  the objective value obtained by simulating $\pmb{\yhat} (\epsilon, N)$  under unseen data), with respect to the data size $N$.

We evaluate the out-of-sample performance for $\epsilon \in \{ 0.1, \ 0.5,\ 1, \ \ldots, \ 100\}$ as follows. First, for each $\epsilon$, we randomly sample 30 data sets of size $N \in $\{5, 10, 50, 100\} from the empirical distribution of surgery duration. Second, we solve the W-DRO in \eqref{DASBA_MILP} for each of the generated data sets and each  $\epsilon$. Third, we fix the first stage variables to $\pmb{\yhat} (\epsilon, N)$, then re-optimize the second-stage of the SP using 10,000  out-of-sample (unseen) data. This is to compute the corresponding out-of-sample overtime and idle time.  For illustrative and brevity  purposes, in this experiment, we focus on an instance of $I=$60 surgeries with Cost 1. We observe similar results for other instances (see \ref{Appx:Effect} for results with $I=$80).

\begin{figure}[t!]
 \centering
  \begin{subfigure}[b]{0.5\textwidth}
          \centering
        \includegraphics[width=\textwidth]{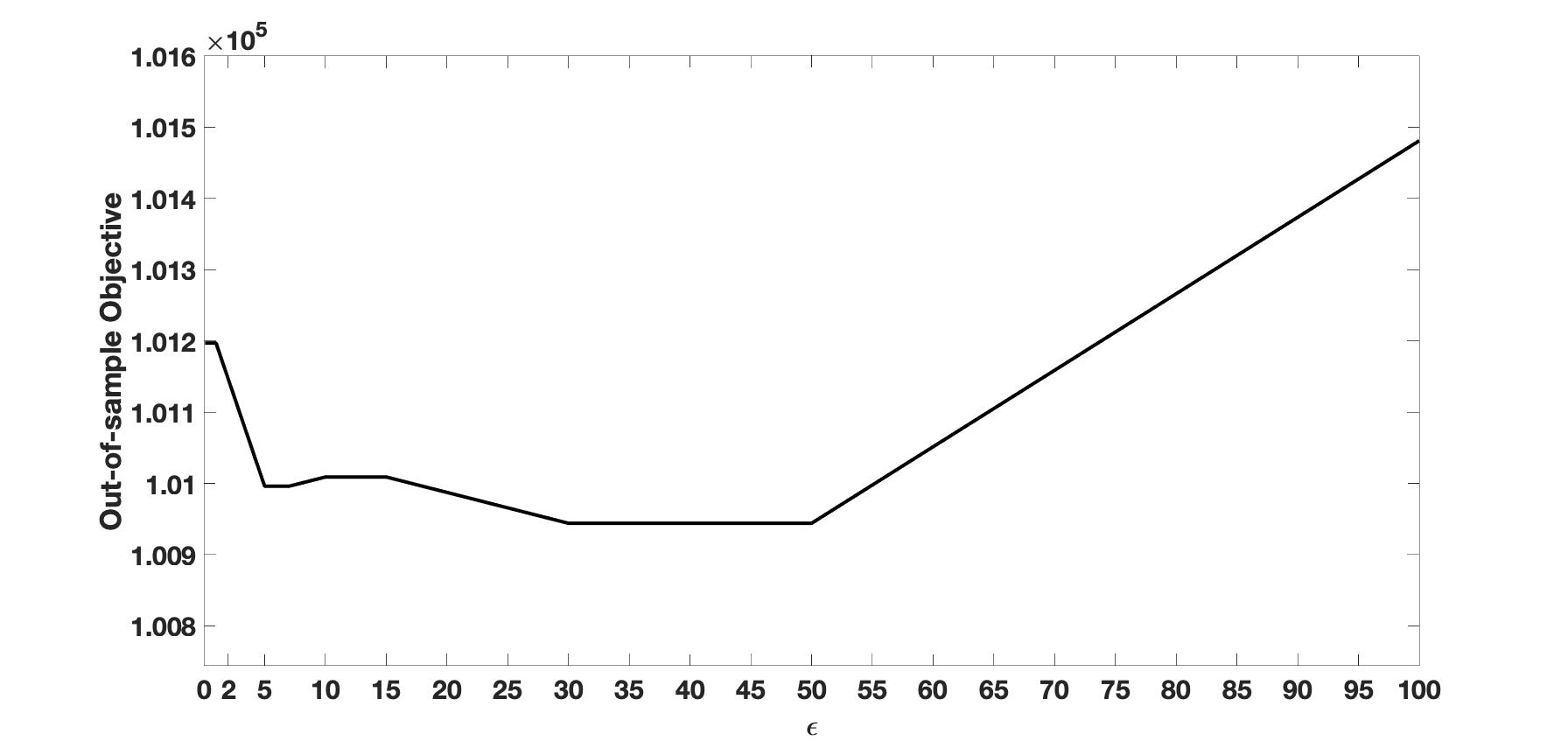}
        \caption{$N=5$}
    \end{subfigure}%
  \begin{subfigure}[b]{0.5\textwidth}
        \centering
        \includegraphics[width=\textwidth]{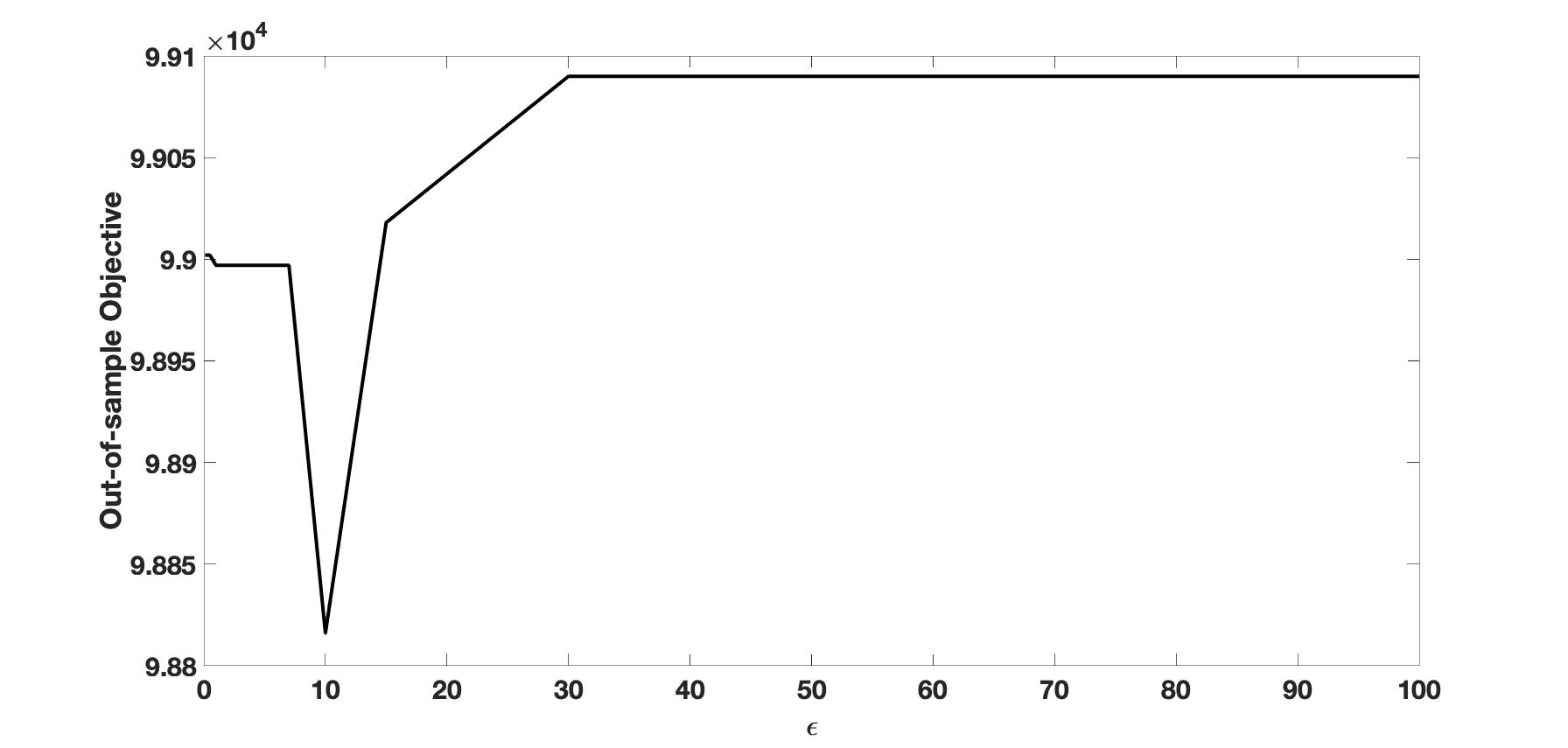}
        \caption{$N=10$}
    \end{subfigure}%
    
      \begin{subfigure}[b]{0.5\textwidth}
        \centering
        \includegraphics[width=\textwidth]{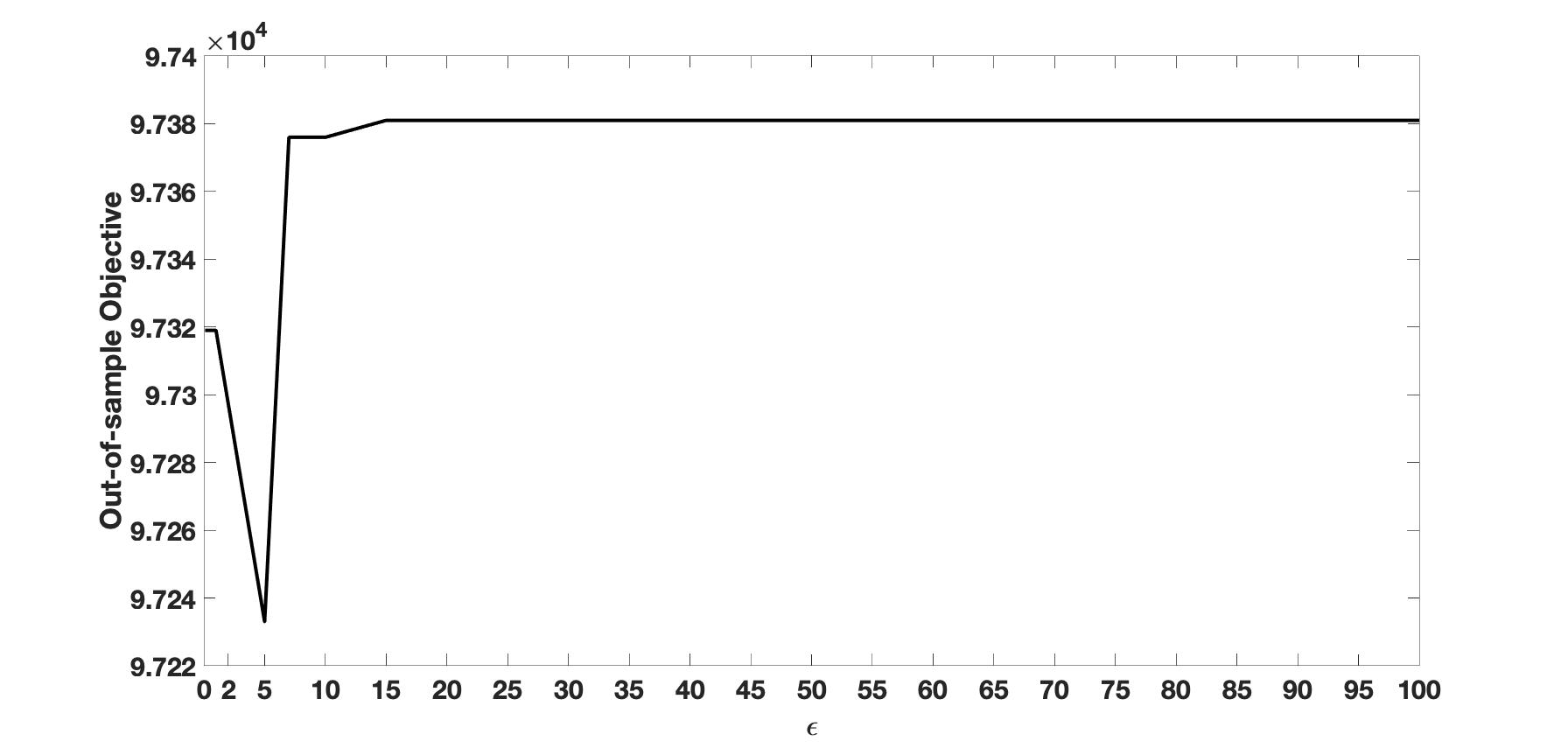}
        \caption{$N=50$}
    \end{subfigure}%
     \begin{subfigure}[b]{0.5\textwidth}
        \centering
        \includegraphics[width=\textwidth]{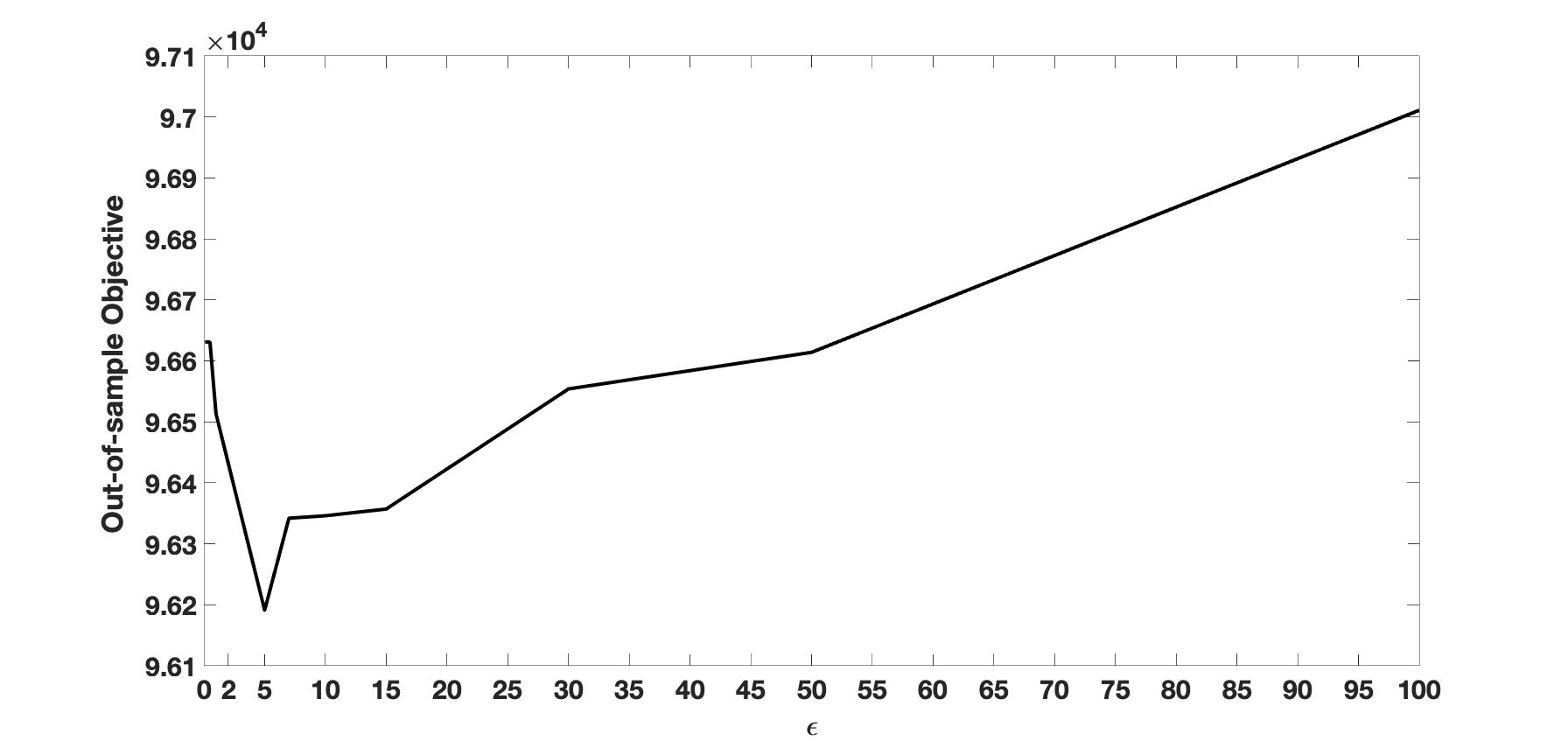}
        \caption{$N=100$}
    \end{subfigure}%
\caption{Out-of-sample performance as a function of the Wasserstein radius}\label{Fig1}
\end{figure}

Figure~\ref{Fig1} illustrates the mean values (over the 30 independent replications) of the out-of-sample cost (i.e., first-stage cost plus the out-of-sample second-stage cost)  as a function of $\epsilon$.  From this figure, we first observe that irrespective of $N$, the out-of-sample cost first improves with $\epsilon$ and then increases (or stabilize) after some (critical) value of $\epsilon$. Such pattern is often observed in the literature (see, e.g., \citealp{esfahani2018data} and \citealp{jiang2019data}), and indicate that there exists a Wasserstein radius $\epsilon^*$ such that the corresponding optimal distributionally robust solutions have the best out-of-sample performance.

Second, we observe that $\epsilon^*$ decreases with the increase in $N$.  Intuitively, a small sample does not have sufficient distributional information, and thus a larger $\epsilon$ produces distributionally robust solutions that better hedge against ambiguity. In contrast, with a larger sample, we may have more information from the data, and thus we can make a less conservative decision using a smaller $\epsilon$.

In practice, it is unlikely that OR managers and other decision-makers have the expertise or the time to conduct the above exercise (or other iterative procedures) to obtain a suitable choice of $\epsilon$ for each surgery list and parameter settings.  In addition, we often have a small data set, which is not enough for the optimization and simulation tasks (validation). Therefore, in the next sections, we derive insights with  $0.1$ (small), $10$ (relatively average), and $100$ (relatively large), which captures different extents of the robustness of the W-DRO model. This approach mimics, to a certain extent, decision-makers' thinking who do not have optimization expertise.

\subsection{Solutions quality under the in-sample distribution}\label{subsection:Out-of-sample}

\noindent In this section, we compare the out-of-sample performance of the W-DRO, M-DRO, and SAA models under unseen data sample from the in-sample data distribution (i.e., perfect information case). We test the out-of-sample performance (i.e.,  the objective value obtained by simulating the optimal solution of a model under a larger unseen data) as follows. First, we sample data sets of sizes $N\in \{5,10, \ldots, 100\}$ from the empirical distribution. Second, using each set of $N$ scenarios, we solve an instance of W-DRO, an instance M-DRO, and an instance of SAA.  Third, we fix the optimal first-stage decisions $\yb$ yielded by each model in the SP model. Then, we solve the second-stage recourse problem in \eqref{2ndstage} using $\yb$ and $N'=10, 000$ out-of-sample (unseen) data to compute the corresponding out-of-sample cost (i.e., first-stage fixed cost plus simulated second-stage cost). We perform these steps 20 times using 20 independent samples of $N$ scenarios. For the W-DRO model, we use $\epsilon=0.1$ (small), $\epsilon=10$ (relatively average), and $\epsilon=100$ (relatively large), which captures different extents of the robustness of the W-DRO model.  For illustrative purposes and brevity, we use an instance of $I=$60 and 80 surgeries under Cost1.  We observe similar results for other instances.

\begin{figure}[t!]
 \centering
  \begin{subfigure}[b]{0.5\textwidth}
          \centering
        \includegraphics[width=\textwidth]{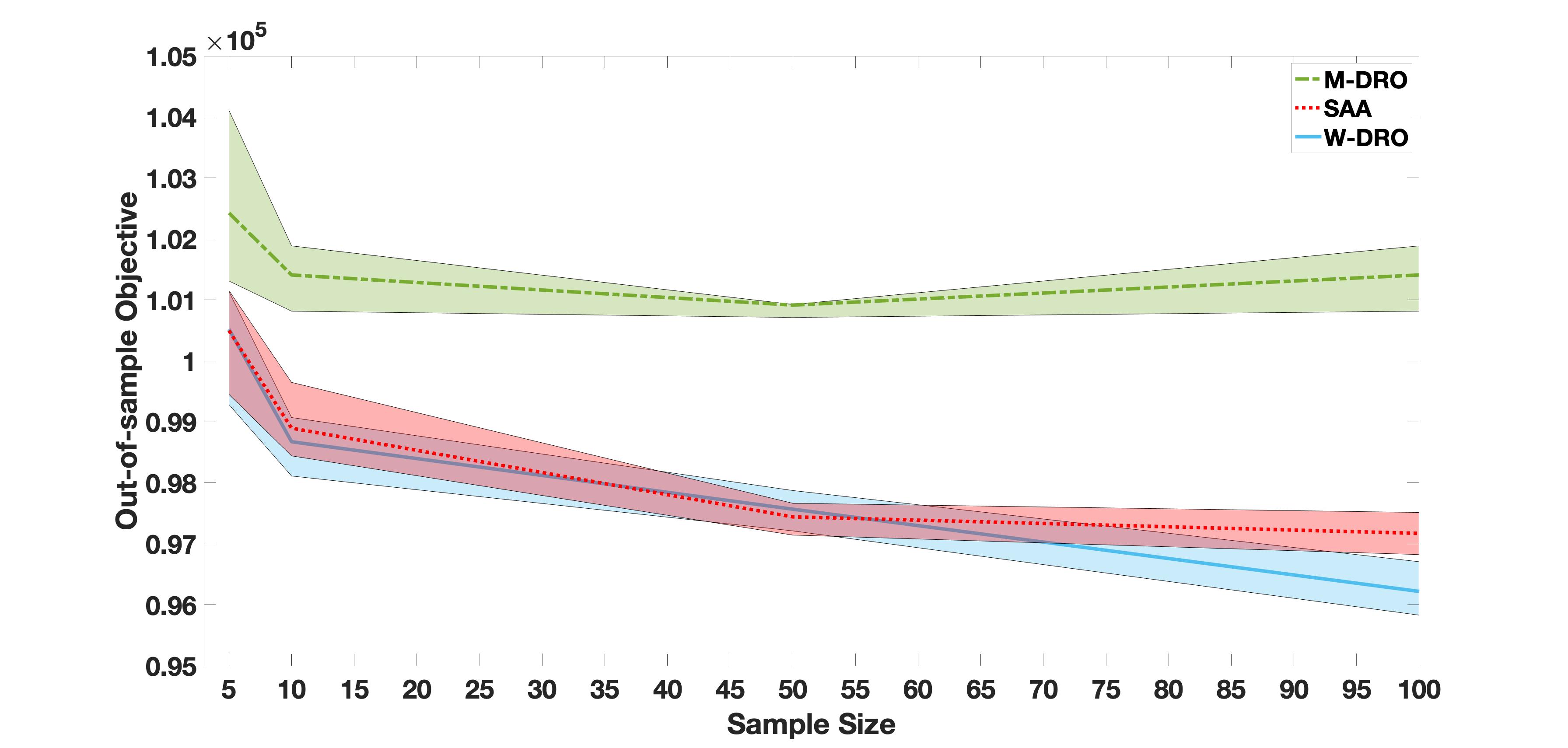}
        \caption{$\epsilon=0.10$}
    \end{subfigure}%
      \begin{subfigure}[b]{0.5\textwidth}
        \centering
        \includegraphics[width=\textwidth]{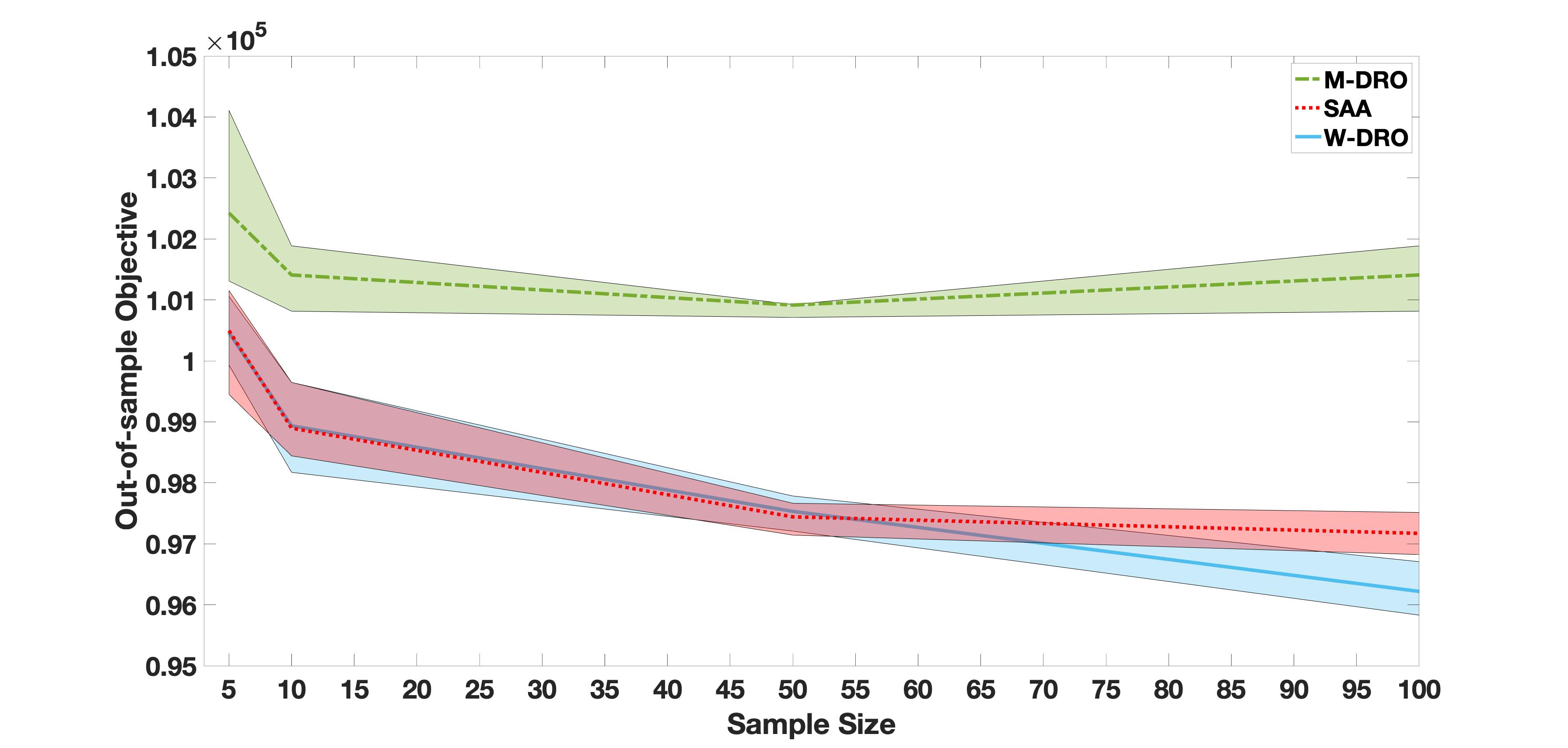}
        \caption{$\epsilon=10$}
    \end{subfigure}%
    
         \begin{subfigure}[b]{0.5\textwidth}
        \centering
        \includegraphics[width=\textwidth]{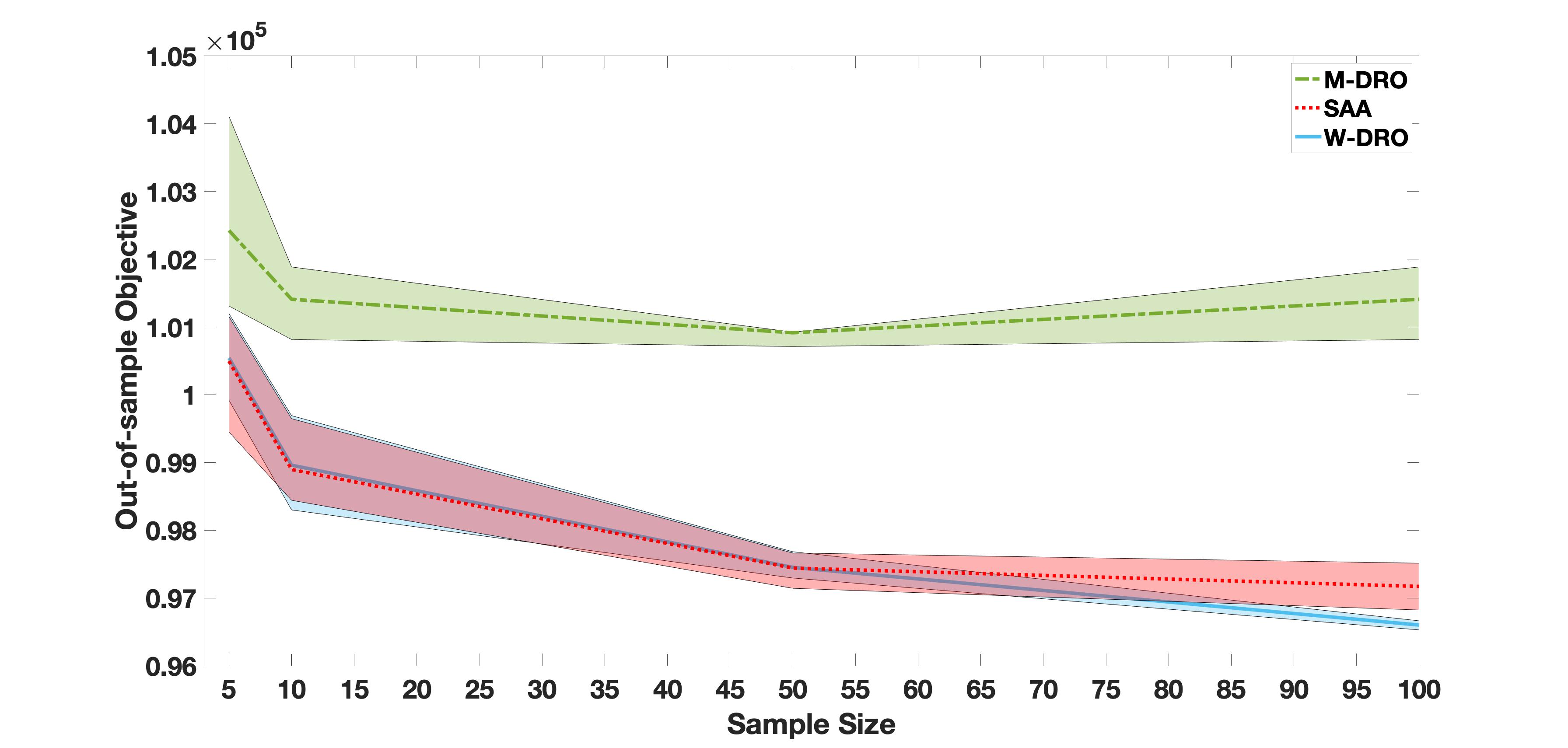}
        \caption{$\epsilon=100$}
    \end{subfigure}%
\caption{Out-of-sample performance under perfect information, $I=60$}\label{Fig2_Out_data_data60}
\end{figure}

\begin{figure}[t!]
 \centering
  \begin{subfigure}[b]{0.5\textwidth}
          \centering
        \includegraphics[width=\textwidth]{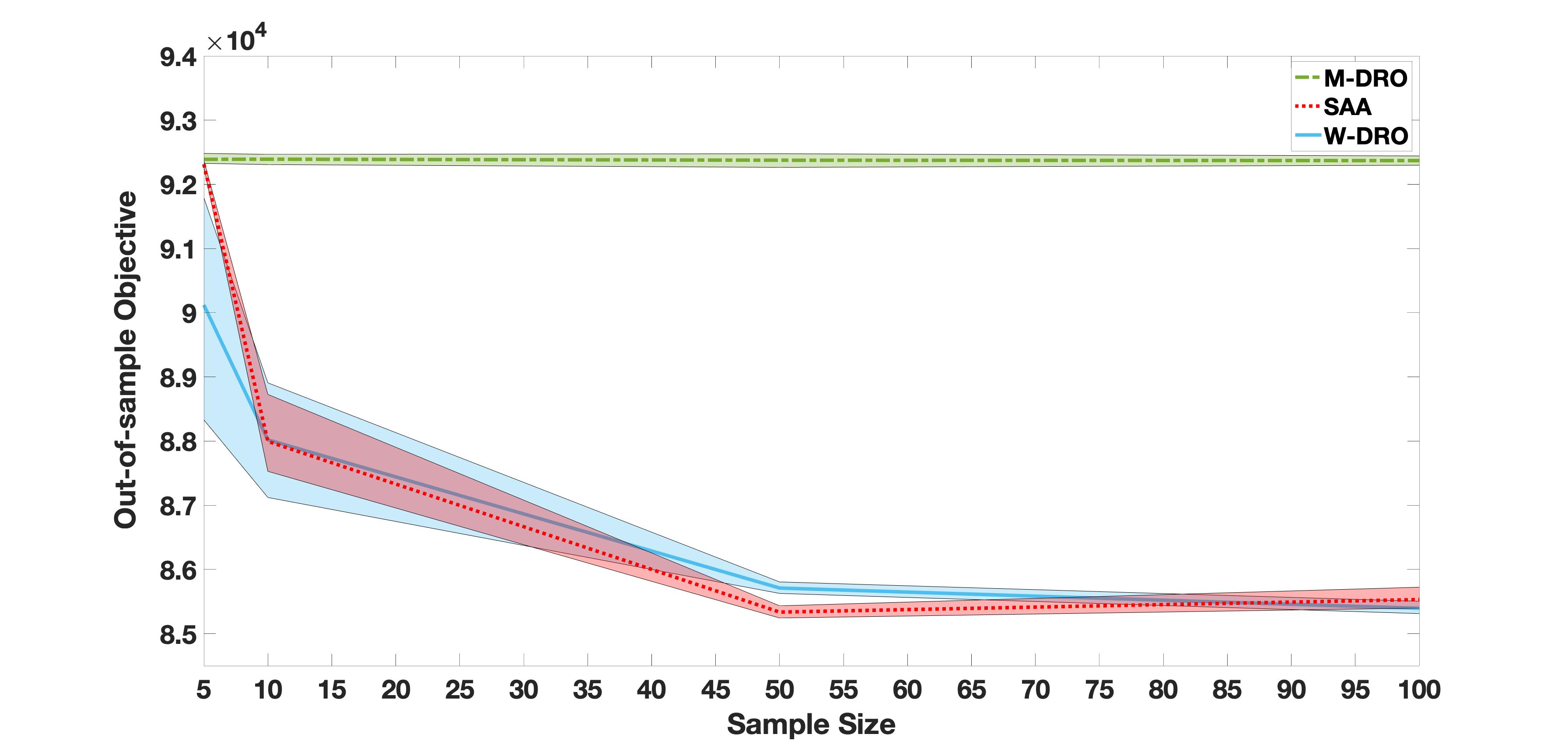}
        \caption{$\epsilon=0.10$}
    \end{subfigure}%
      \begin{subfigure}[b]{0.5\textwidth}
        \centering
        \includegraphics[width=\textwidth]{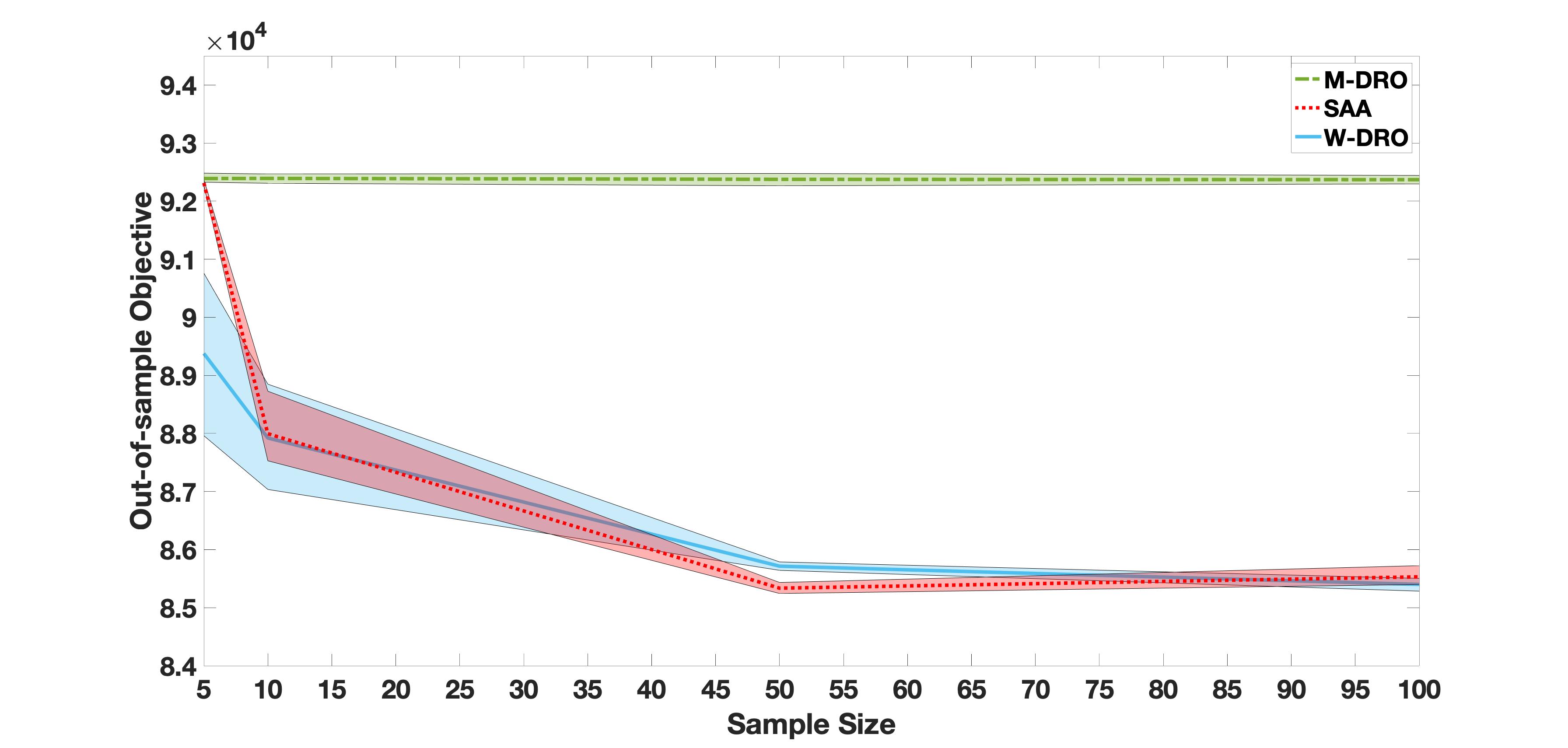}
        \caption{$\epsilon=10$}
    \end{subfigure}%
    
         \begin{subfigure}[b]{0.5\textwidth}
        \centering
        \includegraphics[width=\textwidth]{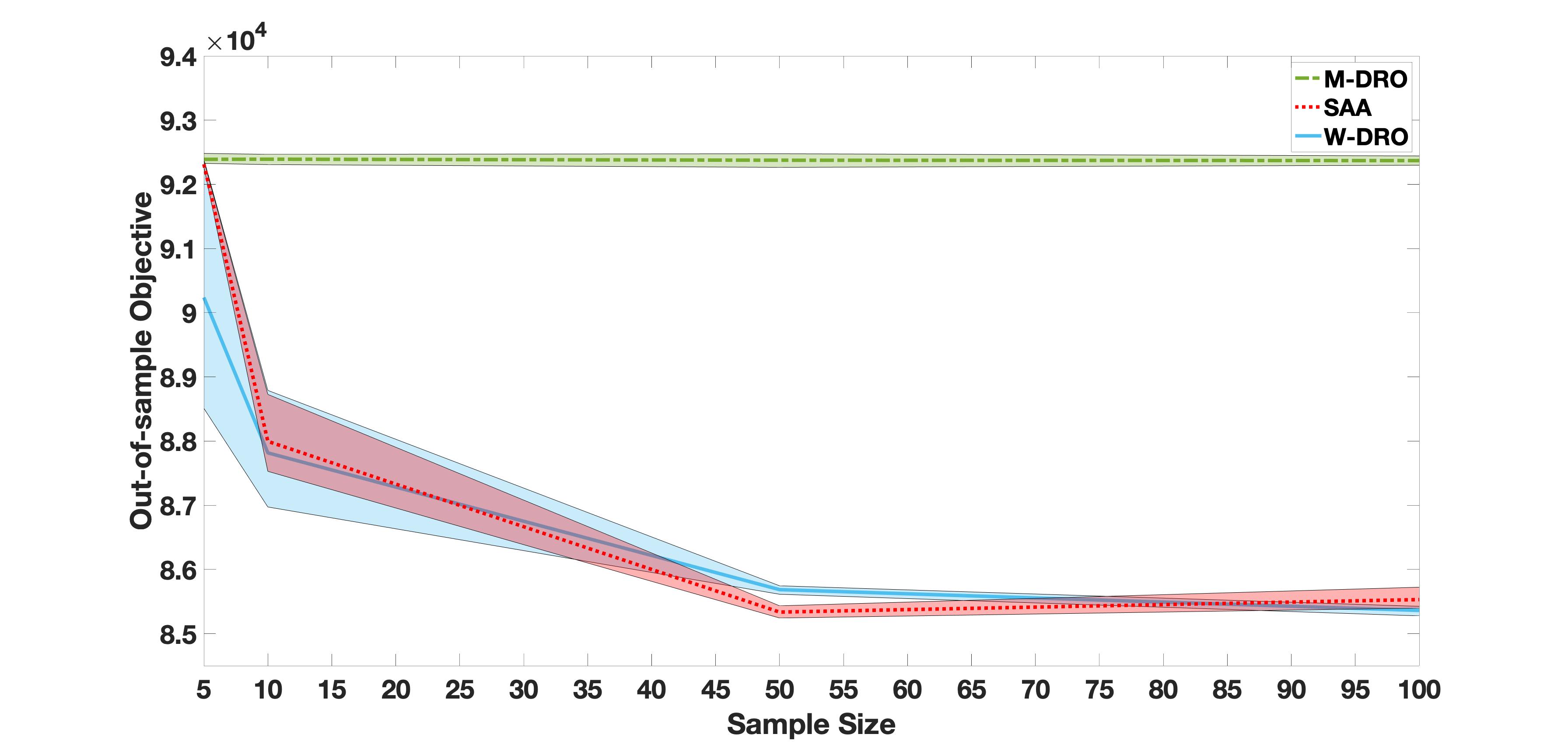}
        \caption{$\epsilon=100$}
    \end{subfigure}%
\caption{Out-of-sample performance under perfect information, $I=80$}\label{Fig2_Out_data_data80}
\end{figure}

Figures \ref{Fig2_Out_data_data60} and \ref{Fig2_Out_data_data80} present the mean (line) and the (shaded) area between the 20\% and 80\% quantiles of the out-of-sample cost as a function of $N$ and $\epsilon$ for $I=60$ and $I=80$,  respectively. We compute these values over the 20 independent simulation runs of each optimal solution and $N$. Note that we do not know the optimal value $Z^*$ of the stochastic optimization problem because we do not know the true (unknown) distribution. However, we know that $Z^*$  should be less than or equal to the bounds yielded by all models. Therefore, for the sake of illustration and convenience of making comparisons, one cane assume that $Z^*$ is the horizontal access in each sub-figure. 

We observe the following from figures~\ref{Fig2_Out_data_data60}--\ref{Fig2_Out_data_data80}. First, the out-of-sample cost of the M-DRO model are significantly larger than those of the W-DRO and SAA model across all values of  $N$ and  $\epsilon$. Second, as $N$ increases (i.e., more distributional information becomes available), the  out-of-sample costs of the W-DRO model decrease and converges to $Z^*$, which is consistent with Theorem 1 (demonstrating the that the W-DRO model enjoys the asymptotic consistency property).  In contrast, the out-of-sample cost of the M-DRO model is relatively stable. These results are consistent with the theoretical results in Theorem 1 (see \ref{Appx:Proof_Thrm1}) demonstrating the asymptotic consistency property of the W-DRO approach. The M-DRO model depends on the mean and support of surgery duration (i.e., descriptive statistics), and thus one cannot guarantee asymptotic consistency. 

Third, we observe that W-DRO slightly outperforms SAA, especially when $N$ is small. \textit{This demonstrates that the W-DRO model is capable to effectively learn distributional information even from a very small amount of data (e.g., 5 and 10 scenarios).} Thus, the proposed W-DRO model is advantageous in ORs with scarce data on surgery duration. Finally,  we observe a lower out-of-sample cost under higher $\epsilon$, which reflects the fact that one can achieve a stronger probabilistic guarantee by increasing $\epsilon$ in W-DRO.

\textcolor{black}{Finally, we attribute the higher out-of-sample cost of the M-DRO to the following. The W-DRO model schedule a larger number of elective surgeries than the M-DRO model to hedge against the cost of rejecting surgeries and excessive OR idle time. In contrast, the M-DRO model is more conservative in the sense that it attempts to hedge against observing extremely long surgery durations and overtime by scheduling fewer surgeries. However, in practice, we may not observe extremely long durations. Accordingly, the M-DRO schedule will result in a higher out-of-sample cost associated with larger idle time and rejecting a larger number of surgeries than the W-DRO schedule.} For example, consider the instance of $I=80$ surgeries with $\epsilon=0.1$, $N=50$, and Cost1. The W-DRO and M-DRO models schedule 74 ad 57 surgery, respectively. The associated average (surgery cost, overtime, idle time) with the W-DRO and M-DRO scheduling decisions are (1196, 27 minutes, 108) and (1638, 26 minutes, 201 minutes), respectively. \textcolor{black}{Note that a lower idle time indicates better utilization of the expensive OR resources. Moreover, scheduling more surgeries indicates better access to surgical care. Thus, a decision-maker may prefer the W-DRO solutions. }

\subsection{Solutions quality under missspecified distributional information}\label{subsection:Out-of-sample2}

\begin{figure}[t!]
 \centering
  \begin{subfigure}[b]{0.5\textwidth}
          \centering
        \includegraphics[width=\textwidth]{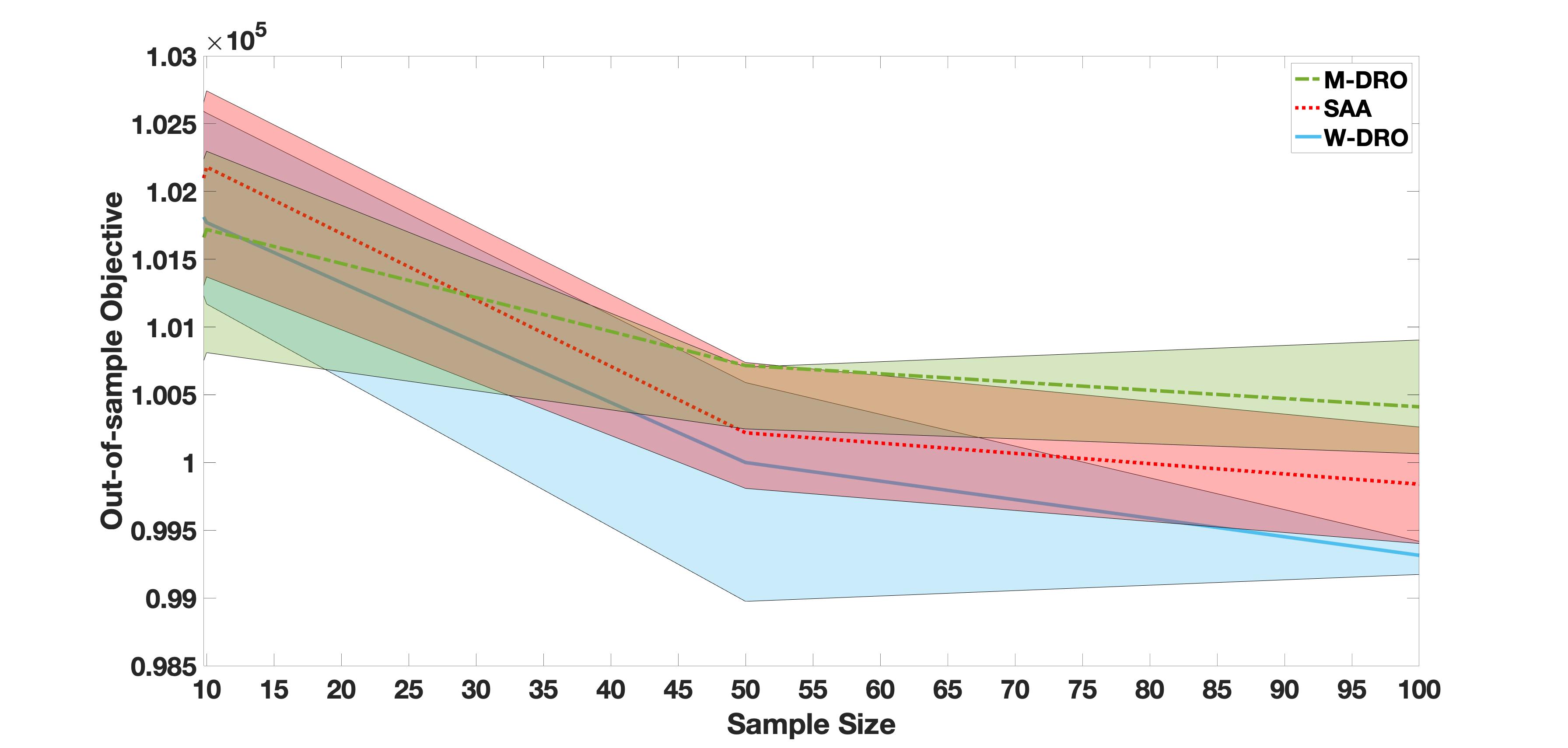}
        \caption{$\epsilon=0.10$}
    \end{subfigure}%
      \begin{subfigure}[b]{0.5\textwidth}
        \centering
        \includegraphics[width=\textwidth]{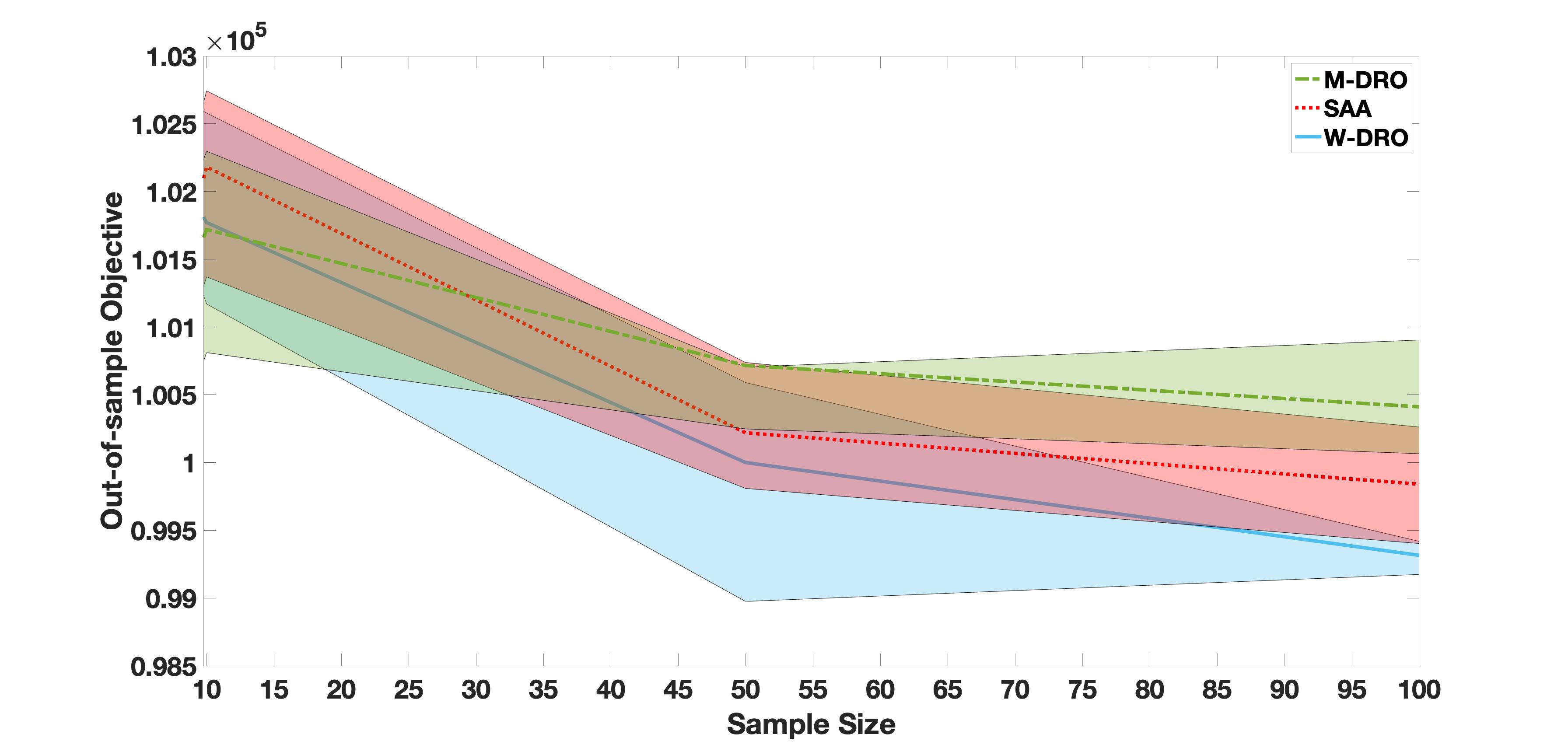}
        \caption{$\epsilon=10$}
    \end{subfigure}%
    
        \begin{subfigure}[b]{0.5\textwidth}
        \centering
        \includegraphics[width=\textwidth]{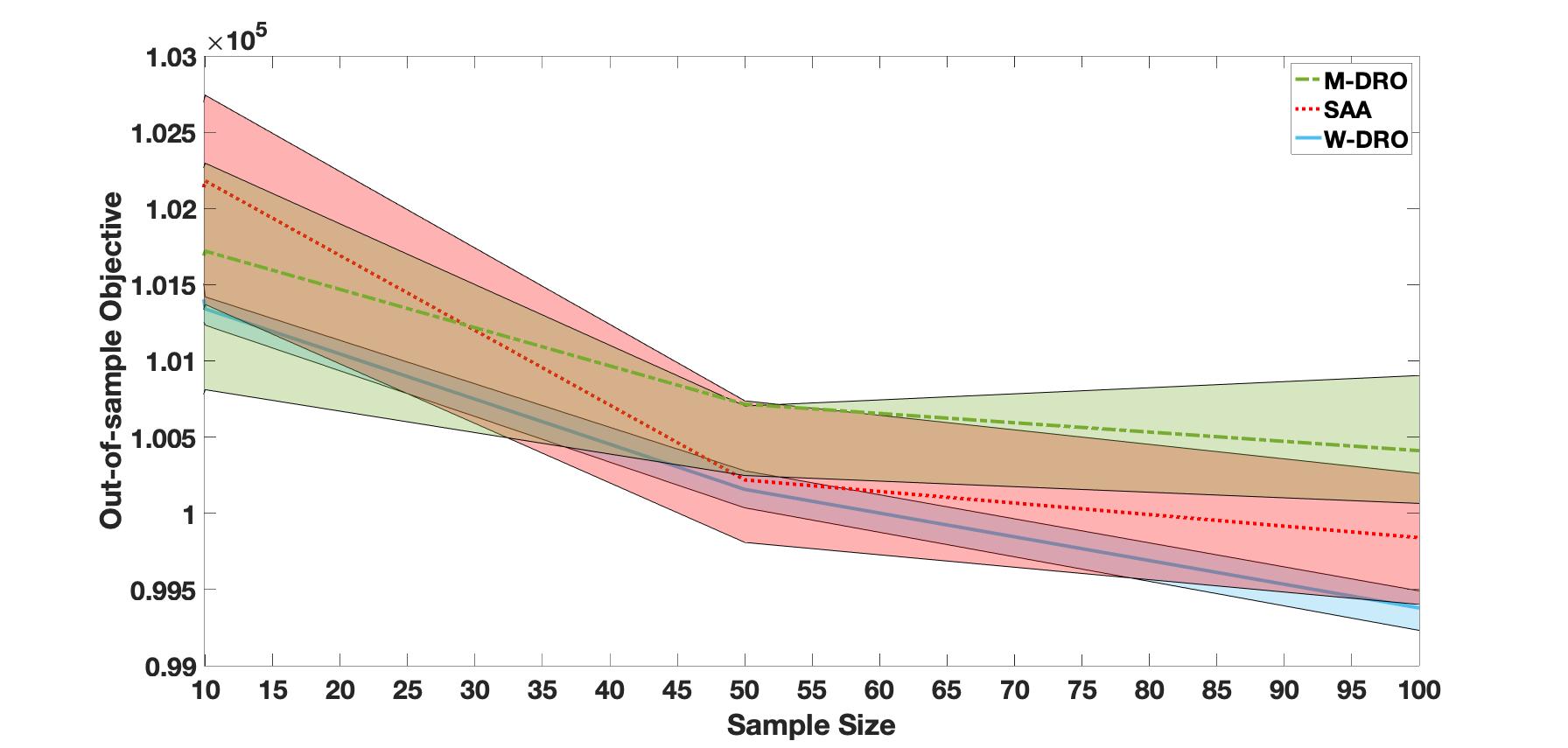}
        \caption{$\epsilon=100$}
    \end{subfigure}%
\caption{Out-of-sample performance under LogN, $I=60$}\label{Fig2_Out_Data_Log_60}
\end{figure}

\noindent \textcolor{black}{Let us now analyze the out-of-sample performance of the models under the case where the true distribution is different than the one used in the optimization.} That is, we investigate the performance of optimal solutions of the models using $N'$ data sample generated from a different distribution than the empirical distribution we used in the optimization. Specifically, we generate the $N'=10,000$ out-of-sample data from a Lognormal (LogN) distribution with the same mean, variance, and support of the in-sample distribution. For brevity, in Figures~\ref{Fig2_Out_Data_Log_60}, we present  the out-of-sample cost under logN for an instance with $I=60$.

Notably, both W-DRO and M-DRO yield lower out-of-sample costs on average and across all quantiles than the SP model when $N$ is small (i.e., limited distributional information). In fact, the M-DRO model yields a slightly lower out-of-sample cost than the W-DRO model when $N$ is small. Second, when $N$ is larger, the W-DRO and SAA models yield lower costs than the M-DRO model. Third, the W-DRO model outperforms the SAA model in most instances, especially when $N$ is small. These results indicate that the DRO approach is effective in an environment where the distribution of random parameters quickly changes, and there is small data or information on such variability. Moreover, these results emphasize the value of modeling uncertainty and distributional ambiguity of surgery duration.

\subsection{\textbf{CPU Time}}\label{sec:CPU}

\noindent In this section, we analyze the solution times of the W-DRO and SAA models. We take advantage of the fact that our problem is separable by block for some instances and thus is decomposable.

 We analyze solution times of the models (i.e., total CPU time needed to construct the OR schedule) with three sizes of surgeries list $I\in \{60, 80, 100\}$ under Cost1 ($c^{\mbox{\tiny o}}=\$26/\text{min}$, $c^{\mbox{\tiny  g}}=c^{\mbox{\tiny o}}/1.5$), and Cost2 ($c^{\mbox{\tiny o}}=\$26/\text{min},$ $c^{\mbox{\tiny g}}=0$). \textcolor{black}{Recall that the W-DRO and SAA models are scenario based. Thus, solution times of these models depend on the sample size, $N$.} Therefore, we next analyze solution times of these two models under $N\in\{5, \ 10, \ 50, \ 100, \ 500\}$. For each $N$, $I$, and cost structure, we generate 10 random instances and solve each using the W-DRO and SAA models.  Table~\ref{table:CPU_cost} presents the minimum (Min), Average (Avg), and maximum (Max) solution times (in seconds) of these instances.

\textcolor{black}{We make the following observations from  Table~\ref{table:CPU_cost}}. First, solution times increase with $I$ and $N$. Second, both models solve all instances fairly quickly under Cost2. Third, solution times are longer under Cost1, which makes sense as this cost structure require more variables (for idle time) and includes more criteria in the objective to optimize. Fourth, under Cost2, solution times of the W-DRO model are approximately similar to the SAA solution times when $N$ is small (e.g., 5 or 10) and longer than those of the SAA model when $N$ is large (e.g., 100 or 500).  However, when $N=500$, solution times of the W-DRO are all less than 5 minutes. Moreover, we have demonstrated that our W-DRO approach is effective in ORs with scares (small) data.

We attribute the difference in solution times between the SAA and W-DRO models to their respective sizes in terms of the number of variables and constraints. The W-DRO model has more variables and constraints than the SAA model. As pointed out by \citep{artigues2015mixed, catanzaro2015improved,  fortz2017compact, junger200950, keha2009mixed,  klotz2013practical, morales2016tight, shehadeh2019analysis}, an increase in model size often suggests an increase in solution time for the LP relaxation and, thus, the model's overall solution time. 

\begin{table}[t!]
 \footnotesize
 \center 
 \footnotesize
   \renewcommand{\arraystretch}{0.5}
  \caption{Solution times in seconds}
\begin{tabular}{llllllllllllllllll}
\hline
\multicolumn{16}{c}{Cost1} \\
\hline
\textbf{Model}  & \textbf{I} & $N$ & \textbf{Min} & \textbf{Avg} & \textbf{Max}&   \textbf{I} & $N$ & \textbf{Min} & \textbf{Avg} & \textbf{Max} & \textbf{I} & $N$ & \textbf{Min} & \textbf{Avg} & \textbf{Max} \\
\hline
W-DRO & 60 & 5 &  0.2&	0.4	&0.7  &  80 & 5 &1	&	1	&	2  &  100 & 5 & 0.4	&	1	&	2\\
SAA & & & 0.1&	0.1	 & 0.2 && & 0.1	&	0.1	&	0.1 & & & 0.2	&	0.3	&	0.5\\
\\
 W-DRO & & 10 & 0.4	&0.5	&0.8 &  & 10 & 1 & 2 & 3 & & 10 & 1	&	2	&	6 \\
SAA  & &     & 0.10&	0.13	&0.2 &&  &  1	&	1	&	2 &&& 1 & 1 & 2\\ 
\\ 
 W-DRO& & 50 & 1	&	3	&	9 &  & 50 & 3	&	6	&	21 & & 50& 4	&	11	&	25 \\
SAA & &  & 0.13	&	0.2	&	1  & & &  1 &	1 	&	3 & & & 1 & 1& 3\\
\\
 W-DRO & & 100 & 3	&	5	&	10 & & 100 & 6	&	10	&	29 & & 100 & 6	&	14	&	63 \\
 SAA &   &  & 0.3	&	1 &	3 & & &  1	&	2	&	3 & & &  1&	2	&	4\\
 \\
 W-DRO & & 500 & 20	&	29	&	58  &  & 500 & 89	&	164	& 267 & & 500 & 105	&	133	&	174	\\
 SAA &   & & 1 	&	2	&	3& & & 2	&	4	&	11 & & & 5	&	6	&	13\\
\hline
\multicolumn{16}{c}{Cost2} \\
\hline
\textbf{Model}  & \textbf{I} & $N$ & \textbf{Min} & \textbf{Avg} & \textbf{Max}&   \textbf{I} & $N$ & \textbf{Min} & \textbf{Avg} & \textbf{Max} & \textbf{I} & $N$ & \textbf{Min} & \textbf{Avg} & \textbf{Max} \\
\hline
W-DRO & 60 & 5 &	 0.3	&	0.4	&	1 &  80 & 5	&	0.2	&	0.3	&	0.4 & 100 & 5 &0.4	&	0.5	&	0.6\\
SAA & & &  0.1	&	0.2	&	0.2  & & &  0.1 & 0.1 & 0.1 & & &0.2	&	0.3	&	0.4 \\
\\
W-DRO & & 10 & 0.4	&	0.5	&	1 & & 10	&	0.25	&	0.30	&	0.4 && 10 & 0.4	&	0.6	&	1.1 \\
SAA & & & 0.1	&	0.2	&	0.2 & & & 0.1 & 0.2 & 0.3 & & & 0.2	&	0.3	&	1\\
\\
W-DRO & & 50 &  0.5	&	0.6	&	1 & & 50 & 0.5	&	0.6	&	0.7 & & 50 & 1	&	1.4	&	1.8\\
SAA &  & & 0.2	&	0.2	&	0.3 & & & 0.2	&	0.3	&	0.5 & & & 1 & 1 & 1\\
\\
W-DRO & & 100 & 0.6	&	0.7	&	1  & & 100 & 0.9	&	1.1	&	1.5 & & 100 & 3	&	4	&	6 \\
SAA &  & & 0.3	&	0.3	&	0.4 & & & 0.3	&	0.4	&	0.5  & & & 1 & 1 & 1\\ 
\\
W-DRO & & 500 & 1.3	&	1.5	&	2.6 & & 500 &  5	&	9&		15 &  & 500 & 16	&	27&	63\\
SAA & & & 1	&	1	&	1.1 & & & 2 & 2& 3 & &&4 & 5 & 7\\
\hline 																	
\end{tabular} 
\label{table:CPU_cost}
\end{table}


\textcolor{black}{Finally, it is worth mentioning that we were able to solve all instances using the M-DRO model fairly quickly. This makes sense, because this model is deterministic, i.e., does not depend on $N$.}
\subsection{\textcolor{black}{Flexible versus dedicated ORs}}\label{sec:Dedicate}

\noindent \textcolor{black}{As mentioned earlier, some hospitals dedicate one or more operating rooms to emergency cases. Dedicating some ORs to emergency cases may have the advantage of having some ORs readily available to handle them promptly without impacting the scheduled elective surgeries.} However, this policy may results in a low utilization of costly OR resources   \citep{xiaoreserved}. In this section, we compare the flexible and dedicated OR policies. In the latter, we solve our W-DRO model assuming that emergency surgeries have some ORs dedicated for them. We compare the optimal number of scheduled surgeries and performance with an instance of $I=80$ surgeries and with the base and twice (double) the base rate of emergency cases described in Section~\ref{subsection:Expt_Setup}.

\begin{table}[t!]
 \footnotesize
 \center 
 \footnotesize
   \renewcommand{\arraystretch}{0.5}
  \caption{Optimal number of scheduled surgery under the flexible and dedicated ORs polices.}
\begin{tabular}{lcccccccccccc}
\hline  
\multicolumn{4}{c}{\textbf{Base Emergency Rate }}  & &  & \multicolumn{4}{c}{\textbf{Double Emergency Rate }} 
\\
\hline
\textbf{Policy} & \multicolumn{1}{c}{\textbf{Cost1}} &  &  \multicolumn{1}{c}{\textbf{Cost2}}  & & &  \multicolumn{1}{c}{\textbf{Cost1}} &  &  \multicolumn{1}{c}{\textbf{Cost2}}   \\
\hline
Flexible  & 79 &  &76 & & &  76 & & 72   \\
Dedicated & 79 & &  79  & & &  79 &  & 79\\
\hline													
\end{tabular} 
\label{table:Solu_Policies}
\end{table}


Table~\ref{table:Solu_Policies} presents the optimal number of scheduled surgery yielded by the flexible and dedicated policies.  \textcolor{black}{From this table, we first observe that the W-DRO model always schedules more elective surgeries under the dedicated policy than under the flexible policy.}  In fact, under the dedicated policy, the model always schedules 79/80 surgeries irrespective of the cost structure and rate of emergency cases.  \textcolor{black}{In contrast, we schedule fewer surgeries under the flexible policy, especially with Cost2 (which ignores the idle time cost) and a higher rate of emergency cases.} These results make sense because by allocating an OR for emergency cases, we can accommodate more elective surgeries, i.e., we can use the capacity of the remaining ORs solely for elective surgeries.

\begin{table}[t!]
 \footnotesize
 \center 
 \footnotesize
   \renewcommand{\arraystretch}{0.5}
  \caption{Simulation performance of optimal schedule of the flexible and dedicated policies. Notation: OT is overtime and UT= }
\begin{tabular}{llcccclcccccccccc}
\hline
\multicolumn{9}{c}{\textbf{Base Emergency Rate }}  \\
\hline
 & \multicolumn{4}{c}{\textbf{Cost1}} &  &  \multicolumn{4}{c}{\textbf{Cost2}} \\ \cline{1-4} \cline{6-9}
Metric	&	Policy	&	OT	&	UT\%	&	&	Metric	&	Policy	&	OT	&	UT\%&	\\
\hline																
Mean 	&	Flexible	&	2.7	&	63	&	&	Mean 	&	Flexible	&	1.3	&	61	\\
	&	Dedicated	&	0.7	&	44	&	&		&	Dedicated	&	0.8	&	44	\\
							&	&								\\
75-q	&	Flexible	&	3.8	&	61	&	&	75-q	&	Flexible	&	2.0	&	61	\\
	&	Dedicated	&	1.0	&	43	&	&		&	Dedicated	&	1.4	&	44	\\
							&	&								\\
95-q	&	Flexible	&	5.9	&	59	&	&	95-q	&	Flexible	&	3.7	&	59	\\
	&	Dedicated	&	2.5	&	41	&	&		&	Dedicated	&	2.9	&	43	\\
\hline	
\multicolumn{9}{c}{\textbf{{Double Emergency Rate }}}  \\
\hline
 & \multicolumn{4}{c}{\textbf{Cost1}} &  &  \multicolumn{4}{c}{\textbf{Cost2}} \\ \cline{2-4} \cline{6-9}		
 Metric	&	Policy	&	OT	&	UT\%	&	&	Metric	&	Policy	&	OT	&	UT\%	&\\
 \hline 
Mean 	&	Flexible	&	4.6	&	69	&	&	Mean 	&	Flexible	&	1.8	&	66	\\
	&	Dedicated	&	0.7	&	44	&	&		&	Dedicated	&	0.7	&	44	\\
							&	&								\\
75-q	&	Flexible	&	5.9	&	68	&	&	75-q	&	Flexible	&	2.6	&	64	\\
	&	Dedicated	&	1.0	&	43	&	&		&	Dedicated	&	1.0	&	43	\\
							&	&								\\
95-q	&	Flexible	&	8.4	&	66	&	&	95-q	&	Flexible	&	4.3	&	62	\\
	&	Dedicated	&	2.5	&	41	&	&		&	Dedicated	&	2.5	&	41	\\
	\hline 													
\end{tabular} 
\label{table:MetricPolicies}
\end{table}

Next, we analyze the simulation performance of these optimal schedules with $N'=10,000$ scenarios from the empirical distribution. Table~\ref{table:MetricPolicies} presents the mean and quantiles of overtime (total OT hours per week) and OR utilization (UT\%=$\frac{\text{available OR time-idle time}}{\text{available OR time}}\times 100\%$, where a larger UT\% indicate less idle time and higher utilization rate). From this table, we observe that there is no clear winner. While the dedicated policy yields significantly lower overtime, it yields a significantly lower utilization than the flexible policy on average and at all quantiles. This makes sense because the flexible policy uses some OR capacity for emergency surgeries (which must be performed), thus mitigating the risk of OR idle time.  In contrast, under the dedicated policy, we perform emergency surgery in different ORs than those dedicated to elective surgery. This could thus mitigate the risk of overtime in ORs dedicated to elective surgery but potentially lead to idle time due to the unused capacity of the ORs dedicated to emergency cases.

These results show that the dedicated policy may improve access by scheduling a larger number of surgery. It may also lead to less overtime but higher idle time than the flexible policy. However, as mentioned in prior literature, dedicating (or opening an OR) for emergency surgery is costly (because the OR is an expensive resource) and may not be a feasible option (e.g., financial constraints, etc).

\subsection{Block allocation example}\label{sec:blockallocations}

\noindent  In this section, we compare the performance of our proposed W-DSBA model for the surgical block allocation (SBA) problem (presented in Section~\ref{Appx:Block}) and its SP counterpart. For simplicity, we use SP-SBA to denote the SP model.

 For illustrative purposes, we use the data related to Gastro (described in Section~\ref{subsection:Expt_Setup}) and parameter settings from prior literature. Specifically, we consider $B=10$ and $B=5$ surgical blocks or ORs, and regular time length of each block is 8 hours. By \cite{min2010scheduling} and \cite{zhang2019two}, $c_b^{\f}=800$, $c_b^{\ov}=780$, and $c_b^{\g}=0$. We set the cost of scheduling or rejecting an elective surgery as described in Section~\ref{subsection:Expt_Setup}. Finally, we consider a list of $I=25$ Gastro surgeries. 

\begin{figure}[t!]
 \centering
  \begin{subfigure}[b]{0.5\textwidth}
          \centering
        \includegraphics[width=\textwidth]{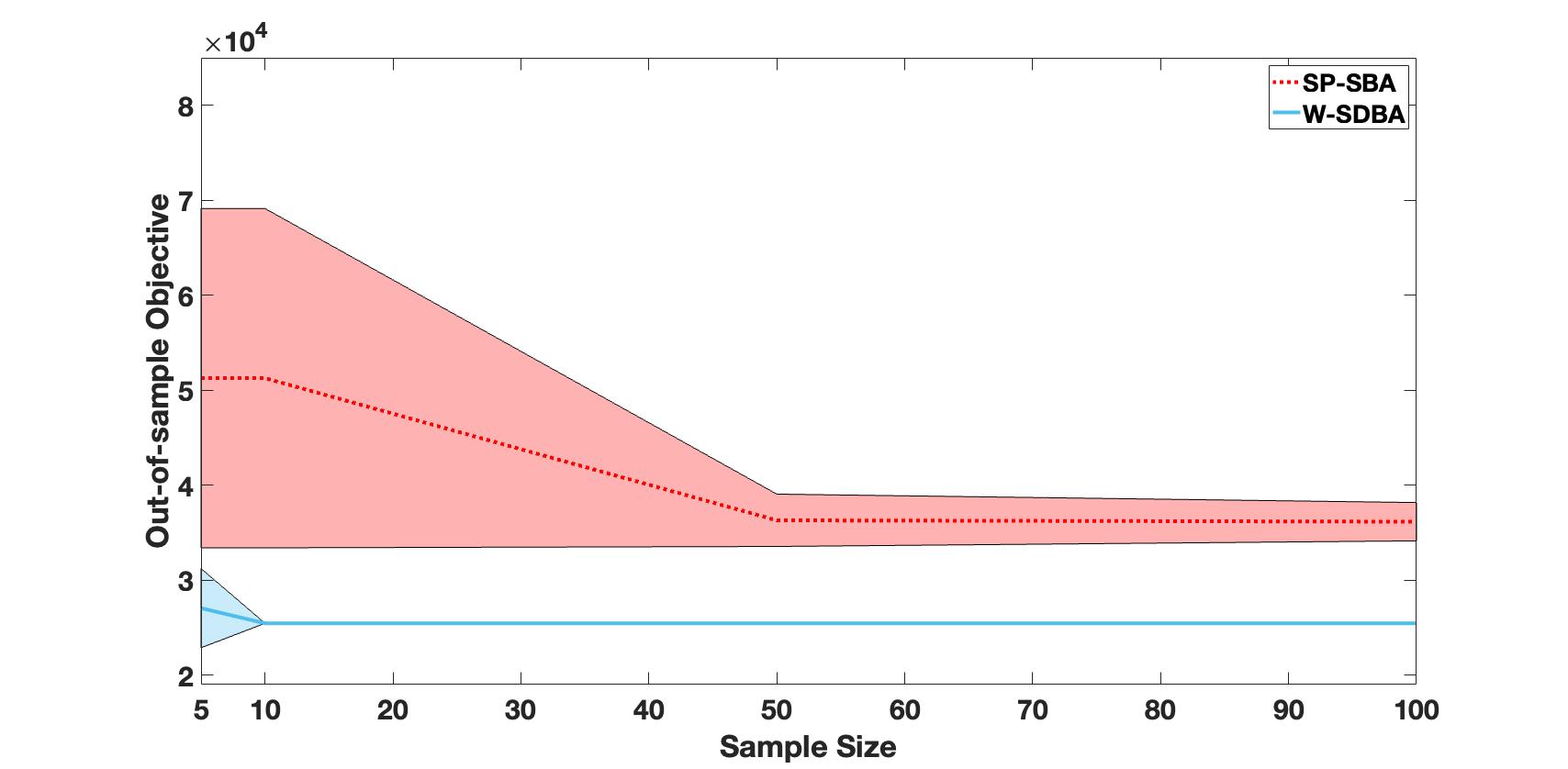}
        \caption{$B=5$}
    \end{subfigure}%
      \begin{subfigure}[b]{0.5\textwidth}
        \centering
        \includegraphics[width=\textwidth]{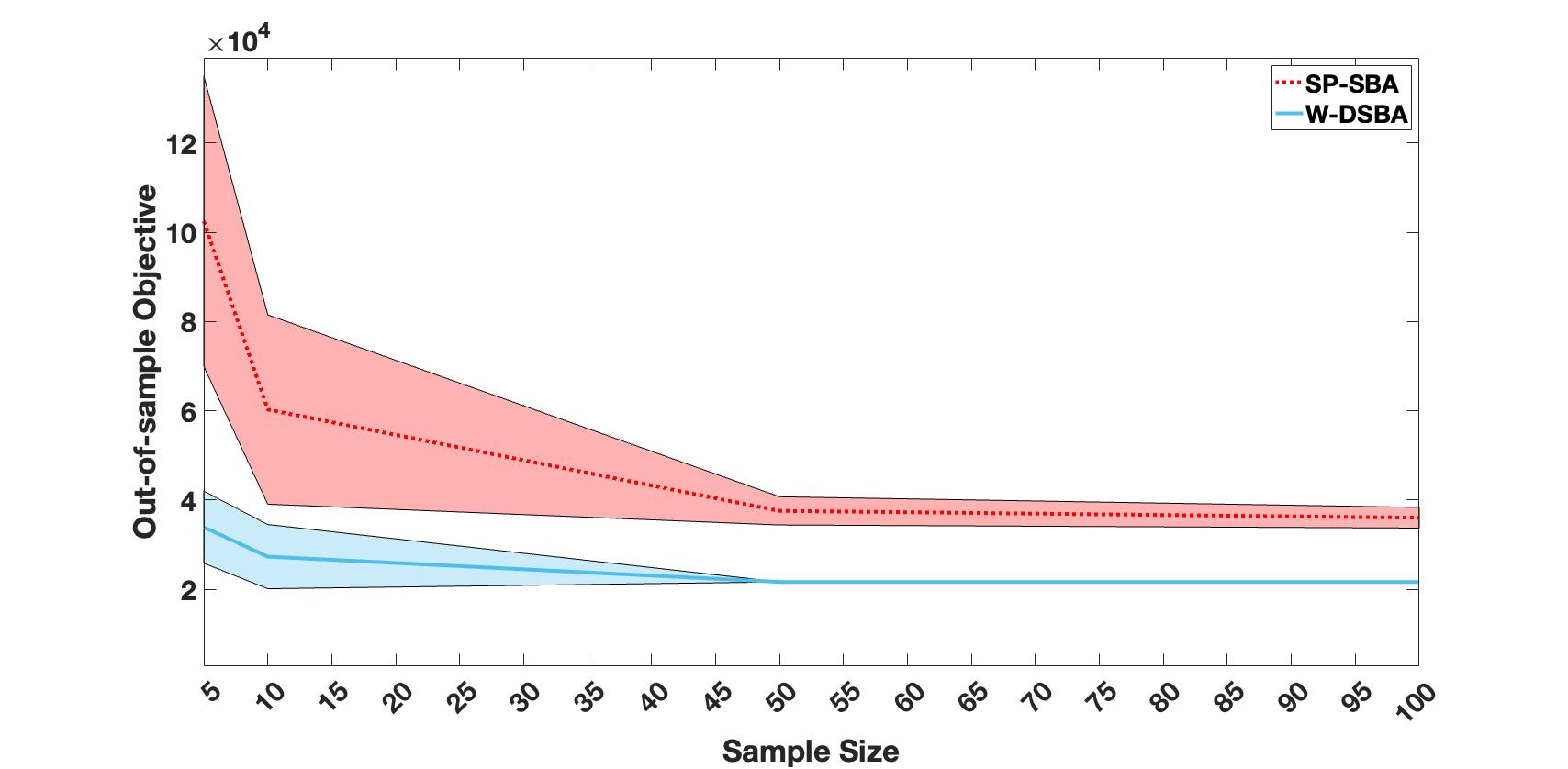}
        \caption{$B=10$}
    \end{subfigure}%
\caption{Out-of-sample performance for SBA}\label{Fig6_Out_Data_SBA}
\end{figure}

\textcolor{black}{Let us first compare the performance of these models under the perfect distributional information case described in Section~\ref{sec:computational}. That is under the case when $N$ (optimization scenarios) and $N'$ (simulation scenarios)  are sampled from the empirical distribution.} Figure~\ref{Fig6_Out_Data_SBA} presents the mean (line) and the (shaded) area between the 20\% and 80\% quantiles of the out-of-sample cost (over 30 optimization and simulation replications). The  out-of-sample cost of DRO model is computed for the best $\epsilon$ determined as described in Section~\ref{subsection:Wass_Eps}.

\textcolor{black}{We observe the following from Figure~\ref{Fig6_Out_Data_SBA}. Notably,  the W-DSBA mode yields substantially lower out-of-sample cost than the SP-SBA model under all values of $B$ and $N$. This is partly because the W-DSBA model schedules fewer surgeries and opens more blocks than the  SP-SBA model to hedge against ambiguity and overtime.} For example, when $B=10$, the W-DSBA model opens  9--10 blocks and schedules 22--25 surgeries, while the SP-SBA model opens 8 blocks and schedules 24--25 surgeries.  When $B=5$, the DRO model opens all blocks and schedules 11--15 surgeries, while the SP-SBA model opens all blocks and schedules 13--16 surgeries. In addition, the W-DSBA decisions result in less overtime, which improve surgical team satisfaction, indirectly mitigate the risk of canceling surgeries due to overtime, and reduces OR cost. For example, when $B=10$ and $N=5$, the SP and DRO models' average OR overtime is 190 and 58 minutes, respectively.  Second, we observe that the out-of-sample cost of the W-DSBA model decrease and converges as $N$ increases.  These results are consistent with the theoretical results in Theorem 1 (see \ref{Appx:Proof_Thrm1}), demonstrating the that W-DSBA  model enjoys the asymptotic consistency property.

\begin{figure}[t!]
 \centering
  \begin{subfigure}[b]{0.5\textwidth}
          \centering
        \includegraphics[width=\textwidth]{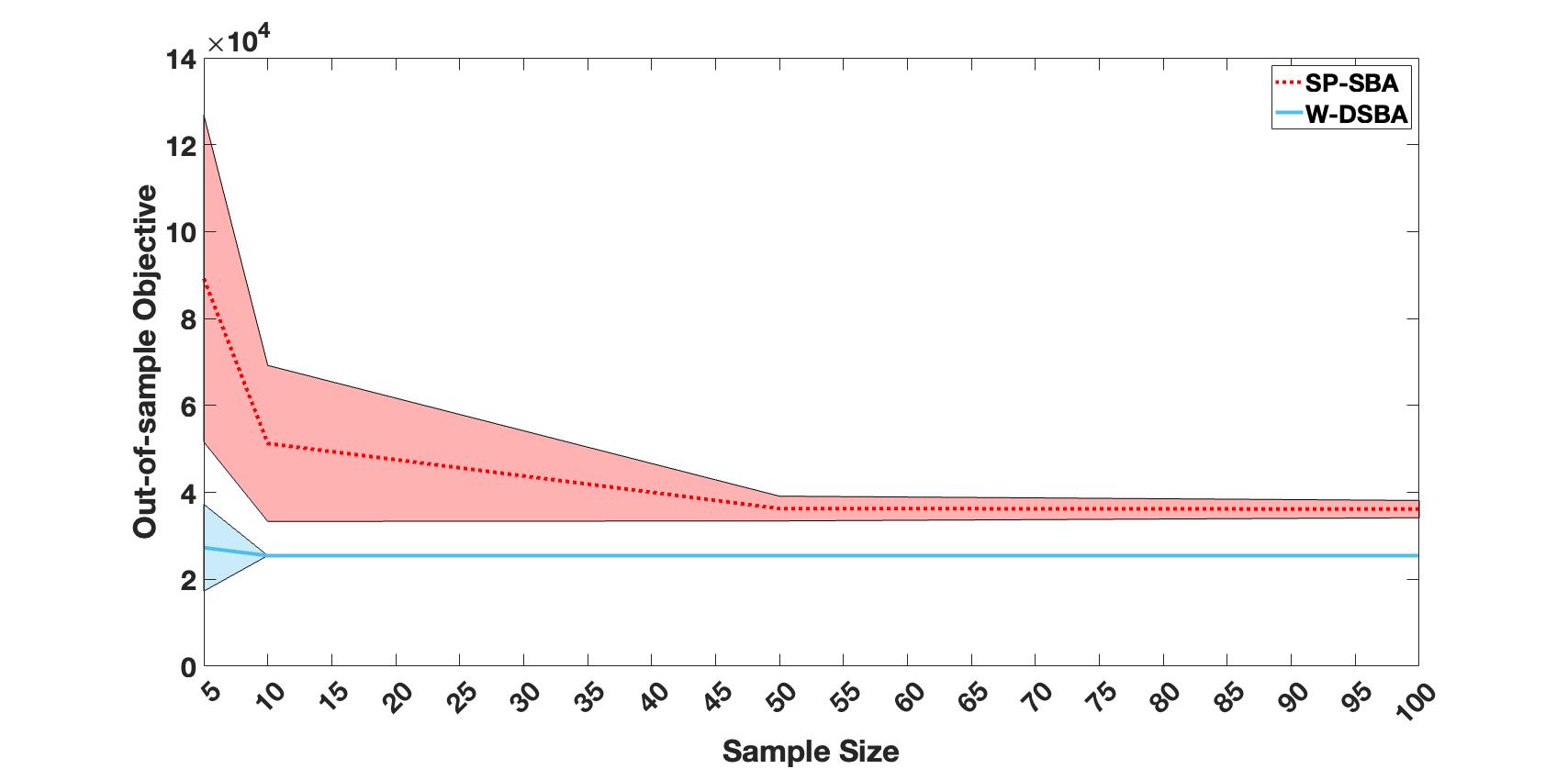}
        \caption{$B=5$}
    \end{subfigure}%
      \begin{subfigure}[b]{0.5\textwidth}
        \centering
        \includegraphics[width=\textwidth]{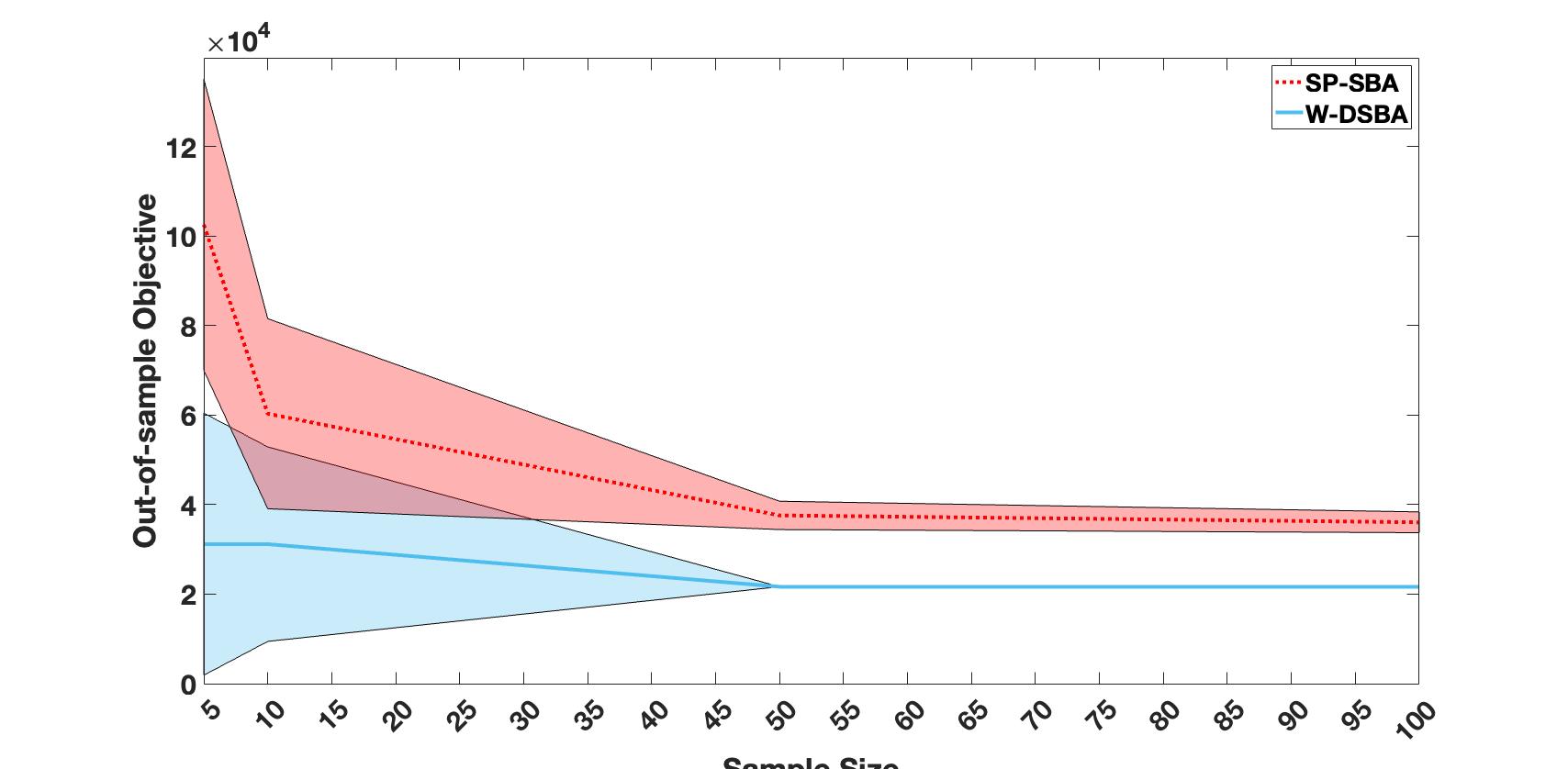}
        \caption{$B=10$}
    \end{subfigure}%
\caption{Out-of-sample performance for SBA missspecified distribution (LogN)}\label{Fig6_Out_Data_SBA_log}
\end{figure}

Next, we analyze the out-of-sample performance under the case when the data we use in the optimization is biased. Specifically, we investigate the out-of-sample performance of the optimal solution yielded by the SP-SBA and W-DSBA models using $N^\prime$ data sample from a LogN distribution defined on the same mean and support of the in-sample data.  Figure~\ref{Fig6_Out_Data_SBA_log} presents the out-of-sample cost under LogN. It is clear that the optimal solutions of the W-DSBA maintains a robust performance and yields a significantly lower out-of-sample cost than the optimal solutions of the SP-SBA model, even when small data set is used in the optimization. \textcolor{black}{These results again confirm that our DRO approach is effective in environments with limited or no data on random parameters.}


\section{Conclusion}\label{Sec:Conclusion}

\noindent In this paper, we proposed DRO models for elective surgery planning in flexible ORs under random surgery duration. Specifically, we proposed a W-DRO model that determines optimal elective surgery assigning decisions to available surgical blocks in multiple ORs over a planning horizon that minimizes the sum of patient-related costs and the worst-case expected cost associated with overtime and idle time of ORs over  all distributions residing in the 1-Wasserstein ambiguity set. We derived an equivalent MILP reformulation of our min-max W-DRO model. \textcolor{black}{In addition, we extended our W-DRO model to determine the number of ORs to open and then surgery assignments to open ORs.  We derived an equivalent MILP of this extension.}

Using real-world surgery data, we conduct extensive numerical experiments comparing our approach with the state-of-the-art approaches, namely the SP approach and a DRO approach based on mean-support ambiguity (M-DRO). Our results demonstrate that our W-DRO approach (1) enjoys both asymptotic consistency and finite-data guarantees, (2) yield robust decisions that have a good out-of-sample performance, especially under limited data and quickly varying distributions of uncertain surgery duration, and (3) is computationally efficient with solution times sufficient for real-world implementation.  \textcolor{black}{In addition, we investigated the pros and cons of the flexible and dedicated ORs policies. In particular, our results indicate that we schedule more elective surgeries under the dedicated policy. However, the dedicated policy yield longer overtime and lower utilization (when emergency cases do not materialize). In contrast, under the flexible policy, we schedule fewer elective surgeries and observe less overtime and better utilization of OR time.  Finally, we show that our proposed W-DABA model for the surgical block allocation problem has superior out-of-sample performance over the SP model for this problem.}

Overall, our results demonstrate that our proposed DRO approach can effectively learn distributional information even from a small amount of data (e.g., 5 and 10 scenarios). Thus, the proposed approach is advantageous in ORs with scarce data on surgery duration.

We suggest the following areas for future research.  First, we want to extend our approach by incorporating other multi-modal sources of uncertainty.  Second, our model can be considered as the first step toward building data-driven and robust OR and surgery planning. We aim to extend our model and build more comprehensive OR and surgery planning models, which consider all relevant organizational and technical constraints and various sources of uncertainties.  Finally, we aim to incorporate other metrics such as surgery waiting time and model the probability of canceling surgery and unexpected ORs or surgical teams' unavailability.

\vspace{4mm}
\noindent \textbf{Acknowledgments}

\vspace{1mm}

\noindent We want to thank all colleagues and practitioners who have contributed significantly to the related literature. \textcolor{black}{We are grateful to the anonymous reviewers for their insightful comments and suggestions that allowed us to improve the paper}. Dr. Karmel S. Shehadeh dedicates her effort in this paper to every little dreamer in the whole world who has a dream so big and so exciting. Believe in your dreams and do whatever it takes to achieve them--the best is yet to come for you.

\vspace{2mm}

\noindent \textbf{References}
\bibliographystyle{elsarticle-harv}
\bibliography{DSA_main_Shehadeh}
\appendix
\newpage

\section{Proof of Proposition~\ref{Prop1}}\label{Proof_Prop1}

\begin{proof} 
Recall that $\hat{\mathbb{P}}^N_\xi=\frac{1}{N}\sum_{n=1}^N \delta_{\hat{\xi}^N}$. The definition of of Wasserstein distance indicates that there exist a joint distribution $\Pi$  of $(\xi, \hat{\xi}$) such that $\mathbb{E}_{\Pi} [||\xi-\hat{\xi}||] \leq \epsilon$. In other words, for any $\mathbb{P}_{\xi}\in\mathcal{P}(\calS)$, we can rewrite any joint distribution $\Pi\in\mathcal{P}(\mathbb{P}_\xi,\hat{\mathbb{P}}^N_\xi)$ by the conditional distribution of $\xi$ given $\hat{\xi}=\hat{\xi}^n$ for $n=1,\dots,N$, denoted as $\mathbb{Q}_\xi^n$. That is, $\Pi=\frac{1}{N}\sum_{n=1}^N (\mathbb{Q}_\xi^n \times \delta_{\hat \xi}^n ) $. Notice that if we find one joint distribution $\Pi\in\mathcal{P}(\mathbb{P}_\xi,\hat{\mathbb{P}}^N_\xi)$ such that $\int ||\xi-\hat{\xi}|| \ d\Pi \leq \epsilon$, then $\text{dist}(\mathbb{P}_\xi,\hat{\mathbb{P}}_\xi^N)\leq\epsilon$. Hence, we can drop the infimum operator in Wasserstein distance and arrive at the following equivalent problem. 
\begin{subequations}\label{AppxInner}
\begin{align}
& \sup_{\mathbb{Q}_\xi^n\in\mathcal{P}(\calS), n\in[N]} \frac{1}{N}\sum_{n=1}^N\int_\calS f(\yb,\xi) \ d\mathbb{Q}_\xi^n\\
& \text{s.t.}  \ \ \ \ \ \ \ \  \frac{1}{N} \sum_{n=1}^N  \int_\calS ||\xi-\hat{\xi}^n|| \ d\mathbb{Q}_\xi^n \leq \epsilon 
\end{align}
\end{subequations}
Using a standard strong duality argument and letting $\rho\geq 0$ be the dual multiplier, we can reformulate problem \eqref{AppxInner} by its dual, i.e.
\begin{align}
&\quad\,\inf_{\rho\geq 0} \sup_{\mathbb{Q}_\xi^n\in\mathcal{P}(\calS),n\in[N]}\left\{\frac{1}{N}\sum_{n=1}^N\int_\calS f(\yb,\xi) \ d\mathbb{Q}_\xi^n+\rho\left[\epsilon-\frac{1}{N} \sum_{n=1}^N  \int_\calS ||\xi-\hat{\xi}^n|| \  d\mathbb{Q}_\xi^n \right] \right\}\nonumber\\
&=\inf_{\rho\geq 0}  \Big\{ \epsilon \rho+ \frac{1}{N}\sum_{n=1}^N \sup_{\mathbb{Q}_\xi^n\in\mathcal{P}(\calS)}\int_\calS \left[f(\yb,\xi)-\rho||\xi-\hat{\xi}^n||\right] d\mathbb{Q}_\xi^n\Big\} \nonumber\\
&=\inf_{\rho\geq 0} \Big\{ \epsilon \rho + \frac{1}{N} \sum_{n=1}^N \sup_{\xi \in \calS}  \{ f(\yb,\xi)-\rho || \xi -\hat{\xi}^n || \} \Big\}.
\end{align}

\end{proof}

\newpage

\section{Proof of Lemma 1 Adapted from \cite{jiang2019data} and Theorem 2 in \cite{fournier2015rate} }\label{Appix:Proof_Lemma1}

\noindent \textbf{Lemma 1.} (Adapted from \cite{jiang2019data} and Theorem 2 in \cite{fournier2015rate}). \textit{Suppose that Assumption A6 holds. Then there exist none-negative constants $c_1$ and $c_2$ such that, for all $N\geq 1$ and $\beta \in (0, \min \{ 1, c_1\} )$,}
\begin{equation*}
\Prob^N \big \{ \text{dist} (\Prob_\xi, \hat{\Prob}_\xi^N) \leq \epsilon_N(\beta) \} \geq 1-\beta
\end{equation*}
\textit{where $\Prob_\xi^n$ represents the product measure of $N$ copies of $\Prob_\xi$ and $\epsilon_N(\beta)= \Big [ \frac{\log(c_1\beta^{-1})}{c_2N}\Big]^{\frac{1}{\max \{ 3p, n\}}}$}
\begin{proof}
 \vspace{2mm}
\noindent For completeness, in what follows, we provide the proof of Lemma 1 of \cite{jiang2019data} as detailed in their paper, which is adapted from Theorem 2 in \cite{fournier2015rate}. 

\vspace{2mm}
First, given that the support is finite and compact (Assumption A6), then there exists $\alpha>p$ and $\gamma>0$ such that $\E_{\Prob_\xi}[\exp\{\gamma||\xi||_p^\alpha\}]<\infty$. By Theorem 2 of \cite{fournier2015rate}, for any $N \geq 1$ and $\epsilon\in(0,\infty)$, there exists positive constants $c$ and $C$ depending only on $p$, $n$, $\alpha$ and $\gamma$ such that  $$\Prob_\xi^N\left(\text{dist}(\Prob_\xi,\hat{\Prob}_\xi^N)\geq \epsilon^{1/p}\right)\leq a(N,\epsilon)\mathbbm{1}_{(\epsilon\leq1)}+b(N,\epsilon),$$ where 
\begin{equation}
a(N, \epsilon)=C\left\{\begin{array}{ll}
\exp \left\{-c N \epsilon^{2}\right\} & \text { if } p>n / 2 \\
\exp \left\{-c N(\epsilon / \log (2+1 / \epsilon))^{2}\right\} & \text { if } p=n / 2 \\
\exp \left\{-c N \epsilon^{n / p}\right\} & \text { if } p \in[1, n / 2)
\end{array}\right.
\end{equation}
\noindent and $b(N,\epsilon)=C\exp\{-cN\epsilon^{\alpha/p}\}\mathbbm{1}_{(\epsilon>1)}$, where $n$ is the dimension of the random vector $\xi$. Second, to derive bounds on $a(N,\epsilon)$  and  $b(N,\epsilon)$,  we discuss the following two cases based on the value of $\epsilon$ 
\begin{itemize}\itemsep0em
\item[Case 1.] $\epsilon\in(0,1]$.  Note that 
$$\left[\frac{\epsilon}{\log(2+1/\epsilon)}\right]^2 \geq \frac{\epsilon^3}{[\log(3)]^2}\,,$$

\noindent It follows that  $\epsilon[\log(2+1/\epsilon)]^2 \leq [\log(3)]^2$.  Thus, we can have the following upper bound on $a(N,\epsilon)$ 
 $$a(N,\epsilon)\leq C\exp\left\{ -\frac{c}{[\log(3)]^2}N\epsilon^{\max\{3,n/p\}} \right\}.$$
\item[Case 2.] $\epsilon>1$. Take $\alpha=\max\{3p,n\} > p$. It follows that  
$$b(N,\epsilon)\leq C\exp\left\{ -cN\epsilon^{\max\{3,n/p\}}\right\}.$$
\end{itemize}
\noindent Summarizing the above two cases and letting $c_1=C$ and $c_2=c/[\log(3)]^2$, we derive
$$\Prob_\xi^N\left(\text{dist}(\Prob_\xi,\hat{\Prob}_\xi^N)\geq \epsilon^{1/p}\right)\leq  c_1\exp\left\{ -c_2N\epsilon^{\max\{3,n/p\}}\right\},$$ 
\noindent By equating the right hand side of the last inequality to $\beta$ and solving for $\epsilon^{1/p}$, we obtain $\epsilon=[(c_2 N)^{-1}\log(c_1\beta^{-1})]^{-\max\{3,n/p\}^{-1}}$. Plugging this expression to the last inequality, we get the desired result. Note that in our ambiguity set $p=1$.

\end{proof}

\section{Proof of Asymptotic Consistency}\label{Appx:Proof_Thrm1}

\noindent \textbf{Theorem 1.} (Asymptotic consistency, adapted from \cite{jiang2019data} and Theorem 3.6 of \cite{esfahani2018data}). \textit{Suppose that Assumption A6 holds. Consider a sequence of confidence levels $\{ \beta_N\}_{N \in \mathbb{R}}$ such that $\sum_{N=1}^\infty \beta_N <\infty $ and $\lim_{N\rightarrow \infty} \epsilon_N(\beta_N)=0$, and let $(\hat{y} (N, \epsilon_N(\beta_N))$ represents an optimal solution to W-DRO with the ambiguity set $\calF (\hat{\Prob}_\xi^N, \epsilon_N(\beta_N))$. Then, $\Prob_\xi^\infty$- almost surely we have $\hat{Z}(N, \epsilon_N(\beta_N))\rightarrow Z^*$ as $N\rightarrow \infty$. In addition, any accumulation points of  $\{\hat{y} (N, \epsilon_N(\beta_N) \}_{N \in \mathbb{N}}$ is an optimal solution of \eqref{SP} $\Prob_\xi^\infty$- almost surely  }
\begin{proof} 
Recall the dual of $f(y, \xib)$. 
\begin{subequations}
\begin{align}
  f(y, \xib ):=  & \max_{\beta} \Bigg \{ \sum_{b \in B} \big ( \sum_{i \in I} d_iy_{i,b}+e_b-\mathcal{L}_b \big )\beta_b  \Bigg \} \nonumber \\
  & \text{s.t.}  \ \  \beta \in \mathcal{B}= \{-c_b^{\g} \leq \beta_b \leq c_b^{\ov}  \}  \nonumber
\end{align} 
\end{subequations}
First, for fixed $\yb \in \calY$, this is an LP in $\beta$. Note that $\mathcal{B}$ is bounded and $\calY$ is a bounded set, which implies that $f(y,\xib)$ is finite and bounded. It follows that we have  $|f(y,\xib)|\leq M(1+\norms{\xib})$, where $M$ is an upper bound on $f(\cdot)$.  Second, we claim that $f(y, \xib)$  is continuous on $\calY$.  For simplicity, write the objective function of the dual problem as $b(y,\xib)^\top \beta$, where the $b$-th entry of $b(y,\xib)$ is the coefficient associated to $\beta_b$. By fundamental theorem of LP, instead of maximizing over the entire $\mathcal{B}$, we can maximize over the set of finite extreme points of $\mathcal{B}$. Hence, $f(y,\xib)$ is the maximum of finite number of linear functions of $y$ and $\xib$ and the continuity follows.  Indeed, we can show that for fixed $\yb \in \calY$, $f(y,\xib)$ is Lipschitz in $\xib$. Let $(\xib^1,\xib^2)\in \calS$. Then,
    \begin{align}
        f(y, \xib^1)- f(\yb, \xib^2) &= \max_{\beta\in\mathcal{B}} b(y, \xib^1)^\top \beta- \max_{\beta\in\mathcal{B}} b(y,\xib^2)^\top \beta \leq \max_{\beta\in\mathcal{B}}[ b(y,\xib^1)-b(y,\xib^2)]^\top \beta\nonumber \\
        &\leq\max_{\beta\in\mathcal{B}} ||\beta||_q ||b(y, \xib^1)-b(y,\xib^2) ||_p,
    \end{align}  
where the last inequality follows from Holder's inequality. Notice that $\max_{\beta\in\mathcal{B}}	 ||\beta||_q $ is bounded since $\mathcal{B}$ is bounded.  Since $p\geq 1$, then for any $b \in B$, $\yb \in \calY$ and $\xib^\textsl{1}, \xib^2\in \calS$
\begin{align}
f_b(y, \xib^1)-f_b(y, \xib^2)&= \max_{\beta\in\mathcal{B}} \Bigg\{  \big ( \sum_{i \in I} d_i^1y_{i,b}+e_b^1-\mathcal{L}_b \big )\beta_b \Bigg \}-  \max_{\beta\in\mathcal{B}} \Bigg\{\big ( \sum_{i \in I} d_i^2y_{i,b}+e_b^2-\mathcal{L}_b \big )\beta_b\Bigg\}\nonumber \\
& \leq  \max_{\beta\in\mathcal{B}}  \Bigg \{\big ( \sum_{i \in I} d_i^1y_{i,b}+e_b^1-\mathcal{L}_b \big )\beta_b- \big ( \sum_{i \in I} d_i^2y_{i,b}+e_b^2-\mathcal{L}_b \big )\beta_b \Bigg\} \nonumber\\
& = \max_{\beta\in\mathcal{B}} \beta_b [ (\pmb{d}^1-\pmb{d}^2)+ (e_b^1-e_b^2)   ] \nonumber\\
& \leq  \max_{\beta\in\mathcal{B}} ||\beta||_q || \xib^1-\xib^2||_p \\
& =C|| \xib^1-\xib^2||_p \\
\end{align}
where $C= \max_{\beta\in\mathcal{B}} \beta_b$ and the last inequality follows from Holder's inequality. This verifies that $f(y, \xib)$ is Lipschitz in $\xib$. Finally, the result follows directly with $p=1$ from Theorem 3.6 of \cite{esfahani2018data}, where the required conditions are shown.
\end{proof}

\section{Proof of Finite-data guarantee}\label{Appix:Proof_Thrm2}

\noindent \textbf{Theorem 2.} (Finite-data guarantee, adapted from \cite{jiang2019data} and Theorem 3.5 in \cite{esfahani2018data}).  \textit{For any $\beta \in (0, 1)$, let $\hat{y}(N, \epsilon_N(\beta_N)),$ represent an optimal solution of W-DRO with ambiguity set $\calF (\hat{\Prob}_\xi^N, \epsilon_N(\beta_N))$. Then,} 
$$\Prob_{\xi}^N \big \{  \E_{\Prob_\mathbf{\xi}} [f (  \hat{y}(N, \epsilon_N(\beta_N),\xi )] \leq \hat{Z} (N, \epsilon_N(\beta_N))\big\} \geq 1- \beta. $$
\begin{proof}
(adapted from \cite{jiang2019data}).   By Assumption A6 on the suppprt $\calS$ and  Lemma 1  all conditions of Theorem 3.5 in \cite{esfahani2018data} are satisfied. Therefore, the conclusions of Theorem 2 hold valid.
\end{proof}
\newpage

\section{Moment Model (M-DRO)}\label{Appx:MoomentModel}

\noindent In this section, we drive a DRO model for DSA based on the known mean and support ($\calS$) of random parameters (we denote this model as M-DRO). The mean-support ambiguity set has been extensively employed in the healthcare scheduling literature for the following two primary reasons. First, from a clinical point of view and as reported in several studies, the mean $\mu$, upper and lower bounds of surgery duration are intuitive statistics that can be easily approximated  based on clinicians or OR manager experiences or from available data. Second, incorporating higher moments often undermines the DRO model's computational tractability, and therefore its applicability in practice.  \textit{We remark that there is no mean-support DRO model for the specific surgery planning in flexible ORs problem that we study in this paper (see Section~\ref{sec2:Literature}).}

\textcolor{black}{Let us introduce some additional notation defining our M-DRO model}.  First, we let $\mu^{\tiny d}$ and $\mu^{\tiny e}$ represent the mean values of $\db$ and $\eb$, respectively, and denote $\mu := \mathbb{E_P}[\xi] = [\mu^{\tiny d}, \mu^{\tiny e}]^\top$ for notational brevity. Then, we consider the following mean-support ambiguity set  $\calF(\calS, \mu) $:
\begin{align}\label{eq:ambiguity}
\calF(\calS, \mu) := \left\{ \mathbb{P} \in \calP(S)  \middle| \begin{array}{l} \int_S d\mathbb{P} = 1\\ \mathbb{E_P}[\xi] = \mu \end{array} \right\},
\end{align}

where $\calP(\calS)$ in $\calF(\calS, \mu)$ represents the set of probability distributions supported on $\calS$ with mean $\mu$. Using the ambiguity set $\calF(\calS, \mu)$, we formulate M-DRO as the following min-max problem:
\begin{align}\label{M_DSA}
&\min \limits_{\yb \in \calY} \Bigg \{ \sum_{i\in I} \sum_{b \in [B]\cup\{b'\}} c_{i,b} y_{i,b} +\sup \limits_{\mathbb{P} \in \calF(\calS,  \mu)} \mathbb{E}_\mathbb{P} [f (\yb,  \xi)] \Bigg \}.
\end{align}

\textcolor{black}{Next, we derive an equivalent reformulation of problem \eqref{M_DSA} that is solvable}. We first consider the inner maximization problem $\sup \limits_{\mathbb{P} \in \calF(\calS,  \mu)} \mathbb{E}_\mathbb{P} [f (\yb,  \xi)] $ for a fixed first-stage decision $\yb \in \calY$, where $\mathbb{P}$ is the decision variable, i.e., we are choosing the distribution that maximizes $\mathbb{E}_\mathbb{P} [f (x, \xi)]$.  For a fixed $\yb \in \calY$, we can formulate this problem as the following linear functional optimization problem.  
\begin{subequations}
\begin{align}
& \max_{\mathbb{P} \geq 0} \ \int_{S}f(\yb,\xi) \ d \mathbb{P}  \\
& \ \text{s.t.} \ \   \int_{S} d_i \ d\mathbb{P}= \mu_i^{\tiny d},  \quad \quad \forall i \in [I], \label{ConInner:duration}\\
& \ \ \ \ \  \ \ \int_{S} e_b\ d\mathbb{P}= \mu_b^{\tiny e},  \quad \quad \forall b \in [B], \label{ConInner:LOS}\\
& \ \ \ \ \  \ \ \int_{S}  d\mathbb{P}= 1. \label{ConInner:Distribution}
\end{align} \label{InnerMax}
\end{subequations}

 Letting  $ \pmb \rho=[\rho_1,\ldots, \rho_I]^\top$, $\pmb \alpha=[\alpha_1,\ldots, \alpha_B]^\top$, and $\theta$ be the dual variable associated with constraints \eqref{ConInner:duration}, \eqref{ConInner:LOS}, and \eqref{ConInner:Distribution}, respectively, we present problem \eqref{InnerMax} in its dual form:
\begin{subequations}
\begin{align}
 \min_{(\pmb \rho, \pmb \alpha), \theta \in \mathbb{R} } &  \left \{ \sum \limits_{i \in I}\mu_i^{\tiny d}\rho_i +\sum_{b \in B}\mu_b^{\tiny e}\alpha_b+ \theta \right \}\label{DualInner:Obj} \\
 \text{s.t.} & \ \ \sum_{i \in I} d_i \rho_i + \sum_{b \in B}e_b \alpha_b + \theta \geq f(\pmb y, \pmb d, \pmb e),  && \forall (\pmb d, \pmb e) \in \calS,  \label{DualInner:PrimalVariabl}
\end{align} \label{DualInnerMax}
\end{subequations}

where $\pmb \rho$, $\pmb \alpha$, and $\theta$ are unrestricted in sign, and constraint \eqref{DualInner:PrimalVariabl} is associated with the primal variable $\mathbb{P}$, and strong duality hold between \eqref{InnerMax} and \eqref{DualInnerMax} (see \cite{bertsimas2005optimal, shapiro2009lectures}).

Note that for fixed ($\pmb \rho, \pmb \alpha, \theta$), constraints \eqref{DualInner:PrimalVariabl} are  equivalent to $\theta \geq \max \limits_{(\pmb{d,e}) \in \calS } \lbrace f(\pmb y, \pmb d, \pmb e)  -  (\sum_{i=1}^P d_i \rho_i + \sum_{b \in B} e_b \alpha_b ) \rbrace $. Since we are minimizing $\theta$ in \eqref{DualInnerMax}, the dual formulation of \eqref{InnerMax} is equivalent to:
\begin{align}\label{M_DSA_Inner_full} 
 \min_{(\rho, \alpha) } & \  \Bigg\{ \sum \limits_{i \in I} \mu_i^{\tiny d}\rho_i +\sum_{b \in B}\mu_b^{\tiny e}\alpha_b+  \max \limits_{(d,e) \in \calS } \big \{ f(y,d,e) -  (\sum_{i \in I} d_i \rho_i + \sum_{b \in B}e_b \alpha_b) \big\}  \Bigg\}
\end{align}
Using the dual formulation of $f(y,d,e) $ in \eqref{Dual_2nd} and the same techniques  we used to derive \eqref{W_LB_Inner}, we derive the following equivalent reformulation of \eqref{M_DSA_Inner_full}.
\begin{subequations}\label{M_DSA_Inner_Final} 
\begin{align}
 \min_{(\pmb{\rho, \alpha}) } & \  \Bigg\{ \sum \limits_{i \in I} \mu_i^{\tiny d}\rho_i y_{i,b}+\sum_{b \in B}\mu_b^{\tiny e}\alpha_b+ \sum_{b \in B}\eta_b  \Bigg\} \\
  \text{s.t.}  & \ \eta_b  \geq \big[\sum_{i \in I} \dU_i y_{i,b}+\eU_b -\mathcal{L}_b\big] (c_b^{\ov})-(\sum_{i \in I} \dU_i \rho_i y_{i,b}+\eU_b \alpha_b), && \forall b \in [B],  \label{C1_InnerM}\\
  & \ \eta_b  \geq \big[\sum_{i \in I} \dU_i y_{i,b}+\eL_b -\mathcal{L}_b\big] (c_b^{\ov})-(\sum_{i \in I} \dU_i \rho_i y_{i,b}+ \eL_b \alpha_b) && \forall b \in [B],\\
  &  \ \eta_b  \geq \big[\sum_{i \in I} \dU_i y_{i,b}+\eU_b -\mathcal{L}_b\big] (-c_b^{\g})-(\sum_{i \in I} \dU_i \rho_i y_{i,b}+ \eU_b \alpha_b), && \forall b \in [B], \\
  & \  \eta_b  \geq \big[\sum_{i \in I} \dU_i y_{i,b}+\eL_b -\mathcal{L}_b\big] (-c_b^{\g})-(\sum_{i \in I} \dU_i \rho_i y_{i,b}+ \eL_b \alpha_b), && \forall b \in [B], \\
  & \eta_b  \geq \big[\sum_{i \in I} \dL_i y_{i,b}+\eU_b -\mathcal{L}_b\big] (c_b^{\ov})-(\sum_{i \in I} \dL_i \rho_i y_{i,b}+ \eU_b \alpha_b) , && \forall b \in [B], \\
  & \ \eta_b  \geq \big[\sum_{i \in I} \dL_i y_{i,b}+\eL_b -\mathcal{L}_b\big] (c_b^{\ov})-(\sum_{i \in I} \dL_i \rho_i y_{i,b}+ \eL_b \alpha_b) , && \forall b \in [B], \\
  &  \ \eta_b  \geq \big[\sum_{i \in I} \dL_i y_{i,b}+\eU_b -\mathcal{L}_b\big] (-c_b^{\g})-(\sum_{i \in I} \dL_i \rho_i y_{i,b}+\eU_b \alpha_b) , && \forall b \in [B], \\
  & \  \eta_b  \geq \big[\sum_{i \in I} \dL_i y_{i,b}+\eL_b -\mathcal{L}_b\big] (-c_b^{\g})-(\sum_{i \in I} \dL_i \rho_i y_{i,b}+ \eL_b \alpha_b), && \forall b \in [B]. \label{C8_InnerM}
\end{align}
\end{subequations}

 Combining the inner maximization problem in the form of \eqref{M_DSA_Inner_Final} with the outer minimization problem in \eqref{M_DSA}, we derive the following equivalent reformulation of the M-DRO model in \eqref{M_DSA} as follow 
\begin{subequations}\label{M_DSA_NLMIP}
\begin{align}
\min &  \Bigg \{ \sum_{i\in I} \sum_{b \in B\cup\{b'\}} c_{i,b} y_{i,b} +\sum \limits_{i \in I}\mu_i^{\tiny d}\rho_i y_{i,b}+\sum_{b \in B}\mu_b^{\tiny e}\alpha_b+\sum_{b \in B}\eta_b \Bigg \} \label{M_ObjDR_NLMIP}\\
\text{s.t. } & \yb \in \calY, \ \eqref{C1_InnerM}-\eqref{C8_InnerM}
\end{align}
\end{subequations} 

 Note that the objective function \eqref{M_ObjDR_NLMIP} and constraints \eqref{C1_InnerM}--\eqref{C8_InnerM} contains the interaction terms $\rho_i y_{i,b}$ with binary variables $y_{i,b}$ and continuous variables $\rho_i$. To linearize,  we define $k_{i,b}=\rho_i y_{i,b}$ and introduce the following McCormick inequalities for variables $k_{i,b}$, for all $i \in [I]$ and $b \in [B]$
\begin{subequations}
\begin{align}
k_{i,b} \geq \underline{\rho}_iy_{i,b} \qquad k_{i,b} \geq \overline{\rho}_i (y_{i,b}-1)+\rho_i \label{M_Mac1} \\
k_{i,b} \leq \overline{\rho}_i y_{i,b} \qquad k_{i,b} \leq \underline{\rho}_i(y_{i,b}-1)+\rho_i \label{M_Mac2}
\end{align}
\end{subequations}

where $\underline{\rho}_i$ and $\overline{\rho}_i$ are respectively valid lower and upper bounds on variables $\rho_i$, for all $i \in [I]$. Accordingly, formulation \eqref{M_DSA_NLMIP} is equivalent to the following MILP:
\begin{subequations}
\begin{align}
\min &  \Bigg \{ \sum_{i\in I} \sum_{b \in B\cup\{b'\}} c_{i,b} y_{i,b} +\sum \limits_{i=1}^I \mu_i^{\tiny d}k_{i,b}+\sum_{b \in B}\mu_b^{\tiny e}\alpha_b+ \sum_{b \in B}\eta_b \Bigg \} \label{M_ObjDR_MIP}\\
\text{s.t. } & \yb \in \calY, \ \eqref{C1_InnerM}-\eqref{C8_InnerM}, \ \eqref{M_Mac1}-\eqref{M_Mac2}
\end{align}
\end{subequations}

\color{black}

\newpage 

\section{The Sample Average Approximation (SAA) Formulation }\label{Appex:SAA}

\begin{subequations}\label{SAA}
\begin{align}
 \min \limits_{\yb \in \calY}  \ \ & \Bigg\{ \sum_{i\in I} \sum_{b \in B\cup\{b'\}} c_{i,b} y_{i,b}+ \sum_{n=1}^N \frac{1}{N} \Big[ \sum_{b \in B} \big( c_b^{\ov}o_b^n+ c_b^{\g}o_g^n \big) \Big] \Bigg\} \\
  \text{s.t.} & \ o_b^n-g_b^n=\sum_{i \in I} d_i^ny_{i,b}+e_b^n-\mathcal{L}_b, && \forall b \in B, \forall n \in [N],\\ 
  			& \ (o_b^n,g_b^n) \geq 0, && \forall b \in B, \forall n \in [N].
\end{align} 
\end{subequations}

\section{Effect of $\epsilon$, $I=80$ Surgeries}\label{Appx:Effect}

\begin{figure}[h!]
 \centering
  \begin{subfigure}[b]{0.5\textwidth}
          \centering
        \includegraphics[width=\textwidth]{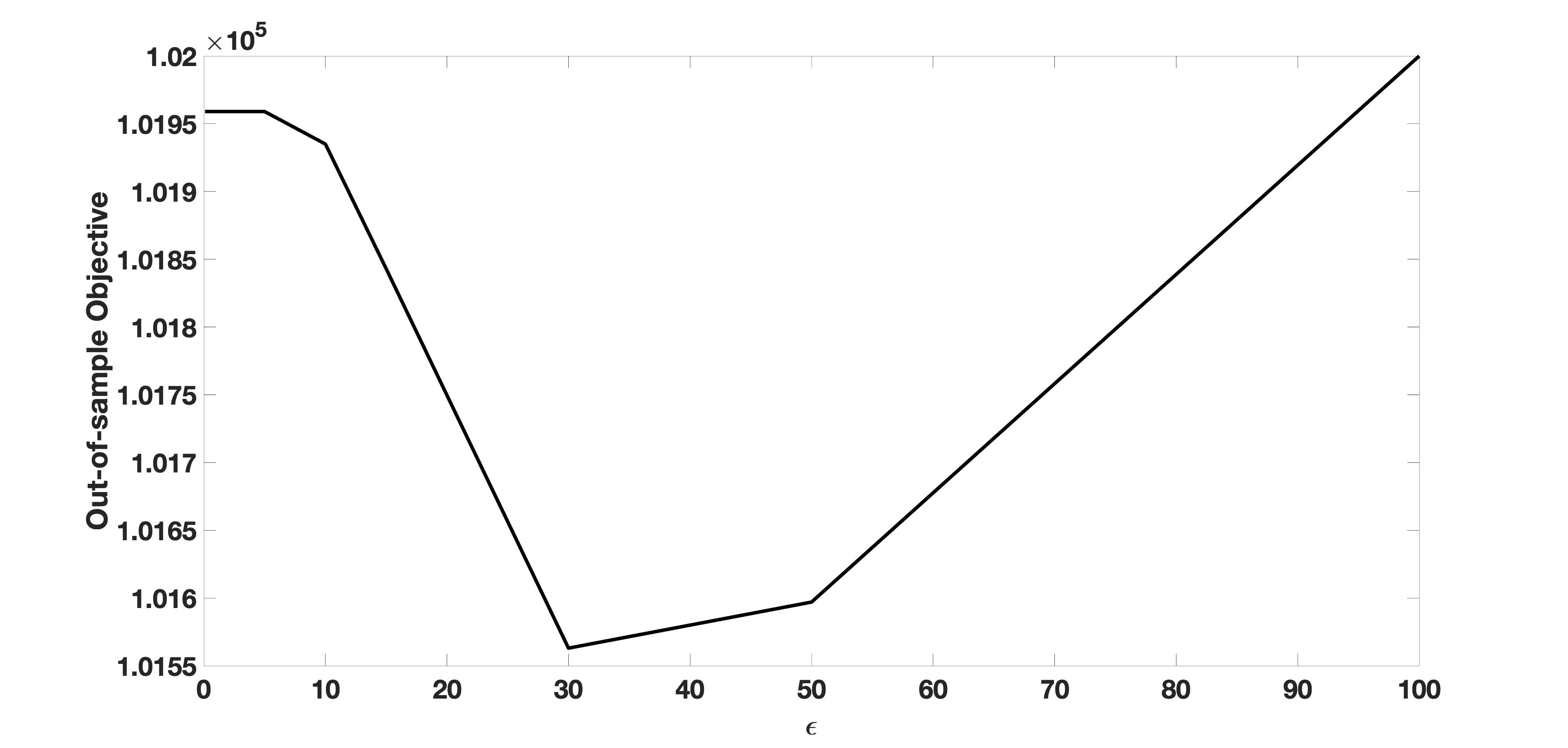}
        \caption{$N=5$}
    \end{subfigure}%
  \begin{subfigure}[b]{0.5\textwidth}
        \centering
        \includegraphics[width=\textwidth]{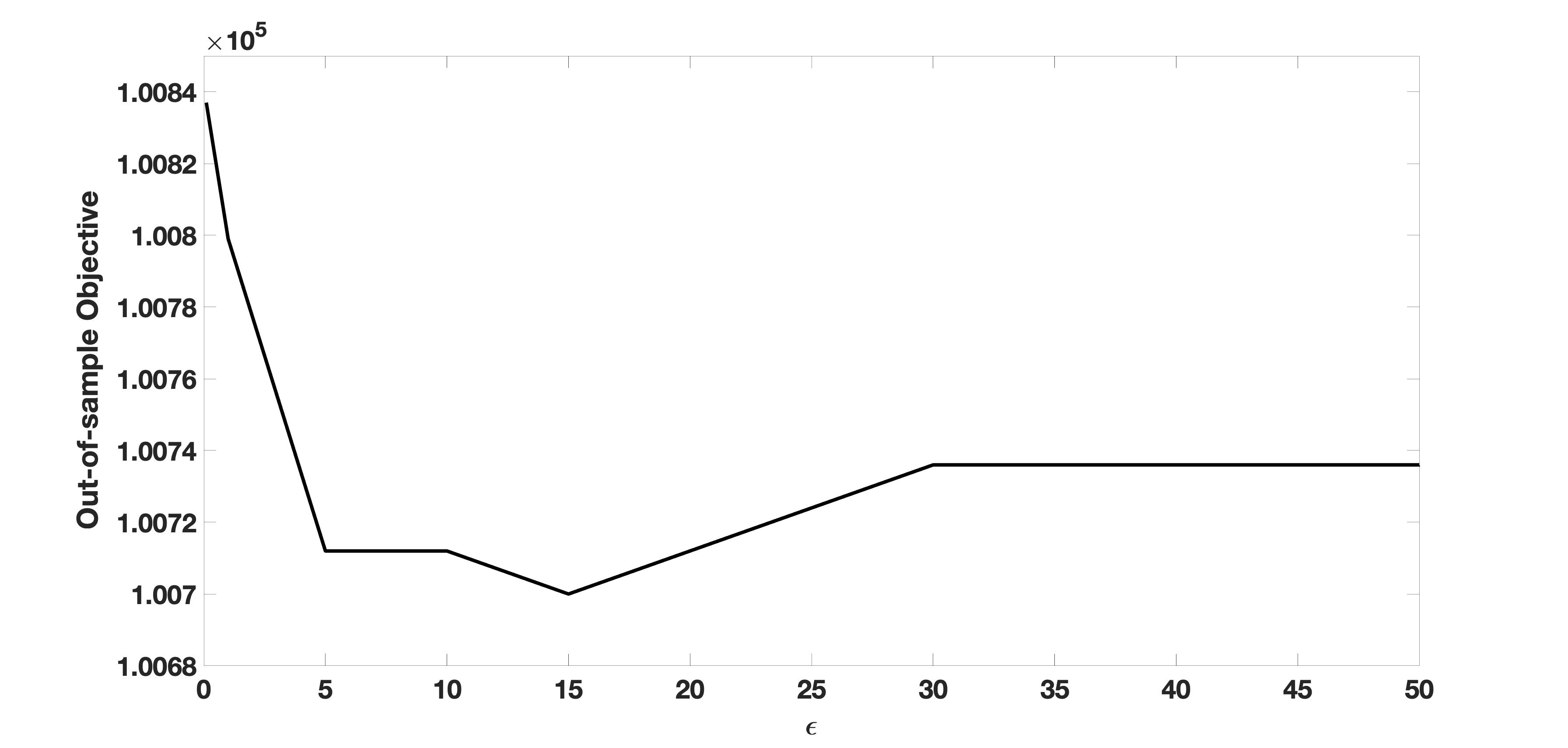}
        \caption{$N=10$}
    \end{subfigure}%
    
      \begin{subfigure}[b]{0.5\textwidth}
        \centering
        \includegraphics[width=\textwidth]{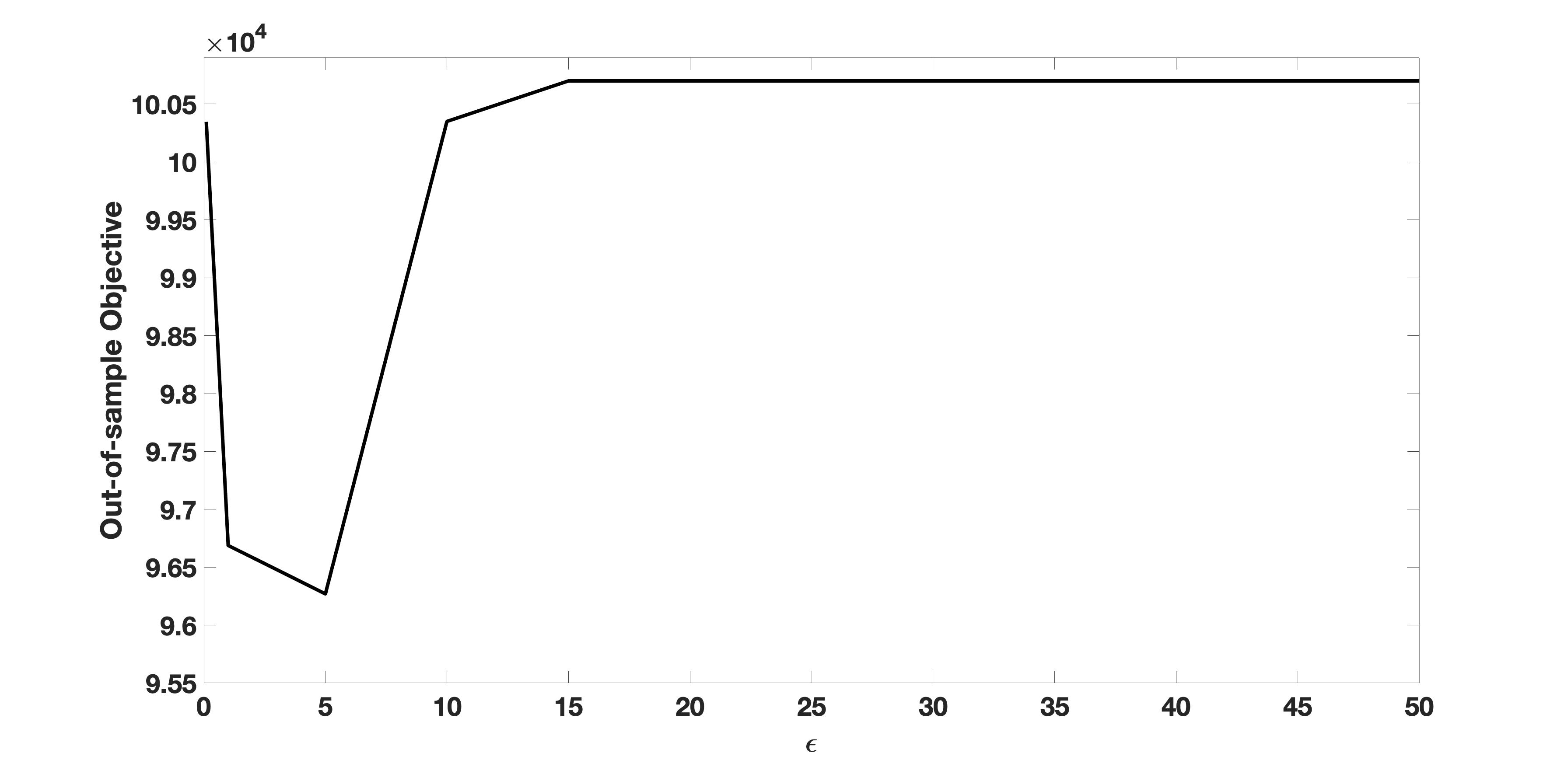}
        \caption{$N=50$}
    \end{subfigure}%
\caption{Out-of-sample performance as a function of the Wasserstein radius for an instance of $I=80$ under Cost1 }\label{Fig1_I80}
\end{figure}

\begin{figure}[h!]
 \centering
  \begin{subfigure}[b]{0.5\textwidth}
          \centering
        \includegraphics[width=\textwidth]{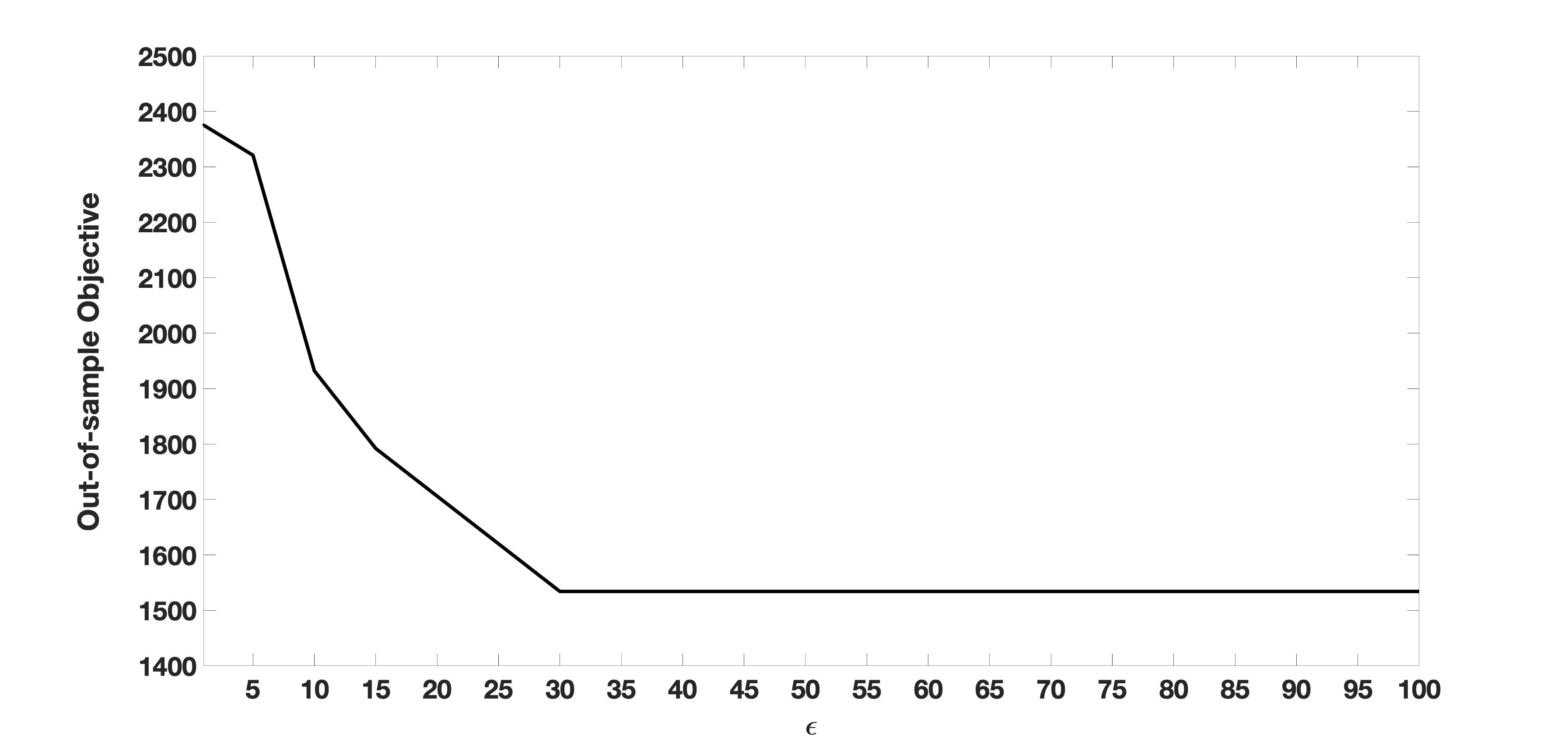}
        \caption{$N=5$}
    \end{subfigure}%
  \begin{subfigure}[b]{0.5\textwidth}
        \centering
        \includegraphics[width=\textwidth]{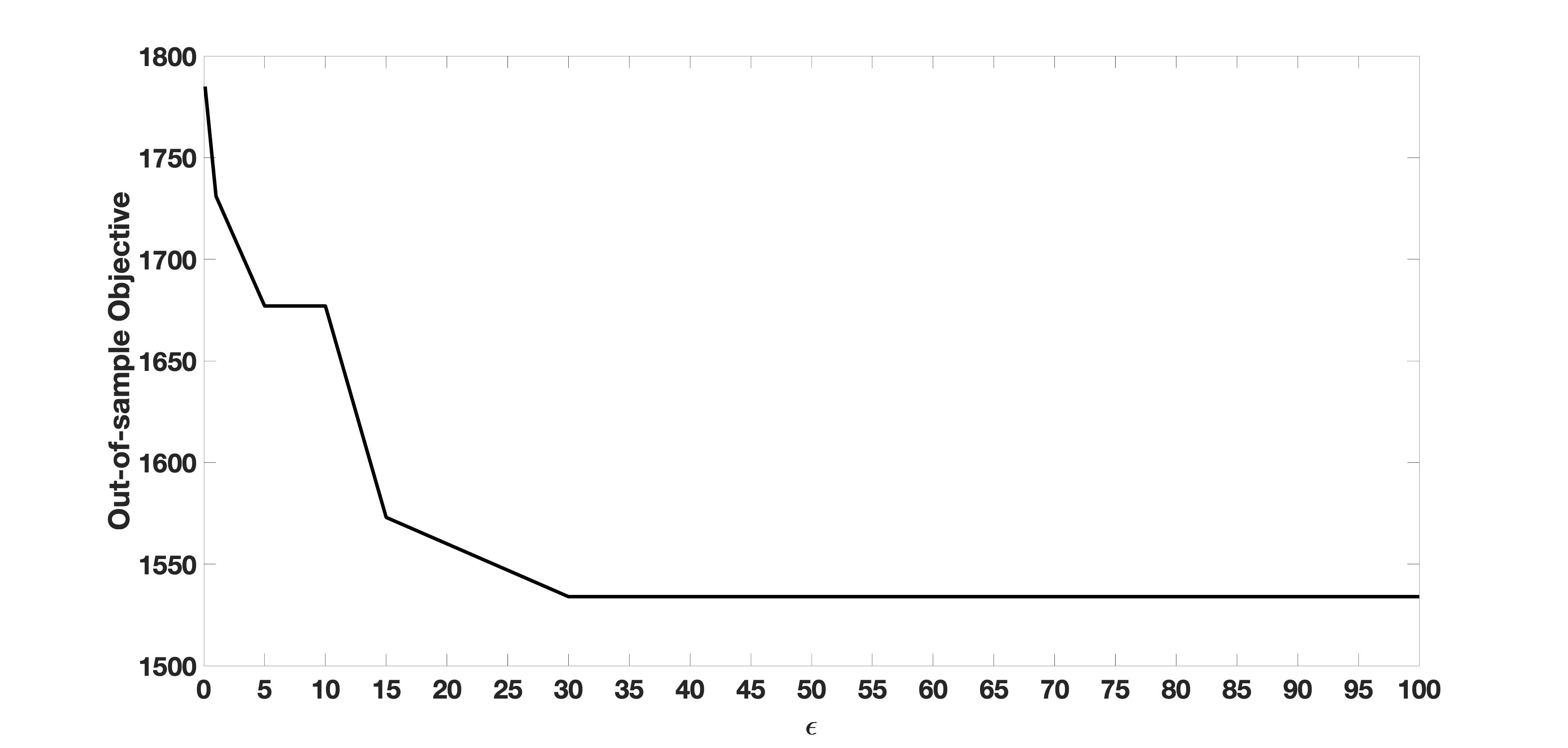}
        \caption{$N=10$}
    \end{subfigure}%
    
      \begin{subfigure}[b]{0.5\textwidth}
        \centering
        \includegraphics[width=\textwidth]{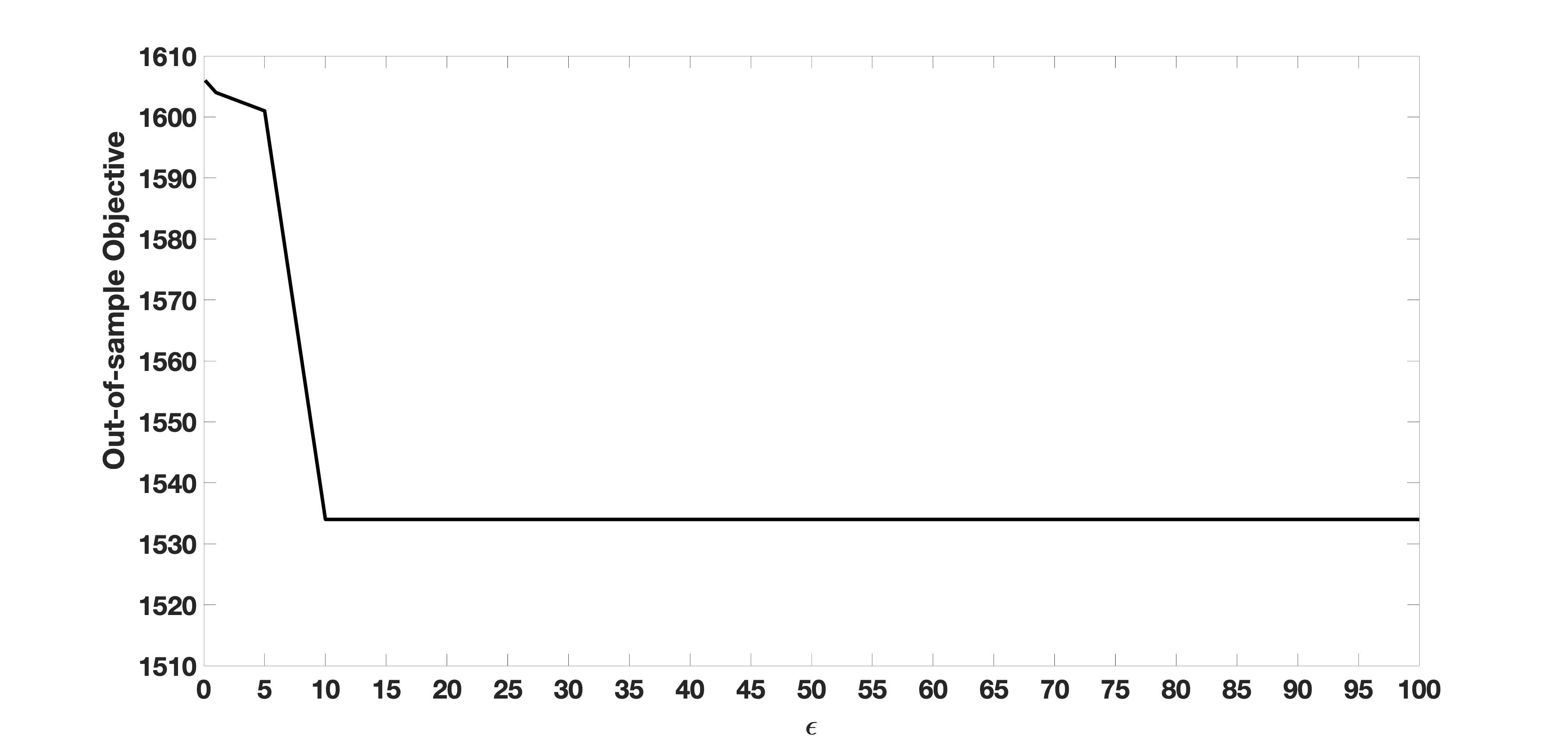}
        \caption{$N=50$}
    \end{subfigure}%
     \begin{subfigure}[b]{0.5\textwidth}
        \centering
        \includegraphics[width=\textwidth]{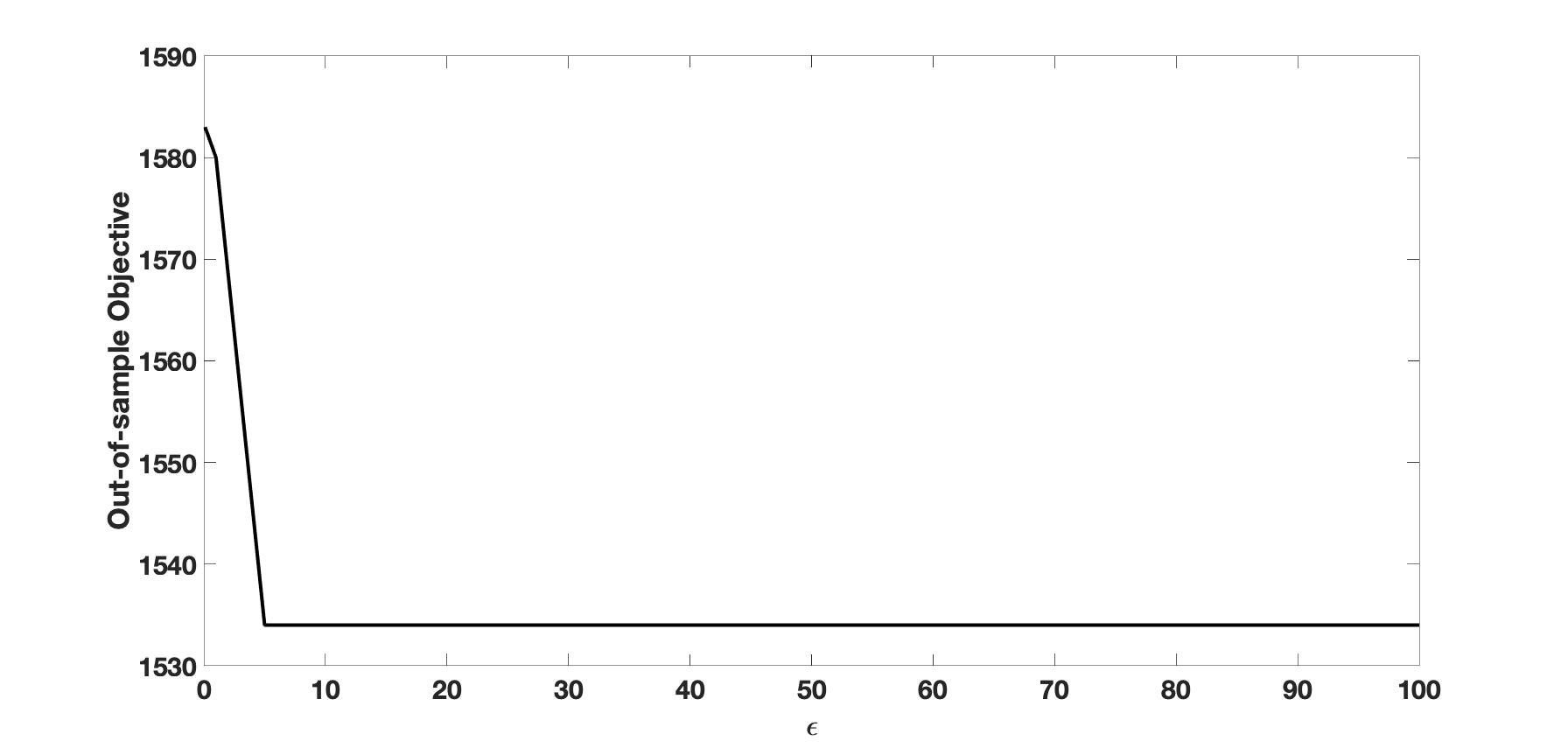}
        \caption{$N=100$}
    \end{subfigure}%
\caption{Out-of-sample performance as a function of the Wasserstein radius for an instance of $I=80$ under Cost2 }\label{Fig1_I80_cost3}
\end{figure}

\end{document}